\documentclass[10pt,francais,a4paper]{smfbook}

\usepackage{smfthm,smfenum}
\usepackage[frenchb]{babel}
\usepackage[latin1]{inputenc}
\usepackage[T1]{fontenc}
\usepackage{t1enc}
\usepackage{xy}
\usepackage{eepic,epic}
\usepackage{epsf,psfig} 
\usepackage{mathrsfs}
\usepackage{graphics,psfrag}
\usepackage{graphicx} 
\xyoption{all}

\def\Ddots{\mathinner{\mkern1mu\raise1pt\hbox{.}\mkern2mu\raise4pt\hbox{.}\mkern
2mu
  \raise7pt\vbox{\kern7pt{\hbox{.}}}\mkern1mu}}

\title{Cohomologie et K-th{\'e}orie {\'e}quivariantes des tours de Bott 
et des vari{\'e}t{\'e}s de drapeaux. Application au calcul de Schubert}

\author[Matthieu Willems ]{Matthieu Willems }

\address{U.F.R. de Math{\'e}matiques \\
 Universit{\'e} Paris 7 Denis Diderot \\ 
  2, place Jussieu, case 7012 \\
75251 Paris Cedex 05 }

\email{willems@math.jussieu.fr}

\begin{document}

\bibliographystyle{smfplain}




\maketitle 

\mainmatter

\tableofcontents

\chapter*{Introduction}

Soient $K$ un groupe de Lie compact connexe et $T$ un tore compact maximal de
$K$ d'alg{\`e}bre de Lie $\mathfrak{t}$. 
On note $X=K/T$ la vari{\'e}t{\'e} de drapeaux associ{\'e}e {\`a} ces donn{\'e}es (plus
g{\'e}n{\'e}ralement, on s'int{\'e}ressera aux vari{\'e}t{\'e}s de drapeaux des groupes de Kac-Moody).

La cohomologie {\'e}quivariante des vari{\'e}t{\'e}s de drapeaux a {\'e}t{\'e} d{\'e}crite par Alberto
Arabia dans \cite{aa} et Kostant et Kumar dans \cite{kkc}. Dans \cite{bs}, Bott
et Samelson expliquent le lien entre la cohomologie ordinaire des vari{\'e}t{\'e}s de
drapeaux et celle des vari{\'e}t{\'e}s de Bott-Samelson. L'id{\'e}e initiale de ce travail
{\'e}tait de mieux comprendre ce lien dans le cadre de la cohomologie
{\'e}quivariante. 

\medskip

Soit $\Gamma$ une vari{\'e}t{\'e} de Bott-Samelson, et soit $g : \Gamma \rightarrow X$ 
l'application naturelle de $\Gamma$ dans la vari{\'e}t{\'e} de drapeaux $X$. Le tore
compact $T$ agit sur $\Gamma$ et $X$, et l'application $g$ est $T$-{\'e}quivariante.
On calcule les restrictions aux points fixes d'une base de la
cohomologie $T$-{\'e}quivariante de $\Gamma$ (on fait ce calcul plus g{\'e}n{\'e}ralement
pour toute tour de Bott), et on d{\'e}crit la structure multiplicative de
$H_T^*(\Gamma)$. En explicitant le morphisme $g^*$,
on retrouve une base de $H_T^*(X)$.

La m{\^e}me strat{\'e}gie nous permet d'expliciter une base de la $K$-th{\'e}orie
$T$-{\'e}quivariante  de $\Gamma$ et d'en d{\'e}duire des r{\'e}sultats sur
$K_T(X)$. 

Gr{\^a}ce {\`a} ces r{\'e}sultats en cohomologie {\'e}quivariante, on donne une m{\'e}thode de
calcul des constantes de structure de $H_T^*(X)$ (calcul de Schubert
{\'e}quivariant).

\bigskip

\bigskip

Dans le chapitre~\ref{notations}, on fixe les notations sur les vari{\'e}t{\'e}s de 
drapeaux d'un groupe de Kac-Moody.

\medskip

Dans le chapitre~\ref{deftb}, on rappelle les d{\'e}finitions et les r{\'e}sultats
de~\cite{tours} sur les tours de Bott et leur lien avec les vari{\'e}t{\'e}s de
Bott-Samelson. Les tours de Bott sont des vari{\'e}t{\'e}s toriques particuli{\`e}res munies
de l'action d'un tore $D$. Une vari{\'e}t{\'e} de Bott-Samelson $\Gamma$ munie de
l'action  du tore
$T$ peut {\^e}tre vue comme une tour de Bott $Y$, et  
l'action de $T$ sur $\Gamma$ s'identifie {\`a} celle d'un sous-tore de $D$ sur $Y$.
On d{\'e}finit des d{\'e}compositions cellulaires de ces vari{\'e}t{\'e}s, 
et on pr{\'e}cise le lien entre les vari{\'e}t{\'e}s de Bott-Samelson et les vari{\'e}t{\'e}s 
de Schubert $\overline{X_w}$, index{\'e}es par les {\'e}l{\'e}ments $w$ du groupe de Weyl $W$.

\medskip

Le chapitre~\ref{cohomologie} est consacr{\'e} {\`a} la cohomologie. On calcule les
restrictions aux points fixes d'une base de la cohomologie $D$-{\'e}quivariante des
tours de Bott (th{\'e}or{\`e}me~\ref{cohomologietours}), 
et on commence {\`a} d{\'e}crire la structure multiplicative de ces
alg{\`e}bres (th{\'e}or{\`e}mes~\ref{produittours} et \ref{carre}). Par restriction, 
on obtient des r{\'e}sultats similaires pour la
cohomologie $T$-{\'e}quivariante des vari{\'e}t{\'e}s de Bott-Samelson
(th{\'e}or{\`e}mes~\ref{cohomologieBS}, \ref{produitBS} et \ref{carreBS}). Gr{\^a}ce {\`a} ces
r{\'e}sultats et {\`a} la description du morphisme $g^*$, on retrouve l'expression des
restrictions aux points fixes de la base $\{\hat{\xi}^w \}_{w \in W}$ de
$H_T^*(X)$ (th{\'e}or{\`e}me~\ref{cohomologieVD}) d{\'e}montr{\'e}e par Sara Billey dans \cite{sb}. 

\medskip

Le chapitre~\ref{ktheorie} est consacr{\'e} {\`a} la K-th{\'e}orie. On construit une base 
de la K-th{\'e}orie $D$-{\'e}quivariante des
tours de Bott, et on explicite les restrictions aux points fixes des {\'e}l{\'e}ments de
cette base (proposition~\ref{propbasektheorie} et
th{\'e}or{\`e}me~\ref{restrictionkb}). Par restriction, on obtient
une base de la K-th{\'e}orie $T$-{\'e}quivariante des vari{\'e}t{\'e}s de Bott-Samelson
(proposition~\ref{propbasekbs} et th{\'e}or{\`e}me~\ref{ktheorieBS}). Ces r{\'e}sultats et
la description du morphisme $g^*$ nous permettent de calculer les restrictions aux
points fixes de la base $\{\hat{\psi}^w \}_{w \in W}$ de $K_T(X)$
d{\'e}finie par Kostant et Kumar
dans \cite{kkk} (th{\'e}or{\`e}me~\ref{ktheorieVD}). Dans le cas fini, on explicite 
la matrice de changement de base entre $\{\hat{\psi}^w \}_{w \in W}$
 et la base $\{*[\mathcal{O}_{\overline{X_w}}]\}_{w \in W}$ de $K_T(X)$, o{\`u}
 $[\mathcal{O}_{\overline{X_w}}]$ 
est d{\'e}fini {\`a}
 partir  du faisceau structural
de $\overline{X_w}$, et o{\`u} $*$ d{\'e}signe la dualit{\'e} en $K$-th{\'e}orie induite par la dualit{\'e} des
fibr{\'e}s vectoriels (th{\'e}or{\`e}me~\ref{changementdebases}). Dans la 
section~\ref{hecke}, on explique le th{\'e}or{\`e}me~\ref{ktheorieVD} {\`a} l'aide des
alg{\`e}bres de Hecke.

\medskip

Dans le chapitre~\ref{schubertequivariant}, on pr{\'e}cise la structure
multiplicative de la cohomologie $T$-{\'e}quivariante des vari{\'e}t{\'e}s de Bott-Samelson
(th{\'e}or{\`e}me~\ref{multiplicationBS}), et on en d{\'e}duit une formule pour calculer
les constantes de structure de $H_T^*(X)$ (th{\'e}or{\`e}me~\ref{formulefinale}),
i.e. les polyn{\^o}mes $p_{u,v}^w \in S(\mathfrak{t}_{\mathbb{C}}^*)$ v{\'e}rifiant : 
$$\hat{\xi}^u \hat{\xi}^v = \sum_{w \in W} p_{u,v}^w \hat{\xi}^w.$$ 

Pour calculer ces polyn{\^o}mes $p_{u,v}^w$,
 on envoie $H_T^*(X)$ dans $H_T^*(\Gamma)$ o{\`u} les
multiplications sont plus faciles.
Le th{\'e}or{\`e}me~\ref{formulefinale} g{\'e}n{\'e}ralise la
formule donn{\'e}e par Haibao Duan pour la cohomologie ordinaire 
\cite{duan}.

\bigskip

\bigskip

Je remercie Alberto Arabia de m'avoir fait comprendre la structure des vari{\'e}t{\'e}s
de Bott-Samelson et de leur cohomologie.

\medskip

En r{\'e}digeant les
textes~\cite{mw}, \cite{mw3}, et \cite{mw2} qui donnent une partie des r{\'e}sultats
des chapitres~\ref{cohomologie} et \ref{ktheorie}, j'ai eu connaissance 
de r{\'e}sultats de William Graham qui prouve le 
th{\'e}or{\`e}me~\ref{ktheorieVD} gr{\^a}ce {\`a} d'autres m{\'e}thodes dans la
pr{\'e}publication~\cite{gra}.

\chapter{Pr{\'e}liminaires et notations} \label{notations}

\section{Alg{\`e}bres de Kac-Moody}

Les d{\'e}finitions et les r{\'e}sultats qui suivent sur les alg{\`e}bres de Kac-Moody
 sont expos{\'e}s dans \cite{km} et \cite{livrekumar}. Soit
 $A=(a_{ij})_{1\leq i,j \leq r}$ une matrice de Cartan g{\'e}n{\'e}ralis{\'e}e (c'est-{\`a}-dire 
telle que $a_{ii}=2$, $-a_{ij} \in \mathbb{N}$ si $i \neq j$, et
  $a_{ij}=0$ si et seulement si $a_{ji}=0$). On choisit un triplet $(\mathfrak{h},
  \mathfrak{\pi}, \mathfrak{\pi^{\vee}})$ (unique {\`a} isomorphisme pr{\`e}s), o{\`u}
  $\mathfrak{h}$ est un $\mathbb{C}$-espace vectoriel de dimension $(2r-{ \rm
 rg}(A))$,
  $\mathfrak{\pi} = \{\alpha_{i}\}_{1 \leq i \leq r} \subset \mathfrak{h}^*$, et 
$\mathfrak{\pi^{\vee}} = \{h_i\}_{1 \leq i \leq r} \subset \mathfrak{h}$ sont
  des ensembles d'{\'e}l{\'e}ments lin{\'e}airement ind{\'e}pendants v{\'e}rifiant
  $\alpha_{j}(h_i)=a_{ij}$.  On notera aussi $h_i$ par $\alpha^{\vee}_{i}$.
 L'alg{\`e}bre de Kac-Moody
  $\mathfrak{g}=\mathfrak{g}(A)$ est l'alg{\`e}bre de Lie sur $\mathbb{C}$ engendr{\'e}e
  par $\mathfrak{h}$ et par les symboles $e_{i}$ et $f_{i}$ ($1 \leq i \leq r$)
  soumis aux relations $[\mathfrak{h},\mathfrak{h}]=0$, $[h,
  e_{i}]=\alpha_{i}(h)e_{i}$, $[h, f_{i}]=-\alpha_{i}(h)f_{i}$ pour tout $h \in
  \mathfrak{h}$ et tout $1 \leq i \leq r$, $[e_{i}, f_{j}]=\delta_{ij}h_{j}$
  pour tout $1 \leq i,j \leq r$, et : 
$$({\rm ad }e_{i})^{1-a_{ij}}(e_{j})=0=({\rm ad }f_{i})^{1-a_{ij}}(f_{j})
   \hspace{0,2 cm} \forall \hspace{0,2 cm}  1 \leq i \neq j \leq r.$$

L'alg{\`e}bre $\mathfrak{h}$ s'injecte canoniquement dans $\mathfrak{g}$. On
l'appelle la sous-alg{\`e}bre de Cartan de $\mathfrak{g}$. On a la
d{\'e}composition suivante : 
$$\mathfrak{g}=\mathfrak{h} \oplus \sum_{\alpha \in
  \Delta_{+}}(\mathfrak{g}_{\alpha} \oplus \mathfrak{g}_{-\alpha}),$$
o{\`u} pour $\lambda \in \mathfrak{h}^*$, $\mathfrak{g}_{\lambda} = \{ x \in
\mathfrak{g} \: {\rm tels \: que }\: [h, x]=\lambda(h)x, \forall h \in
\mathfrak{h} \}$, et o{\`u} on d{\'e}finit $\Delta_{+}$ par $\Delta_{+} = \{ \alpha \in
\sum_{i=1}^{r}\mathbb{N}\alpha_{i} \: {\rm tels \: que } \: \alpha \neq 0
\:{\rm et } \: \mathfrak{g}_{\alpha} \neq 0 \}$. On pose $\Delta=\Delta_{+} \cup
\Delta_{-}$ o{\`u} $\Delta_{-} = -\Delta_{+}$. On appelle $\Delta_{+}$
(respectivement $\Delta_{-}$) l'ensemble des racines positives (respectivement
n{\'e}gatives). Les racines $\{\alpha_{i}\}_{1 \leq i \leq r}$ sont appel{\'e}es les
racines simples. On d{\'e}finit une sous-alg{\`e}bre de Borel $\mathfrak{b}$ de $\mathfrak{g}$ par
$\mathfrak{b}=\mathfrak{h} \oplus \sum_{\alpha \in \Delta_{+}}
\mathfrak{g}_{\alpha}$.

\medskip

Au couple $(\mathfrak{g}, \mathfrak{h})$, on associe le groupe de Weyl 
$W\subset { \rm Aut}(\mathfrak{h}^*)$, engendr{\'e} par les r{\'e}flexions simples 
$\{s_{i}\}_{1  \leq i \leq r}$ d{\'e}finies par :  
$$ \forall \lambda \in \mathfrak{h}^*, \,\,
 s_{i}(\lambda)=\lambda-\lambda(h_{i})\alpha_{i} .$$

Si on note $S$ l'ensemble des r{\'e}flexions simples, le couple $(W,S)$ est un
syst{\`e}me de Coxeter. On a donc une notion d'ordre de Bruhat qu'on notera 
$u \leq v$ et une notion de longueur
qu'on notera $l(w)$. On notera $1$ l'{\'e}l{\'e}ment neutre de $W$ et dans le cas fini
(i.e. $W$ fini $\Leftrightarrow
\mathfrak{g}$ de dimension finie), 
on note $w_{0}$ le plus grand {\'e}l{\'e}ment de $W$.

On obtient une repr{\'e}sentation de $W$ dans $\mathfrak{h}$ par dualit{\'e}. Plus
pr{\'e}cis{\'e}ment, pour tout $1 \leq i \leq r$, on a : 
$$\forall h \in \mathfrak{h}, \, s_i(h)=h-\alpha_i(h)h_i.$$

Le groupe de Weyl pr{\'e}serve $\Delta$. On pose $R=W\pi$, c'est l'ensemble des
racines r{\'e}elles. On pose $R^{+}=R\cap
\Delta_{+}$, et pour $\beta = w\alpha_{i} \in R^{+}$, on pose
$s_{\beta}=ws_{i}w^{-1} \in W$ (qui est ind{\'e}pendant du choix du couple $(w,
\alpha_{i})$ v{\'e}rifiant $\beta = w\alpha_{i}$) et 
$\beta^{\vee} = wh_{i} \in \mathfrak{h}$.

Pour $(v,w) \in W^2$, on note $v \rightarrow w$ si les deux 
conditions suivantes sont v{\'e}rifi{\'e}es : 
\begin{enumerate}
\item[$(i)$] \mbox{il existe} $\beta \in R^+$  \mbox{tel
que} \,  $w=vs_{\beta}$,
 \item[$(ii)$] $l(w)=l(v)+1$,
\end{enumerate}
et  pour un tel couple, on pose
$\beta(v,w)=\beta$ (cette racine est bien uniquement d{\'e}termin{\'e}e par la donn{\'e}e de
$v$ et de $w$) et $\beta^{\vee}(v,w)= \beta^{\vee}$.

Pour tout {\'e}l{\'e}ment $w$ de
$W$, on d{\'e}finit l'ensemble $\Delta(w)$ des inversions de $w$ par
 $\Delta(w)=\Delta_{+} \cap w^{-1}\Delta_{-}$.

\medskip

On fixe un r{\'e}seau $\mathfrak{h}_{\mathbb{Z}} \subset \mathfrak{h}$ tel que :

\smallskip

\begin{enumerate}

\item[$(i)$] $\mathfrak{h}_{\mathbb{Z}} \otimes_{\mathbb{Z}}\mathbb{C}=\mathfrak{h}$,

\item[$(ii)$] $h_{i} \in \mathfrak{h}_{\mathbb{Z}}$ pour tout $1 \leq i \leq r$,

\item[$(iii)$] $\mathfrak{h}_{\mathbb{Z}}/ \sum_{i=1}^{r}\mathbb{Z}h_{i}$ est sans
torsion,

\item[$(iv)$] $\alpha_{i} \in \mathfrak{h}_{\mathbb{Z}}^* =
{ \rm Hom}(\mathfrak{h}_{\mathbb{Z}}, \mathbb{Z})$ ($\subset \mathfrak{h}^*$) pour tout 
$1 \leq i \leq r$.
 
\end{enumerate}

\medskip

On choisit des poids fondamentaux $\rho_{i} \in \mathfrak{h}_{\mathbb{Z}}^*$ 
($1 \leq i \leq r$) qui v{\'e}rifient $\rho_{i}(h_{j})=\delta_{i, j}$, pour tout 
$1 \leq i,j \leq r$. On pose $\rho=\sum_{i=1}^{r}\rho_{i}$.

\section{Groupes de Kac-Moody et vari{\'e}t{\'e}s de drapeaux}

On note $G=G(A)$ le groupe de Kac-Moody associ{\'e} {\`a} $\mathfrak{g}$ par Kac et
Peterson dans \cite{kp}.  On notera $e$ l'{\'e}l{\'e}ment neutre de $G$. 
Dans le cas fini, $G$ est un groupe de Lie semi-simple
complexe connexe et simplement connexe. On note $H \subset B \subset G$ les sous-groupes de $G$
associ{\'e}s respectivement {\`a} $\mathfrak{h}$ et $\mathfrak{b}$. Soit $K$ la forme
unitaire standard de $G$ et $T=K \cap H$ le tore maximal de $K$ associ{\'e} {\`a}
$\mathfrak{h}$. On notera $\mathfrak{t} \subset \mathfrak{h}$ 
l'alg{\`e}bre de Lie de $T$. Les racines
$\alpha_i$ et les poids fondamentaux $\rho_i$ appartiennent {\`a}  $i\mathfrak{t}^*$

 Soit $N_{G}(H)$ 
le normalisateur de $H$ dans $G$, le groupe quotient
$N_{G}(H)/H$ s'identifie {\`a} $W$. On pose $X=G/B=K/T$. C'est une vari{\'e}t{\'e} de
drapeaux g{\'e}n{\'e}ralis{\'e}e. On fait agir $H$ sur $X$ par
multiplication {\`a} gauche, ce qui induit une action de $T$ sur $X$. L'ensemble des
points fixes de $T$ dans $X$ s'identifie {\`a} $W$. Pour toute
racine simple $\alpha$, on note $G_{\alpha}$ le sous-groupe connexe de $G$
associ{\'e} {\`a} la sous-alg{\`e}bre de Lie $\mathfrak{g}_{-\alpha}\oplus\mathfrak{h}\oplus 
\mathfrak{g}_{\alpha} \subset \mathfrak{g}$ et on pose $K_{\alpha}=K \cap 
 G_{\alpha} $. Pour $w
\in W$, on d{\'e}finit $C(w)=B \cup BwB$ et pour toute racine simple $\alpha$, on
d{\'e}finit le sous-groupe $P_{\alpha}$ de $G$ par $P_{\alpha}=C(s_{\alpha})$. Le
groupe $K_{\alpha} \subset P_{\alpha} $ est un sous-groupe compact maximal 
de $ P_{\alpha} $ et $K_{\alpha}\cap B =T$.
On a la d{\'e}composition de Bruhat $G=\bigsqcup_{w\in W}BwB$ et
si on pose $X_{w}=BwB/B$, $X=\bigsqcup_{w\in W}X_{w}$. Pour tout $w \in W$,
 la cellule de Schubert $X_{w}$ est isomorphe {\`a}
$\mathbb{R}^{2l(w)}$. On obtient ainsi  une d{\'e}composition cellulaire T-invariante de $X$ o{\`u} toutes les
cellules sont de dimension paire. 

Pour tout $w \in W$, la vari{\'e}t{\'e} de Schubert $\overline{X_{w}}$ est l'adh{\'e}rence 
de la cellule $X_{w}$. 
C'est une sous-vari{\'e}t{\'e} irr{\'e}ductible et $T$-invariante de $X$ de dimension r{\'e}elle
$2l(w)$. Les vari{\'e}t{\'e}s de Schubert ne sont pas lisses en g{\'e}n{\'e}ral. Pour tout $w
\in W$, on a : 
$$\overline{X_{w}}=\bigsqcup_{w' \leq w }X_{w'}.$$

\section{Le mono{\"\i}de \underline{$W$} }

On d{\'e}finit le mono{\"\i}de $\underline{W}$ comme le mono{\"\i}de engendr{\'e} par les {\'e}l{\'e}ments
$\{\underline{s}_{i}\}_{1 \leq i \leq r}$ soumis aux relations
$\underline{s}_{i}^2=\underline{s}_{i}$ et
aux  relations de tresses  de $W$ :  
$$ \left\{ \begin{array}{cc}
     \underline{s}_{i}^2=\underline{s}_{i} & \\ 
  \underbrace{\underline{s}_{i}\, \underline{s}_{j} \cdots }_{m_{i,j}\, { \rm termes}}=
  \underbrace{\underline{s}_{j}\, \underline{s}_{i} \cdots }_{m_{i,j} \, { \rm
      termes}  }
  & { \rm si }\,
  m_{i,j}<\infty \, ,
   \end{array} \right.$$
o{\`u} $m_{i,j}$ est l'ordre de  $s_is_j$ dans $W$.

\medskip

D'apr{\`e}s l'{\'e}tude g{\'e}n{\'e}rale des alg{\`e}bres de Hecke (voir \cite{hum}), l'ensemble 
$\underline{W}$ s'identifie {\`a} l'ensemble
$W$. Pour un {\'e}l{\'e}ment $w$ de $W$, on notera $\underline{w}$ l'{\'e}l{\'e}ment correspondant
dans $\underline{W}$ d{\'e}fini par $\underline{w}=\underline{s}_{i_{1}} \cdots
\underline{s}_{i_{l}}$ si $w=s_{i_{1}} \cdots s_{i_{l}}$ est une d{\'e}composition
r{\'e}duite de $w$, et pour $\underline{v} \in \underline{W}$, on notera $v$ l'{\'e}l{\'e}ment
associ{\'e} dans $W$.

\medskip

Dans $\underline{W}$, on a les relations suivantes :

\begin{equation} \label{hecke1}  \left\{ \begin{array}{ll}
      \underline{w}\, \underline{s}_{i}=
\underline{ws_{i}}
 & {\rm si } \, \, ws_{i}>w, \\ 
\underline{w} \, \underline{s}_{i}=\underline{w}
 & {\rm si } \, \, ws_{i}<w.  \end{array} \right.
\end{equation}

\begin{equation} \label{hecke2}
 \left\{ \begin{array}{ll} \underline{s}_{i}\, \underline{w}=\underline{s_{i}w}
 & {\rm si } \, \, s_{i}w>w, \\ 
\underline{s}_{i} \, \underline{w}=\underline{w}
 & {\rm si } \, \, s_{i}w<w.  \end{array} \right.
\end{equation}

\chapter{Tours de Bott, vari{\'e}t{\'e}s de Bott-Samelson et vari{\'e}t{\'e}s de Schubert} \label{deftb}

Soit $N \geq 1$ un entier naturel. On pose $\mathcal{E}=\{0,1\}^N$. Pour $\epsilon 
=(\epsilon_1, \epsilon_2, \ldots, \epsilon_N) \in \mathcal{E}$, on note $\pi_+(\epsilon)$
 l'ensemble des entiers $i \in \{1,2,\ldots, N\}$ tels que
$\epsilon_{i}=1$ et $\pi_{-}(\epsilon)$ l'ensemble des entiers $i 
\in \{1,2,\ldots, N\} $
tels que $\epsilon_{i}=0$.
On appelle longueur de $\epsilon$, not{\'e}e
$l(\epsilon)$, le cardinal de $\pi_{+}(\epsilon)$. 
Pour $1\leq i \leq N$, on note $(i) \in \mathcal{E}$
 l'{\'e}l{\'e}ment de $\mathcal{E}$ d{\'e}fini par $(i)_{j}=\delta_{i,j}$.
 On d{\'e}finit  l'{\'e}l{\'e}ment $(\bf{1})$ de $\mathcal{E}$ par $(\bf{1}\rm)_{j}=1$ 
pour tout $j$.  On munit 
$\mathcal{E}$ d'une structure de groupe en
 identifiant $\{0,1\}$ avec $\mathbb{Z}/2\mathbb{Z}$. Pour tout entier $1 \leq n
 \leq N$, on pose $(\overline{n})=(1)+(2)+ \cdots + (n) \in \mathcal{E}$.

 On d{\'e}finit un ordre partiel sur $\mathcal{E}$ par :

$$\epsilon \leq \epsilon' \Leftrightarrow \pi_{+}(\epsilon) \subset
\pi_{+}(\epsilon').$$

\section{Tours de Bott} \label{tours}

Les d{\'e}finitions et les r{\'e}sultats des sections~\ref{211} et ~\ref{212}
sont expos{\'e}s dans \cite{tours}.

\subsection{D{\'e}finition} \label{211}

Les tours de Bott sont des vari{\'e}t{\'e}s complexes compactes et lisses construites de
la mani{\`e}re suivante : 

Soit $\mathbf{L}_2$ un fibr{\'e} en droites holomorphe sur $\mathbb{C}P^1$. On pose
$Y_2=\mathbb{P}(\mathbf{1} \oplus \mathbf{L}_2)$, o{\`u} $\mathbf{1}$ est le fibr{\'e}
en droites 
trivial au dessus de $\mathbb{C}P^1$. La vari{\'e}t{\'e} $Y_2$ est un fibr{\'e} au dessus de
$Y_1=\mathbb{C}P^1$ de fibre $\mathbb{C}P^1$; c'est une surface de
Hirzebruch. On peut it{\'e}rer ce processus {\`a} l'aide de fibr{\'e}s en droites
not{\'e}s $\mathbf{L}_2,  \mathbf{L}_3, \ldots , \mathbf{L}_N$. 
A chaque {\'e}tape, la vari{\'e}t{\'e} $Y_j$ est un
fibr{\'e} au dessus de $Y_{j-1}$ de fibre $\mathbb{C}P^1$. On obtient alors le
diagramme suivant (o{\`u} pour tout $2 \leq j \leq N$, $\mathbf{L}_j$ est un fibr{\'e}
en droites au dessus de $Y_{j-1}$) :  
$$ \begin{array}{clccccc}
 & &  & & & \mathbb{P}(\mathbf{1}\oplus \mathbf{L}_N)=&Y_N \\ 
 &  &  & &  &  \downarrow \pi_N &  \\
  &  & & & &  Y_{N-1}&  \\
 & &   & &\Ddots &  \\
  & & \mathbb{P}(\mathbf{1}\oplus  \mathbf{L}_2)=& Y_2&   & &  \\
  & & \downarrow \pi_2 & &  & &  \\
  & \mathbb{C}P^1 \,\,\, = & \!\! Y_1 & & & &  \\
  &\,\, \downarrow \pi_1  & & &  & &  \\
 \{un \, point\}= &\,\,  Y_0 & & &  & &  
\end{array}$$

\bigskip

A chaque {\'e}tape, on a deux sections particuli{\`e}res $s_j^0 : Y_{j-1} \rightarrow
Y_j$ et $s_j^{\infty} : Y_{j-1} \rightarrow Y_j$ d{\'e}finies par 
$s_j^0(x)=(x,[1,0])$ et $s_j^{\infty}(x)=(x,[0,1])$.

\medskip

Par d{\'e}finition, une tour de Bott de dimension $N$ 
est une famille $\{Y_j, \pi_j, s_j^0,s_j^{\infty}\}_{1 \leq j \leq N}$ issue
d'un diagramme du type pr{\'e}c{\'e}dent.

 On dit que deux tours de Bott 
$\{Y_j, \pi_j, s_j^0,s_j^{\infty}\}_{1 \leq j \leq N}$ et 
$\{{Y'}_j, {\pi'}_j, {s'}_j^0,{s'}_j^{\infty}\}_{1 \leq j \leq N}$
sont isomorphes s'il existe $N$ diff{\'e}omorphismes holomorphes $\{F_j : 
Y_j \rightarrow {Y'}_j \}_{1 \leq j \leq N}$ qui commutent avec les applications 
$\pi_j$, $s_j^0$, $s_j^{\infty}$ et ${\pi'}_j$, ${s'}_j^0$,
${s'}_j^{\infty}$.

\begin{exem}

$\mathbb{C}P^1 \times \cdots \times \mathbb{C}P^1$ (N fois) est une tour de Bott
de dimension $N$.

\end{exem}

\subsection{Classes d'isomorphisme des tours de Bott}  \label{212}

On se donne une liste d'entiers $C=\{c_{i,j}\} _{1 \leq
  i < j \leq N}$. On consid{\`e}re $\mathbb{R}^N$ muni de sa base canonique
  $(e_{1}, e_{2}, \ldots , e_{N})$, et on d{\'e}finit $N$ {\'e}l{\'e}ments $v_{1}, v_{2},
  \ldots , v_{N}$ de $\mathbb{R}^{N}$ par les formules suivantes :

\begin{eqnarray*}
v_{N} & = &-e_{N}, \\
v_{N-1}&=&-e_{N-1}-c_{N-1,N}e_{N},\\
  &\vdots &  \\
v_{1}&=&-e_{1}-c_{1,2}e_{2}- \cdots -c_{1,N}e_{N}.
\end{eqnarray*}

\bigskip

On d{\'e}finit l'{\'e}ventail $\Sigma_{C}$ de $\mathbb{R}^N$ comme la r{\'e}union de tous
les c{\^o}nes engendr{\'e}s par les vecteurs de sous-ensembles $\Lambda$ de $\{ e_{1}, e_{2}, 
\ldots , e_{N}, v_{1}, v_{2}, \ldots , v_{N} \}$ tels que si $e_{i} \in
\Lambda$, alors $v_{i} \notin \Lambda$. On note alors 
$Y_{C}$ la vari{\'e}t{\'e} torique associ{\'e}e {\`a} l'{\'e}ventail $\Sigma_{C}$ (voir 
\cite{cox}, ou \cite{livreaudin} chapitre $6$), c'est-{\`a}-dire le quotient de
  $(\mathbb{C}^2 \setminus {(0,0)})^{N}$ par l'action {\`a} droite de
  $(\mathbb{C}^*)^N$ 
o{\`u} le $i$-{\`e}me
 facteur de $(\mathbb{C}^*)^N$ agit sur  $(\mathbb{C}^2 \setminus {(0,0)})^{N}$ par : 
\begin{equation}\label{defaction}
\begin{array}{cc}(z_{1}, w_{1}, \ldots , z_{i-1}, w_{i-1}, z_{i}, w_{i},z_{i+1}, w_{i+1},
 \ldots,  z_{N}, w_{N} )a_{i}=\\
(z_{1}, w_{1}, \ldots ,
z_{i-1}, w_{i-1}, z_{i}a_{i}, w_{i}a_{i},  z_{i+1},
w_{i+1}a_{i}^{c_{i,i+1}},  \ldots ,z_{N}, w_{N}a_{i}^{c_{i,N}}   ).
\end{array}
\end{equation}
 
\medskip

On obtient ainsi une vari{\'e}t{\'e} complexe de dimension $N$. La vari{\'e}t{\'e} $Y_C$ est
 compacte car l'{\'e}ventail
 $\Sigma_{C}$ est complet dans $\mathbb{R}^N$ (i.e. la r{\'e}union des c{\^o}nes de
 $\Sigma_{C}$ est {\'e}gale {\`a} $\mathbb{R}^N$), et lisse car l'{\'e}ventail $\Sigma_{C}$
 est r{\'e}gulier (i.e. les c{\^o}nes de $\Sigma_{C}$ sont engendr{\'e}s par des {\'e}l{\'e}ments du
 r{\'e}seau $\mathbb{Z}^N \subset \mathbb{R}^N$ qui peuvent {\^e}tre compl{\'e}t{\'e}s en une
 base de $\mathbb{Z}^N$). 

On notera 
$[z_{1}, w_{1}, \ldots ,
z_{N}, w_{N}]$ la classe de $(z_{1}, w_{1}, \ldots , z_{N}, w_{N})$ dans $Y_C$. 

Soit $\epsilon \in \mathcal{E}$. On note $\{i_{1}<
i_{2}<\cdots <i_k \}$ les {\'e}l{\'e}ments de $\pi_+(\epsilon)$. On d{\'e}finit 
 alors une liste d'entiers
$C(\epsilon)=\{d_{i,j}(\epsilon)
\} _{1 \leq
  i < j \leq k}$ par $d_{l,m}(\epsilon)=c_{i_l,i_m}$. En particulier, pour tout entier
$1\leq n \leq N$, on pose $C_n=C((\overline{n}))=\{ c_{i,j} \}_{1 \leq i<j \leq
    n}$.  

Pour tout $2 \leq n \leq N$, 
$Y_{C_n}$ est un fibr{\'e} au dessus de $Y_{C_{n-1}}$ de fibre $\mathbb{C}P^1$. En effet,
on d{\'e}finit un fibr{\'e} en droites  $\mathbf{L}(C_{n-1},c_{1,n},c_{2,n}, \ldots,
c_{n-1,n}) $ sur 
$Y_{C_{n-1}}$ par $\mathbf{L}(C_{n-1},c_{1,n},c_{2,n}, \ldots, c_{n-1,n})=
(\mathbb{C}^2 \setminus {(0,0)})^{n-1}\times_{(\mathbb{C}^*)^{n-1}}\mathbb{C}$, o{\`u}
le $i$-{\`e}me facteur de $(\mathbb{C}^*)^{n-1}$ agit par : 
$$((z_{1}, w_{1}, \ldots , z_{n-1}, w_{n-1}),v)a_i=
((z_{1}, w_{1}, \ldots , z_{n-1}, w_{n-1})a_i,a_i^{c_{i,n}}v).$$

Ici l'action de $a_i$ sur $ (z_{1}, w_{1}, \ldots , z_{n-1}, w_{n-1}) $ est donn{\'e}e
par : 
$$ (z_{1}, w_{1}, \ldots , z_{i-1}, w_{i-1}, z_{i}, w_{i},z_{i+1}, w_{i+1},
 \ldots,  z_{n-1}, w_{n-1} )a_{i}=  $$
$$(z_{1}, w_{1}, \ldots ,
z_{i-1}, w_{i-1}, z_{i}a_{i}, w_{i}a_{i},  z_{i+1},
w_{i+1}a_{i}^{c_{i,i+1}},  \ldots ,z_{n-1}, w_{n-1}a_{i}^{c_{i,n-1}}   ).$$ 

\medskip

On v{\'e}rifie imm{\'e}diatement qu'on a bien $\mathbb{P}(\mathbf{1}\oplus
\mathbf{L}_n)=Y_{C_n}$, o{\`u} $\mathbf{1}$ d{\'e}signe le fibr{\'e} en droites trivial au dessus de 
$Y_{C_{n-1}} $, et o{\`u} on note $\mathbf{L}_n$ au lieu de
$\mathbf{L}(C_{n-1},c_{1,n},c_{2,n}, \ldots, c_{n-1,n})$. 
Si on d{\'e}finit $\pi_n : Y_{C_n} \rightarrow Y_{C_{n-1}}$ par 
$$\pi_n([z_{1}, w_{1}, \ldots  ,z_{n-1}, w_{n-1}  , z_{n}, w_{n}  ])
= [z_{1}, w_{1}, \ldots  ,z_{n-1}, w_{n-1} ], $$
la vari{\'e}t{\'e} $Y_{C}$ est alors construite {\`a} l'aide de fibrations
successives de fibres $\mathbb{C}P^1$ selon le diagramme suivant :

$$ \begin{array}{clccccc}
 & &  & & & \mathbb{P}(\mathbf{1}\oplus \mathbf{L}_N)=&Y_{C_N} = Y_C \\ 
 &  &  & &  &  \downarrow \pi_N &  \\
  &  & & & &  Y_{C_{N-1}}&  \\
 & &   & &\Ddots &  \\
  & & \mathbb{P}(\mathbf{1}\oplus  \mathbf{L}_2)=& Y_{C_2}&   & &  \\
  & & \downarrow \pi_2 & &  & &  \\
  & \mathbb{C}P^1 \,\,\, = & \!\! Y_{C_1} & & & &  \\
  &\,\, \downarrow \pi_1  & & &  & &  \\
 \{un \, point\}= &\,\,  Y_0 & & &  & &  
\end{array}$$

\bigskip

Si on d{\'e}finit $s_j^0 : Y_{j-1} \rightarrow
Y_j$ et $s_j^{\infty} : Y_{j-1} \rightarrow Y_j$ comme dans la section
pr{\'e}c{\'e}dente, alors la famille $\{Y_{C_j}, \pi_j, s_j^0, s_j^{\infty}  \}_{1 \leq
  j \leq N}$
est donc une tour de Bott de dimension $N$, et on obtient ainsi une application : 
\begin{equation} \label{isomtours}
\mathbb{Z}^{N(N-1)/2} \rightarrow \{ \mbox{ classes  d'isomorphisme 
de  tours  de   Bott 
   de  dimension } \, N \},
\end{equation}
qui {\`a} $C \in \mathbb{Z}^{N(N-1)/2}$ associe $Y_C$,
et pour toute tour de Bott $Y_C$ dans l'image de~\ref{isomtours}, une application : 
\begin{equation} \label{isomfibr{\'e}s}
\mathbb{Z}^{N} \rightarrow \{ \mbox{ classes  d'isomorphisme 
de  fibr{\'e}s  en   droites 
   holomorphes  sur }\,  Y_C \},
\end{equation}
qui {\`a} $(m_1,m_2, \ldots, m_N) \in \mathbb{Z}^{N}$ associe le fibr{\'e} 
$\mathbf{L}(C,m_1,m_2, \ldots, m_N)$ sur $Y_C$.

\smallskip

Le r{\'e}sultat suivant est alors prouv{\'e} dans \cite{tours} : 

\begin{prop}

Les applications~\ref{isomtours} et \ref{isomfibr{\'e}s} sont des bijections.

\end{prop}

\begin{rema}

Dans la proposition pr{\'e}c{\'e}dente, il s'agit de classes d'isomorphisme de tours de
Bott et non de vari{\'e}t{\'e}s complexes ayant une telle structure. Par exemple, pour
tout entier $k$, les tours de Bott $Y_{\{k \}}$ et $Y_{\{-k\}}$ repr{\'e}sentent la m{\^e}me
vari{\'e}t{\'e} complexe.

\end{rema}

Dans toute la suite on suppose donn{\'e}e une liste d'entiers $C$ et on note $Y$ 
(au lieu de $Y_{C}$) la tour de Bott associ{\'e}e {\`a} $C$. 

On fait agir $D_{\mathbb{C}}=(\mathbb{C}^*)^N$ (d'alg{\`e}bre de Lie
$\mathfrak{d}_{\mathbb{C}} \simeq \mathbb{C}^N$) 
sur $Y$ par : 

$$(e^{\lambda_{1}(d)},e^{\lambda_{2}(d)}, \ldots , e^{\lambda_{N}(d)})[z_{1}, w_{1}, 
z_{2}, w_{2}, \ldots ,z_{N},
w_{N}]=$$
$$[z_{1}, e^{-\lambda_{1}(d)} w_{1},z_{2}, e^{-\lambda_{2}(d)} w_{2},
 \ldots ,z_{N}, e^{-\lambda_{N}(d)} w_{N}],$$
o{\`u} $\lambda_{i} \in \mathfrak{d}_{\mathbb{C}}^*$ est d{\'e}finie par 
$\lambda_{i}((d_{1},d_{2}, \ldots , d_{N}))=d_{i}$. 

\medskip

Soit $S^3=\{(z,w) \in \mathbb{C}^2 : |z|^2+|w|^2=1 \}$. La vari{\'e}t{\'e} $Y$
s'identifie au quotient de $(S^3)^N$ par $(S^1)^N$ o{\`u} l'action de $(S^1)^N$ sur 
$(S^3)^N$  est
donn{\'e}e par la formule~\ref{defaction}.


 L'action de $D_{\mathbb{C}}$ sur $Y$ 
induit une action de $D=(S^1)^N$ (d'alg{\`e}bre de Lie $\mathfrak{d} \simeq
i\mathbb{R}^N \subset\mathfrak{d}_{\mathbb{C}} $) sur $Y$. 
Les $\lambda_i$ sont dans $i\mathfrak{d}^*$.

\subsection{D{\'e}composition cellulaire} On d{\'e}finit une d{\'e}composition cellulaire de
$Y$ index{\'e}e par $\mathcal{E}$ de la mani{\`e}re suivante : 

 Pour $\epsilon \in \mathcal{E}$, on note $Y_{\epsilon}
\subset Y$ l'ensemble des classes $ [z_{1}, w_{1}, \ldots ,z_{N}, w_{N}] $
 qui v{\'e}rifient pour tout entier $i$ compris entre
$1$ et $N$ : 

$$\left\{ \begin{array}{ll} w_{i}=0
 & { \rm si} \hspace{0,15 cm} \epsilon_{i} =0, \\ 
 w_{i}\neq 0
 & { \rm si }\hspace{0,15 cm} \epsilon_{i} =1.
\end{array}\right.$$

\smallskip

On v{\'e}rifie imm{\'e}diatement que cette d{\'e}finition est bien compatible avec
l'action~\ref{defaction} de $(\mathbb{C}^*)^{N}$ sur
$(\mathbb{C}^2\setminus{(0,0)})^N$. 

\medskip

 Pour $\epsilon \in \mathcal{E}$ et $1 \leq k<l \leq N$, on pose : 

$$\displaystyle{ c_{k,l} (\epsilon) = \!\!\!\!\!\!\!\!\!\!\!\!\!\!\!\!
\sum_{ {\tiny \begin{array}{cc} {i_{0}=k<i_{1}< \cdots < i_{m}=l}  \\ 
 m>0, i_{j} \in \pi_{+}(\epsilon) \end{array}}}}\!\!\!\!\!\!\!\!
\!\!\!\!\!\!\!\!(-1)^m c_{i_{0},i_{1}}c_{i_{1},i_{2}}
\cdots c_{i_{m-1},i_{m}}.$$

\medskip

 Pour $\epsilon \in \mathcal{E}$ et $i \in \{ 1, 2, \ldots , N \}$, on d{\'e}finit
 $\lambda_{i}(\epsilon) \in \mathfrak{d}_{\mathbb{C}}^* $ par:

$$\lambda_{i}(\epsilon)=(-1)^{\epsilon_{i}+1} \big( \lambda_{i}+\sum_{j<i , j \in
  \pi_{+}(\epsilon)}\!\!\!\!\!\!c_{j,i}(\epsilon)\lambda_{j} \big) ,$$
o{\`u}, par convention, $\sum_{\emptyset}=0$.

\medskip

\medskip

\medskip

On d{\'e}montre alors facilement la proposition suivante : 
\begin{prop}  \label{decompositiontours}
    
\indent

\begin{enumerate}

\item[$(i)$] Pour tout $\epsilon \in \mathcal{E}$, $Y_{\epsilon}$ est un
espace affine complexe de dimension $l(\epsilon)$ stable sous l'action 
lin{\'e}aire du tore $D_{\mathbb{C}}$. 

\item[$(ii)$] Pour tout $\epsilon \in \mathcal{E}$, $\overline{Y_{\epsilon}} = 
\coprod_{\epsilon' \leq \epsilon} Y_{\epsilon'}$.

\item[$(iii)$] $Y = \coprod_{\epsilon \in
  \mathcal{E}} Y_{\epsilon}$.

\item[$(iv)$] Pour tout $\epsilon  \in \mathcal{E}$, la sous-vari{\'e}t{\'e}
$\overline{Y_{\epsilon}}$ s'identifie {\`a} la vari{\'e}t{\'e} $Y_{C(\epsilon) }$ 
 et est donc une sous-vari{\'e}t{\'e} irr{\'e}ductible lisse de $Y$.

\end{enumerate}

\end{prop}

De plus, nous allons avoir besoin du lemme suivant dont la d{\'e}monstration est imm{\'e}diate.

\begin{lemm} \label{pointsfixestours}

\indent

\begin{enumerate}
    
\item[$(i)$] L'ensemble $Y^{D}$ des points fixes de $Y$ sous l'action de
$D$ est constitu{\'e} des $2^N$ points :

$$[z_{1}, w_{1},  z_{2}, w_{2}, \ldots, z_{N}, w_{N}], 
\hspace{0,15 cm} o\grave{u} \hspace{0,15 cm}
(z_{i},w_{i}) \in \{ (1,0), (0,1) \}.$$

On identifiera donc $Y^{D}$ avec $\mathcal{E}$ en identifiant $(1,0)$
avec $0$ et $(0,1)$ avec $1$. Le point fixe $\epsilon \in Y^{D}$ est
l'unique point fixe de $Y_{\epsilon}$.

\item[$(ii)$] Soit $(\epsilon,\epsilon') \in \mathcal{E}^2$, alors : 

$$ \epsilon \in \overline{Y_{\epsilon'}} \Leftrightarrow \epsilon \leq \epsilon',$$
et dans ce cas si on note $T_{\epsilon'}^{\epsilon}$ l'espace tangent {\`a}
$\overline{Y_{\epsilon'}}$ en $\epsilon$, les poids de la repr{\'e}sentation de
$D$ dans $T_{\epsilon'}^{\epsilon}$ induite par l'action de
$D$ sur
$Y$ sont les $\{ \lambda_{i}(\epsilon) \}_{i \in \pi_{+}(\epsilon') }$. 

\end{enumerate}

\end{lemm}

\section{Vari{\'e}t{\'e}s de Bott-Samelson} \label{BS}

On utilise les notations du chapitre~\ref{notations}.

\subsection{D{\'e}finition} \label{def}

Consid{\'e}rons une suite de $N$ racines
simples $\mu_{1}$, \ldots , $\mu_{N}$ non n{\'e}cessairement distinctes. 
On d{\'e}finit : 

$$\Gamma(\mu_{1}, \ldots ,\mu_{N})=P_{\mu_1} \times_{B} 
P_{\mu_2} \times_{B} \cdots \times_{B} P_{\mu_N}/B,$$
comme l'espace des orbites de $B^N$ dans $P_{\mu_1} \times P_{\mu_2} 
\times \cdots \times P_{\mu_N}$,
sous l'action {\`a} droite de $B^N$ d{\'e}finie par : 
$$(g_{1}, g_{2}, \ldots , g_{N})
(b_{1}, b_{2}, \ldots , b_{N}) = 
(g_{1}b_{1},b_{1}^{-1} g_{2}b_{2}, \ldots ,b_{N-1}^{-1} g_{N}b_{N}),\,
 b_{i} \in B, \, g_{i} \in P_{\mu_i}.$$

\medskip

On obtient ainsi  une vari{\'e}t{\'e}
projective irr{\'e}ductible et lisse. On notera 
$[g_{1}, g_{2}, \ldots , g_{N}]$ la classe de $(g_{1}, g_{2},
\ldots , g_{N})$ dans $\Gamma(\mu_{1}, \ldots ,\mu_{N})$. On note $g_{\mu_{i}}
\in P_{\mu_i}$ un repr{\'e}sentant quelconque de la r{\'e}flexion de 
$N_{P_{\mu_i}}(H)/H \simeq \mathbb{Z}/2\mathbb{Z}$.

\medskip

On d{\'e}finit de m{\^e}me $\Gamma^K( 
 \mu_{1}, \ldots ,\mu_{N} ) $ comme l'espace des orbites de $T^N$ dans 
$K_{\mu_1} \times K_{\mu_2} \times \cdots \times K_{\mu_N}$,
sous l'action {\`a} droite de $T^N$ d{\'e}finie par : 
$$(k_{1}, k_{2}, \ldots , k_{N})
(t_{1}, t_{2}, \ldots , t_{N}) = 
(k_{1}t_{1},t_{1}^{-1} k_{2}t_{2}, \ldots ,t_{N-1}^{-1} k_{N}t_{N}),\, t_{i} \in
T, \, k_{i} \in K_{\mu_i}.$$

\medskip

Les inclusions $K_{\mu_i} \subset P_{\mu_i}$ et $T \subset B$ induisent une application : 
$$i : 
\Gamma^K( \mu_{1}, \ldots ,\mu_{N} ) \rightarrow \Gamma(  \mu_{1}, \ldots
,\mu_{N} ). $$

Comme pour tout $i$ les inclusions $K_{\mu_i} \subset P_{\mu_i}$ et $T \subset
B$ induisent un isomorphisme de vari{\'e}t{\'e}s $C^{\infty}$ : $K_{\mu_i}/T \cong 
 P_{\mu_i}/B$, l'application $i$
est un isomorphisme de vari{\'e}t{\'e}s $C^{\infty}$. 

  On notera la vari{\'e}t{\'e} $\Gamma(\mu_{1}, \ldots ,\mu_{N} ) $ par $\Gamma$ 
lorsqu'il n'y aura pas d'ambigu{\"\i}t{\'e}.

\medskip

On d{\'e}finit une action {\`a} gauche de $B$ sur $\Gamma$ par : 
$$b[g_{1}, g_{2}, \ldots ,g_{N}]=[bg_{1}, g_{2}, \ldots ,g_{N}],\hspace{0,1 cm} b \in B,
\hspace{0,1 cm}  g_{i} \in P_{\mu_i}.$$

Par restriction, on obtient ainsi une action de $H$ et de $T$.

\subsection{D{\'e}composition cellulaire}

 Pour $\epsilon \in \mathcal{E}$, on note $\Gamma_{\epsilon}
\subset \Gamma$ l'ensemble des classes $[g_{1},
g_{2}, \ldots , g_{N}]$ qui v{\'e}rifient pour tout entier $i$ compris entre
$1$ et $N$ : 

$$\left\{ \begin{array}{ll} g_{i} \in B
 & { \rm si} \hspace{0,15 cm} \epsilon_{i} =0, \\ 
 g_{i} \notin B
 & { \rm si }\hspace{0,15 cm} \epsilon_{i} =1.
\end{array}\right.$$

On v{\'e}rifie imm{\'e}diatement que cette d{\'e}finition est bien compatible avec
l'action de $B^N$.

\medskip
 
 Pour $\epsilon \in \mathcal{E}$ et $1 \leq i\leq N$, on d{\'e}finit : 
$$\displaystyle{v_{i}(\epsilon) =\prod_{ {\tiny \begin{array}{cc}  1\leq k \leq i, \\ 
 k \in \pi_{ +}(\epsilon) \end{array}} 
  }\!\!\!\!\!\!s_{\mu_{k}}},$$
o{\`u}, par convention, $\prod_{\emptyset}=1$. On pose $v(\epsilon)=v_{N}(\epsilon)$. 

\smallskip

Pour $i \leq j$, on
d{\'e}finit {\'e}galement : 
$$\displaystyle{v_{i}^j(\epsilon) =\prod_{ {\tiny \begin{array}{cc}  i\leq k \leq j, \\ 
 k \in \pi_{ +}(\epsilon) \end{array}}}
\!\!\!\!\!\!  s_{\mu_{k}}}.$$ 

Ce sont des {\'e}l{\'e}ments de $W$.

\medskip

De plus, on pose 
$\alpha_{i}(\epsilon)=v_{i}(\epsilon)\mu_{i}$, c'est une racine.

\medskip

On d{\'e}finit de m{\^e}me : $$\displaystyle{\underline{v}(\epsilon) = 
\prod_{ {\tiny \begin{array}{cc}  1\leq k \leq N, \\ 
 k \in \pi_{ +}(\epsilon) \end{array}} } 
 \!\!\!\!\!\!  \underline{s_{\mu_{k}}} \in \underline{W}}.$$ 

\bigskip
\medskip

On d{\'e}montre alors facilement la proposition suivante : 

\begin{prop} \label{structureBS}
    
\indent

\begin{enumerate}

\item[$(i)$] Pour tout $\epsilon \in \mathcal{E}$, $\Gamma_{\epsilon}$ est un
espace affine complexe de dimension $l(\epsilon)$ stable sous l'action de $B$,
et cette action induit une action lin{\'e}aire du tore $H$ sur $\Gamma_{\epsilon}$. 

\item[$(ii)$] Pour tout $\epsilon \in \mathcal{E}$, $\overline{\Gamma_{\epsilon}} = 
\coprod_{\epsilon' \leq \epsilon} \Gamma_{\epsilon'}$.

\item[$(iii)$] $\Gamma = \coprod_{\epsilon \in
  \mathcal{E}} \Gamma_{\epsilon}$.

\item[$(iv)$] Pour tout $\epsilon  \in \mathcal{E}$, $\overline{\Gamma_{\epsilon}}$ s'identifie
{\`a} la vari{\'e}t{\'e} $\Gamma(\mu_{i}, i \in \pi_{+}(\epsilon))$ et est donc une
sous-vari{\'e}t{\'e} irr{\'e}ductible lisse de $\Gamma$.

\end{enumerate}   

\end{prop}

\begin{rema}

Dans \cite{gauss}, St{\'e}phane Gaussent d{\'e}finit d'autres d{\'e}compositions cellulaires
des vari{\'e}t{\'e}s de Bott-Samelson.

\end{rema}

Soit $\Gamma^T$ l'ensemble des points fixes de $T$ dans $\Gamma$, on peut
identifier $\Gamma^T$ avec $\mathcal{E}$ gr{\^a}ce au lemme suivant :

\begin{lemm} \label{pointsfixesBS}

\indent 
    
\begin{enumerate}

\item[$(i)$] L'ensemble $\Gamma^T$ est
constitu{\'e} des $2^N$ points :

$$[g_{1}, g_{2}, \ldots, g_{N}], \,\, \mbox{o{\`u}} \,\,\,
g_{i} \in \{ e, g_{\mu_{i}} \}.$$

On identifiera donc $\Gamma^T$ avec $\mathcal{E}$ en identifiant $e$
avec $0$ et $g_{\mu_{i}}$ avec $1$.

\item[$(ii)$] Pour tout $\epsilon \in \mathcal{E}$, le point fixe $\epsilon$ est l'unique
point fixe de $\Gamma_{\epsilon}$.

\end{enumerate}

\end{lemm}

\subsection{Une autre structure complexe}

On va d{\'e}finir une structure complexe sur  $  \Gamma^K( 
 \mu_{1}, \ldots ,\mu_{N} )$ diff{\'e}rente de celle d{\'e}duite de l'isomorphisme de
 vari{\'e}t{\'e}s $C^{\infty}$ entre $  \Gamma^K(\mu_{1}, \ldots ,\mu_{N} )$  et $\Gamma( 
 \mu_{1}, \ldots ,\mu_{N} )$. Cette structure permettra d'identifier les
 vari{\'e}t{\'e}s de Bott-Samelson {\`a} des tours de Bott en tant que vari{\'e}t{\'e}s complexes.

La projection $\mathfrak{b} \rightarrow \mathfrak{h}$ induit un homomorphisme de
groupe de Lie $\Theta :  B \rightarrow H$. On d{\'e}finit alors $\Gamma^{\infty}(
\mu_{1}, \ldots ,\mu_{N} )$ comme l'espace des orbites de $B^N$ dans $P_{\mu_1} 
\times P_{\mu_2} \times \cdots \times P_{\mu_N}$,
sous l'action {\`a} droite de $B^N$ d{\'e}finie par : 
\begin{equation} \label{infty}
\begin{array}{cc}(g_{1}, g_{2}, \ldots , g_{N})
(b_{1}, b_{2}, \ldots , b_{N}) = \\
(g_{1}b_{1},
\Theta(b_{1})^{-1} g_{2}b_{2}, \ldots ,\Theta(b_{N-1})^{-1} g_{N}b_{N}),\hspace{0,1
  cm} b_{i} \in
B, \hspace{0,1 cm} g_{i} \in P_{\mu_i}. \end{array} \end{equation}

\medskip

Comme pr{\'e}c{\'e}demment, la vari{\'e}t{\'e} $\Gamma^{\infty}(\mu_{1}, \ldots ,\mu_{N} )$ est
isomorphe (en tant que vari{\'e}t{\'e} $C^{\infty}$) {\`a} $\Gamma^{K}(\mu_{1}, \ldots
,\mu_{N} )$ et donc aussi {\`a} $\Gamma(\mu_{1}, \ldots ,\mu_{N} )$. On fait agir
$H$ sur  $\Gamma^{\infty}(\mu_{1}, \ldots ,\mu_{N} )$ de la m{\^e}me mani{\`e}re que sur 
$\Gamma(\mu_{1}, \ldots ,\mu_{N} )$

\subsection{Les vari{\'e}t{\'e}s de Bott-Samelson sont des tours de Bott} \label{toursbott}

Pour identifier les vari{\'e}t{\'e}s de Bott-Samelson {\`a} des tours de Bott, on proc{\`e}de
par r{\'e}currence sur $N \geq 1$.

 Pour $N=1$, on va expliciter un isomorphisme entre $\Gamma_{\mu}$
et $X=\mathbb{C}P^1$. On note $k \in \{1,2, \ldots, r\}$
l'indice tel que $\mu=\alpha_k$. On d{\'e}finit alors
un homomorphisme $\mathfrak{sl}(2,\mathbb{C}) \rightarrow \mathfrak{p}_{\mu}$ en
envoyant $\left(\begin{array}{cc}  0 & 1 \\
                          0 & 0 \end{array}\right)$ sur $e_k$, 
$\left(\begin{array}{cc}  0 & 0 \\
                          1 & 0 \end{array}\right)$ sur $f_k$, et 
 $\left(\begin{array}{cc}  1 & 0 \\
                          0 & -1 \end{array}\right)$ sur $h_k$. On
                      obtient ainsi un homomorphisme $SL(2,\mathbb{C})
                      \rightarrow P_{\mu}$ et un diff{\'e}omorphisme holomorphe :
$$ SL(2,\mathbb{C})/B_{ SL(2,\mathbb{C}) } \cong P_{\mu}/B,$$
 o{\`u} $B_{
  SL(2,\mathbb{C})}$ d{\'e}signe le sous-groupe de $ SL(2,\mathbb{C}) $
des matrices triangulaires sup{\'e}rieures.

 De plus,
  l'application $\left(\begin{array}{cc}  a & b \\
                          c & d \end{array}\right) \mapsto [a,c]$ d{\'e}finit un
                      diff{\'e}omorphisme holomorphe : 
$$SL(2,\mathbb{C})/B_{ SL(2,\mathbb{C}) }
                      \cong \mathbb{C}P^1 . $$

On obtient ainsi bien un
                      diff{\'e}omorphisme $\phi^{\infty}_{\mu} : \Gamma_{\mu}
 \cong \mathbb{C}P^1$.

\medskip

On montre maintenant que $\Gamma^{\infty}(\mu_1, \mu_2, \ldots, \mu_N)$ peut
s'{\'e}crire $\mathbb{P}(\mathbf{1} \oplus \mathbf{L}_N)$ o{\`u} $\mathbf{L}_N$ est un
fibr{\'e} en droites sur  $\Gamma^{\infty}(\mu_1, \mu_2, \ldots, \mu_{N-1})$.

On note $\mathbb{C}_{\mu_N}$ l'espace vectoriel $\mathbb{C}$ muni de l'action
de $B$ triviale sur la partie unipotente de $B$ et d{\'e}finie pour $h \in H$ 
par la multiplication par $e^{\mu_N}(h)$. On munit $\mathbb{C}_{\mu_N}$ d'une
action de $B^{N-1}$ en faisant agir seulement la derni{\`e}re composante de
$B^{N-1}$.
On pose alors $\mathbf{L}_N= 
(P_{\mu_1}\times P_{\mu_2} \times
 \cdots \times P_{\mu_{N-1}})\times_{B^{N-1}} \mathbb{C}_{\mu_N}$, o{\`u}
l'action de $B^{N-1}$  sur $P_{\mu_1}\times P_{\mu_2}
\times \cdots \times P_{\mu_{N-1}}$ est donn{\'e}e par l'{\'e}quation~\ref{infty} (en
rempla{\c c}ant $N$ par $N-1$).

On a alors : 
$$\mathbb{P}(\mathbf{1}\oplus \mathbf{L}_N)=(P_{\mu_1}\times P_{\mu_2}
\times \cdots \times P_{\mu_{N-1}} )
\times_{B^{N-1}}\mathbb{P}(\mathbb{C}\oplus
\mathbb{C}_{\mu_N})$$
et : 
$$\Gamma^{\infty}(\mu_1, \mu_2, \ldots, \mu_N)=(P_{\mu_1}\times P_{\mu_2}
\times \cdots \times P_{\mu_{N-1}} )
\times_{B^{N-1}}(P_{\mu_N}/B),$$
o{\`u} l'action de $B^{N-1}$  sur $P_{\mu_1}\times P_{\mu_2}
\times \cdots \times P_{\mu_{N-1}}$  est donn{\'e}e par l'{\'e}quation~\ref{infty} (en
rempla{\c c}ant $N$ par $N-1$), et o{\`u} $(b_1,b_2,\dots, b_{N-1}) \in B^{N-1}$ agit 
sur $P_{\mu_N}/B$ par 
multiplication {\`a} gauche par $(\Theta(b_{N-1}))^{-1}$.

Il suffit maintenant de trouver un diff{\'e}omorphisme holomorphe $B^{N-1}$-{\'e}quivariant
entre $P_{\mu_N}/B$ et $\mathbb{P}(\mathbb{C}\oplus\mathbb{C}_{\mu_N})$. Comme
seule la derni{\`e}re composante de $B^{N-1}$ agit sur $P_{\mu_N}/B$ et sur 
$\mathbb{P}(\mathbb{C}\oplus\mathbb{C}_{\mu_N})$, et comme la
partie unipotente de cette derni{\`e}re composante agit trivialement sur $P_{\mu_N}/B$ et sur
$\mathbb{P}(\mathbb{C}\oplus\mathbb{C}_{\mu_N})$, il suffit de trouver un
diff{\'e}omorphisme $H$-{\'e}quivariant : $P_{\mu_N}/B \rightarrow \mathbb{P}
(\mathbb{C}\oplus\mathbb{C}_{\mu_N})$, o{\`u} $h \in H$ agit sur $P_{\mu_N}/B$ par 
multiplication {\`a} gauche par $h^{-1}$ et sur $\mathbb{P}
(\mathbb{C}\oplus\mathbb{C}_{\mu_N})$ par multiplication de la deuxi{\`e}me
composante par $e^{\mu_N}(h)$. 
On montre facilement que le diff{\'e}omorphisme $\phi^{\infty}_{\mu_N}$ convient.

\medskip

On a ainsi identifi{\'e} les vari{\'e}t{\'e}s de Bott-Samelson {\`a} des tours de Bott et le
r{\'e}sultat suivant, prouv{\'e} dans \cite{tours}, nous donne la liste d'entiers 
associ{\'e}e {\`a} cette construction : 

\begin{prop}

La construction pr{\'e}c{\'e}dente d{\'e}finit un diff{\'e}omorphisme holomorphe $\phi^{\infty}$ entre 
 $\Gamma^{\infty}
(\mu_1,\mu_2, \ldots ,\mu_N)$ et la
tour de Bott $Y_C$, avec $C=\{c_{j,k} \}_{1 \leq j < k \leq N}$, o{\`u} $c_{j,k}$ est d{\'e}fini par 
$c_{j,k}=\mu_{k}(\mu^{\vee}_j)$.

La liste $C$ est donc un ensemble de nombres de Cartan d{\'e}pendant du choix de la
suite $\mu_1,\mu_2, \ldots ,\mu_N$.

\end{prop}

Soit $\tau : \mathfrak{d}_{\mathbb{C}}^* \rightarrow \mathfrak{h}^*$ l'application d{\'e}finie par
$\tau(\lambda_{i}) = \mu_{i}$. Elle envoie $\mathfrak{d}^*$ dans $\mathfrak{t}^*$.
Soit $\mathfrak{s}^* \subset \frak{t}^*$ l'image de  $\mathfrak{d}^*$ par 
$\tau$ et soit 
$\mathfrak{s}_{\mathbb{C}}^* \subset \frak{h}^*$ 
l'image de $\mathfrak{d}_{\mathbb{C}}^*$. On a les deux suites suivantes (o{\`u} les
premi{\`e}res fl{\`e}ches sont surjectives et les deuxi{\`e}mes injectives) :  
 $$\xymatrix{ \mathfrak{d}^* \ar@{->>}[r]^{\tau }\,\,\, & 
\, \frak{s}^*\,\,
 \ar@{^{(}->}[r]&   \mathfrak{t}^* }, $$
$$\xymatrix{ \mathfrak{d}_{\mathbb{C}}^* \ar@{->>}[r]^{\tau }\,\,\, & 
\, \mathfrak{s}_{\mathbb{C}}^*\,\,
 \ar@{^{(}->}[r]&   \mathfrak{h}^* }.
$$

On en d{\'e}duit la suite suivante sur les tores complexes : 
$$\xymatrix{ H \ar@{->>}[r]^{ }\,\,\, & 
\, S_{\mathbb{C}}\,\,
 \ar@{^{(}->}[r]^{\gamma }&   D_{\mathbb{C}} }. $$

On continue {\`a} noter $\gamma$ le morphisme de $H$ dans $D_{\mathbb{C}}$ ainsi
d{\'e}fini.

On a alors la suite suivante sur les tores compacts : 
$$\xymatrix{ T \ar@{->>}[r]^{ }\,\,\, & 
\, S\,\,
 \ar@{^{(}->}[r]^{\gamma }&   D }. $$

 L'action de $h \in H$ sur $\Gamma^{\infty}$ {\'e}tant donn{\'e}e par la
formule : 
$$h[g_1,g_2, \ldots, g_N]=[hg_1h^{-1},hg_2h^{-1}, \ldots, ,hg_Nh^{-1}],$$
le tore $H$ agit sur $\Gamma^{\infty}$ via son image $S_{\mathbb{C}}$.
 De plus, d'apr{\`e}s la construction de $\phi^{\infty}$ : 
$$\forall (h,x) \in H \times \Gamma^{\infty}, \,
\phi^{\infty}(hx)= \gamma(h)\phi^{\infty}(x),  $$
et l'action de $H$ sur $\Gamma^{\infty}$ s'identifie donc {\`a} celle d'un sous-tore de 
$D_{\mathbb{C}}$ sur $Y_{C}$.

\bigskip

On note $\phi : \Gamma \rightarrow Y_C$ l'isomorphisme (de vari{\'e}t{\'e}s
$C^{\infty}$) induit par $\phi^{\infty}$. On v{\'e}rifie, gr{\^a}ce {\`a} la construction 
de $\phi^{\infty}$, que $\phi$ 
 envoie $\Gamma_{\epsilon}$ sur $Y_{\epsilon}$ et le point $\epsilon \in \Gamma^T$ sur 
le point $\epsilon \in Y^{D}$.

Le tore $T$ agit sur 
$\Gamma$ via son image $S$, et on a :
$$\forall (t,x) \in T \times \Gamma, \,
\phi(tx)= \gamma(t)\phi(x).  $$

L'action de $T$ sur $\Gamma$ s'identifie donc {\`a} celle d'un sous-tore de 
$D$ sur $Y_{C}$.

\begin{rema}

L'action holomorphe de $D_{\mathbb{C}}$ sur $\Gamma^{\infty}$ ne 
d{\'e}finit pas en g{\'e}n{\'e}ral une action holomorphe
de $D_{\mathbb{C}}$ sur $\Gamma$.

\end{rema}

\section{Vari{\'e}t{\'e}s de Bott-Samelson et vari{\'e}t{\'e}s de Schubert}

Soit $\mu_{1}, \ldots , \mu_{N}$ une suite quelconque de $N$ racines
simples. On d{\'e}finit une
application $g_{\mu_{1}, \ldots, \mu_{N}}$ de $\Gamma(\mu_{1}, \ldots, \mu_{N})$ 
dans $X$ par
multiplication (c'est-{\`a}-dire $g_{\mu_{1}, \ldots, \mu_{N}}([g_{1},\ldots ,g_{N}]) = g_{1}*\cdots *g_{N}
\hspace{0.1cm}  [B]$, o{\`u} $*$ d{\'e}signe la multiplication dans le groupe $G$). 
Cette application est $T$-{\'e}quivariante. 
Nous aurons besoin du lemme suivant :

\begin{lemm} \label{image}

Soit $\mu_{1}, \ldots , \mu_{N}$ une suite quelconque de $N$ racines simples et
soit $\underline{w}=\underline{s_{\mu_{1}}} \cdots \underline{s_{\mu_{N}}}$,
alors l'image de l'application $g_{\mu_{1}, \ldots, \mu_{N}}$ est
{\'e}gale {\`a} $\overline{X_{w}}$.

\end{lemm}

\begin{proof}

On pose $ X_{\mu_{1}, \ldots, \mu_{N}}=g_{\mu_{1}, \ldots, \mu_{N}}(\Gamma
(\mu_{1}, \ldots, \mu_{N}))$. Les vari{\'e}t{\'e}s $\Gamma$ {\'e}tant compactes, $
X_{\mu_{1}, \ldots, \mu_{N}}$ est ferm{\'e}e.

Pour d{\'e}montrer ce lemme, on utilisera les relations suivantes, valables pour
tout $v \in W$ (voir
\cite{livrekumar}, chapitre $5$ sur les syst{\`e}mes de Tits) :
\begin{equation} \label{cellules}
 \left\{ \begin{array}{ll} C(s_{i})C(v)=C(s_{i}v)
 & {\rm si } \hspace{0.2 cm} s_{i}v>v, \\ 
C(s_{i})C(v)=C(v)\cup C(s_{i}v)
 & {\rm si } \hspace{0.2 cm} s_{i}v<v.  \end{array} \right.
\end{equation}

On a en particulier : 
$$ \left\{ \begin{array}{ll} C(s_{i}v) \subset C(s_{i})C(v)
 & {\rm si } \hspace{0.2 cm} s_{i}v>v, \\ 
C(v)\subset C(s_{i})C(v)
 & {\rm si } \hspace{0.2 cm} s_{i}v<v.  \end{array} \right.$$

\smallskip

En utilisant ces derni{\`e}res relations et les relations \ref{hecke2}, on constate
que  $BwB \subset   C(w) \subset 
P_{\mu_{1}}P_{\mu_{2}} \cdots
P_{\mu_{N}} $, et on
a donc $X_{w}\subset X_{\mu_{1}, \ldots, \mu_{N}}$, d'o{\`u} $\overline{X_{w}}\subset 
X_{\mu_{1}, \ldots, \mu_{N}}$.

\medskip

Pour d{\'e}montrer l'inclusion r{\'e}ciproque, on proc{\`e}de par r{\'e}currence sur $N$. Le
r{\'e}sultat {\'e}tant trivial pour $N=1$, on le suppose vrai pour toute suite de $N-1$
racines simples. On note $w'$ l'{\'e}l{\'e}ment de $W$ tel que $\underline{w'}=
\underline{s_{\mu_{2}}} \cdots \underline{s_{\mu_{N}}}$. Par hypoth{\`e}se de
r{\'e}currence $P_{\mu_{2}}P_{\mu_{3}} \cdots
P_{\mu_{N}} \subset \coprod_{v \leq w'} BvB$. Il suffit donc de montrer que
pour tout $v \leq w', C(s_{\mu_{1}})C(v) \subset \coprod_{u \leq w}BuB$.

Distinguons deux cas. 

Tout d'abord si $s_{\mu_{1}}w'>w'$, alors, d'apr{\`e}s les
  relations~\ref{hecke2}, $w=s_{\mu_{1}}w'$ et le r{\'e}sultat est alors une
  cons{\'e}quence imm{\'e}diate des relations~\ref{cellules}, en remarquant que si $v\leq
  w'$ alors $s_{\mu_{1}}v \leq w$.

Si $s_{\mu_{1}}w'<w'$, alors, d'apr{\`e}s les relations \ref{hecke2}, $w=w'$. Soit
  $v$ un {\'e}l{\'e}ment de $W$ tel que $v \leq w$.

Si $s_{\mu_{1}}v<v$, alors, d'apr{\`e}s les relations~\ref{cellules}, 
$C(s_{\mu_{1}})C(v) = C(v) \cup C(s_{\mu_{1}}v) \subset \coprod_{u \leq
  w}BuB$. 

Si $s_{\mu_{1}}v>v$, $C(s_{\mu_{1}})C(v) = C(s_{\mu_{1}}v)$ et il faut
  donc montrer que $s_{\mu_{1}}v \leq w$. Cela provient du fait que
  $w=s_{\mu_{1}}x$ avec $s_{\mu_{1}}x>x$ et $v \leq x$ car aucune
  d{\'e}composition r{\'e}duite de $v$ ne commence par $s_{\mu_{1}}$ puisque
  $s_{\mu_{1}}v>v$.

\end{proof}

De plus, le r{\'e}sultat suivant est prouv{\'e} dans \cite{hbs} et \cite{demazure} 
dans le cas fini (voir
\cite{livrekumar}, proposition $7.1.15$ pour la g{\'e}n{\'e}ralisation aux vari{\'e}t{\'e}s de
Schubert des groupes de Kac-Moody) :

\begin{prop} \label{desingularisation}

Si $w=s_{\mu_{1}} \cdots s_{\mu_{N}}$ est une d{\'e}composition r{\'e}duite de $w$,
l'application $g_{\mu_{1}, \ldots, \mu_{N}} : \Gamma(\mu_{1}, \ldots, \mu_{N})
\rightarrow \overline{X_{w}}$ est  une 
d{\'e}singularisation de $\overline{X_{w}}$. Cette application est un isomorphisme 
au dessus de la cellule de Schubert $X_w$.

\end{prop}

\medskip

\begin{rema}

M{\^e}me quand la vari{\'e}t{\'e} $\overline{X_{w}}$ est lisse, l'application  $g_{\mu_{1}, \ldots,
  \mu_{N}}$ permet de comprendre plus
  facilement la g{\'e}om{\'e}trie de $\overline{X_{w}}$, les vari{\'e}t{\'e}s de Bott-Samelson
  ayant une structure plus simple que les vari{\'e}t{\'e}s de Schubert.

\end{rema}

\chapter{Cohomologie {\'e}quivariante} \label{cohomologie}

\section{Pr{\'e}liminaires}


 Soit $U_{\mathbb{C}}$ un tore complexe 
 d'alg{\`e}bre de Lie $\mathfrak{u}_{\mathbb{C}}$
et soit $U   \subset U_{\mathbb{C}}$ le tore compact maximal de
 $U_{\mathbb{C}}$.

  On note 
 $\mathfrak{u} \subset \mathfrak{u}_{\mathbb{C}} $ l'alg{\`e}bre de Lie de $U$, et
 $S(\mathfrak{u}_{\mathbb{C}}^{*})$ l'alg{\`e}bre sym{\'e}trique de
 $\mathfrak{u}_{\mathbb{C}}^{*}$. C'est l'alg{\`e}bre des fonctions polynomiales {\`a}
 coefficients complexes sur $\mathfrak{u}$.

\medskip

Soit $M$ un espace topologique muni d'une action continue de $U$. Notons
$\mathbb{E}U \rightarrow \mathbb{B}U$ le fibr{\'e} universel de $U$ et $ \mathbb{E}U
\times^{U}M$ l'espace topologique obtenu en quotientant  $ \mathbb{E}U\times M$
par l'action de $U$ d{\'e}finie par $u(p,m)=(pu^{-1},um)$ pour tous $u \in
\mathbb{E}U$, $m \in M$ et $u \in U$. Pour tout anneau commutatif $A$, on
d{\'e}finit la cohomologie $U$-{\'e}quivariante de $X$ {\`a} coefficients dans $A$, not{\'e}e
$H_{U}^{*}(M,A)$, comme la cohomologie singuli{\`e}re de $\mathbb{E}U\times^{U} M$ {\`a}
coefficients dans $A$. C'est une $H_{U}^{*}(pt,A)$-alg{\`e}bre. Dans toute la suite 
on prendra pour $A$ le corps des nombres complexes, et on notera $H_{U}^{*}(M)$
au lieu de $H_{U}^{*}(M,\mathbb{C})$. 
Dans ce cas, $H_{U}^{*}(pt)$
s'identifie {\`a} $S(\mathfrak{u}_{\mathbb{C}}^{*})$.

Toute application $g : M \rightarrow N$  continue  et $U$-{\'e}quivariante entre deux
espaces topologiques munis d'une action continue de $U$ induit un morphisme
de $S(\mathfrak{u}_{\mathbb{C}}^{*}) $-alg{\`e}bres $g^* : H_U^*(N) \rightarrow H_U^*(M)$.

\subsection{Int{\'e}gration et restriction aux points fixes} \label{integrationrestriction}

On suppose que $M$ est un $CW$ complexe orient{\'e} de dimension paire $n$, 
ne poss{\'e}dant pas de cellule de codimension $1$. Alors, la fibration
$U$-{\'e}quivariante triviale : $M \rightarrow \{ pt \}$ et l'int{\'e}gration sur les
fibres (voir \cite{aa} pour ces d{\'e}finitions) permettent de d{\'e}finir un
homomorphisme de  $S(\mathfrak{u}_{\mathbb{C}}^{*})$-modules gradu{\'e}s : 
$$\int_{M} : H_U^*(M) \rightarrow S(\mathfrak{u}_{\mathbb{C}}^{*}),$$
de degr{\'e} $-n$.

 Si $N$ est une sous-vari{\'e}t{\'e} $U$-invariante de $M$ admettant une telle structure,
on peut d{\'e}finir le morphisme de $S(\mathfrak{u}_{\mathbb{C}}^{*})$-modules
gradu{\'e}s : 
$$\int_{N} :  H_U^*(M) \rightarrow S(\mathfrak{u}_{\mathbb{C}}^{*}), $$
 obtenu en restreignant $\mu$ {\`a} $N$ puis
en int{\'e}grant sur $N$.

\medskip

Si on note $M^U$ l'ensemble des points fixes de $U$ dans $M$, l'inclusion
$M^U\subset M$ fournit un morphisme $i_U^* : H_{U}^*(M) \rightarrow
H_{U}^*(M^U) $ appel{\'e} restriction aux points fixes. Si $M^U$ est discret, 
$H_{U}^*(M^U) $ s'identifie {\`a} la $S( \mathfrak{u}_{\mathbb{C}}^{*} )$-alg{\`e}bre
des fonctions sur $M^U$ {\`a} valeurs dans
$S(\mathfrak{u}_{\mathbb{C}}^*)$ de degr{\'e} born{\'e}
munie de l'addition et de la multiplication point par point. On notera cette
alg{\`e}bre $F_{b}(M^U;S(\mathfrak{u}_{\mathbb{C}}^*))$.

\begin{rema}

Si $M^U$ est fini, comme c'est le cas pour les tours de Bott et les vari{\'e}tes de
drapeaux dans le cas fini, toute fonction sur $M^U$ {\`a} valeurs dans
$S(\mathfrak{u}_{\mathbb{C}}^*)$ est de degr{\'e} born{\'e}.

\end{rema}

Si $M$ admet une structure de $CW$-complexe $U$-{\'e}quivariant (non n{\'e}cessairement
fini) orient{\'e} $M=\coprod_{f \in \mathcal{F}}M_f$, ne comportant aucune cellule de
dimension impaire et seulement un nombre fini de cellules en chaque dimension
paire, et si $M^U$ est discret, les propositions $2.5.1$ et $2.6.1$ de \cite{aa} 
nous donnent la proposition suivante :

\begin{prop} \label{propbasecohomologie}

\indent

\begin{enumerate}

\item[$(i)$] La cohomologie $U$-{\'e}quivariante de $M^{U}$ s'identifie {\`a}
$F_{b}(M^U;S(\mathfrak{u}_{\mathbb{C}}^*))$.

\item[$(ii)$] La restriction aux points fixes $i_{U}^*$ : $H_{U}^*(M)
 \rightarrow F(M^U;S(\mathfrak{u}_{\mathbb{C}}^*))$ est injective.
    
\item[$(iii)$] La cohomologie $U$-{\'e}quivariante de $M$ est un 
$S(\mathfrak{u}_{\mathbb{C}}^*)$-module libre
qui admet comme base la famille $\{{\hat{\sigma}}^{U}_{f}\}_{f \in
\mathcal{F}}$ d'{\'e}l{\'e}ments homog{\`e}nes de degr{\'e} $ dim (M_f)$ 
caract{\'e}ris{\'e}s par : 

$$\int_{\overline{M_{f'}}}{\hat{\sigma}}^{U}_{f}=\delta_{f', f},$$
o{\`u} $\delta$ d{\'e}signe le symbole de Kronecker.

\end{enumerate}

\end{prop}

\subsection{Formule de localisation}

 Si $M$ est une vari{\'e}t{\'e} diff{\'e}rentiable munie d'une action diff{\'e}rentiable de $U$, 
la cohomologie $U$-{\'e}quivariante de 
$M$ s'identifie {\`a} la cohomologie du complexe des formes
diff{\'e}rentielles $\mu(Y)$ {\`a} coefficients complexes sur $M$ d{\'e}pendant
polynomialement de $Y \in \mathfrak{u}$ et v{\'e}rifiant une condition de
$U$-{\'e}quivariance {\'e}vidente muni de la diff{\'e}rentielle $D=d-2i\pi\, \iota(Y)$ 
(voir \cite{aa}).

 De plus, si $M$ est compacte, connexe, et orient{\'e}e,
pour toute 
forme $\mu(Y)$, on peut d{\'e}finir le
polyn{\^o}me $\int_{M}\mu(Y) \in S( \mathfrak{u}_{\mathbb{C}}^{*} )$ 
qui induit un morphisme en cohomologie. Ce morphisme
s'identifie {\`a} l'int{\'e}gration topologique d{\'e}finie dans la
section~\ref{integrationrestriction} (voir \cite{aa}), et 
on notera donc encore $\int_{M} : H_{U}^*(M) \rightarrow 
S( \mathfrak{u}_{\mathbb{C}}^{*} )  $ ce
morphisme de $S( \mathfrak{u}_{\mathbb{C}}^{*} )$-modules.

\medskip

On suppose que $M$ est une vari{\'e}t{\'e} diff{\'e}rentiable compacte  
de dimension $2n$ munie d'une structure presque complexe pr{\'e}serv{\'e}e par l'action
de $U$.
On suppose de plus que $M^U$ est fini. 
En chaque point fixe $m \in M^U$, l'action de $U$ sur $M$ induit
une action de $U$ sur $T_mM$, l'espace tangent {\`a} $M$ en $m$ muni de sa structure
complexe. On note
$(\alpha_1^m, \ldots ,\alpha_n^m) \in (i\mathfrak{u}^*)^n 
 \subset (\mathfrak{u}_{\mathbb{C}}^*)^n$ les poids de la
repr{\'e}sentation de $U$ dans $T_mM$. Dans ce cas, la formule de
localisation (voir \cite{bv}) s'{\'e}crit de la mani{\`e}re suivante :

\begin{prop} \label{formuledelocalisation}

Pout toute forme $\hat{\sigma}^U \in H_U^*(M)$, on a : 

$$\int_{M}\hat{\sigma}^U=\sum_{m\in M^U}\frac{i_U^*(\hat{\sigma}^U)(m)}
{\prod_{1 \leq i \leq n}\alpha_i^m  }.$$

\end{prop}

\subsection{Evaluation {\`a} l'origine et restriction}

La projection canonique $  \mathbb{E}U \times M \rightarrow
\mathbb{E}U\times^{U}M$ induit un homomorphisme de $\mathbb{C}$-alg{\`e}bre 
$v_0 : H_U^*(M) \rightarrow
H^*(M)$ appel{\'e} {\'e}valuation {\`a} l'origine. 

On suppose  que  $M$  admet une d{\'e}composition
cellulaire  $U$-{\'e}quivariante  orient{\'e}e $M=\coprod_{f \in \mathcal{F}}M_f$, ne comportant 
aucune cellule de dimension impaire et seulement un nombre fini de cellules en
chaque dimension paire. On note $\hat{\sigma}_{f}^U$ la base de
 $H_U^*(M)$ (en tant que $S(\mathfrak{u}_{\mathbb{C}}^*)$-module) donn{\'e}e par la
proposition~\ref{propbasecohomologie}.

La d{\'e}composition  $M=\coprod_{f \in \mathcal{F}}M_f$ fournit une base $y_f$ sur
$\mathbb{C}$ de la
cohomologie ordinaire de $M$. La proposition $2.5.1$ et le lemme $5.5.1$ de
\cite{aa} nous donnent le r{\'e}sultat suivant : 

\begin{prop} \label{evaluation}

L'{\'e}valuation {\`a} l'origine $v_0$ est l'homomorphisme de $\mathbb{C}$-alg{\`e}bres : 
$ H_U^*(M) \rightarrow H^*(M)$ d{\'e}fini par :
$$\forall f
\in \mathcal{F}, \, \mbox{et } \forall h \in S(\mathfrak{u}_{\mathbb{C}}^*), \,
v_0(h\hat{\sigma}^U_f)=h(0)y_f.$$

\end{prop}

Soit $V_{\mathbb{C}}$ un sous-tore de $U_{\mathbb{C}}$ d'alg{\`e}bre de Lie
  $\mathfrak{v}_{\mathbb{C}} \subset\mathfrak{u}_{\mathbb{C}} $ et soit
  $V=V_{\mathbb{C}}\cap U$ le tore compact maximal de $V$. 
On a alors un morphisme de restriction
$p_V^U : H_U^*(M) \rightarrow H_V^*(M)$. L'int{\'e}gration commutant {\`a} la
  restriction, on obtient facilement la proposition suivante : 

\begin{prop} \label{pUV}

L'application  $p_V^U$ est l'homomorphisme de 
$S(\mathfrak{u}_{\mathbb{C}}^*)$-alg{\`e}bres 
$ H_U^*(M) \rightarrow H_V^*(M)$ d{\'e}fini par :  
$$\forall f \in \mathcal{F},\, \mbox{et } \forall h \in 
S(\mathfrak{u}_{\mathbb{C}}^*), \,\,\,\, p_V^U(h\hat{\sigma}_{f}^U)=
\tilde{p}_V^U(h)\hat{\sigma}_{f}^V
,$$
o{\`u} l'application $\tilde{p}_V^U : S(\mathfrak{u}_{\mathbb{C}}^*) \rightarrow 
S(\mathfrak{v}_{\mathbb{C}}^*)$ est d{\'e}duite de l'inclusion 
$\mathfrak{v}_{\mathbb{C}} \subset\mathfrak{u}_{\mathbb{C}} $.

\end{prop}

\section{Cohomologie {\'e}quivariante des tours de Bott}

On reprend les notations de la section \ref{tours}. Soit $N \geq 1$ un entier
  naturel, $\mathcal{E}= \{0,1\}^N$, et soit $C=\{c_{i,j}\} _{1 \leq
  i < j \leq N}$ une liste d'entiers. Soit $Y=Y_{C}$ la tour de Bott associ{\'e}e {\`a}
  cette liste.

\subsection{Restrictions aux points fixes}

La d{\'e}composition $Y=\coprod_{\epsilon \in
  \mathcal{E}}Y_{\epsilon}$ munit $\Gamma$ d'une structure de
  $CW$-complexe orient{\'e} $D$-{\'e}quivariant fini o{\`u} toutes les cellules sont de
  dimension
 paire; de plus l'ensemble des points fixes de l'action de $D$ sur
  $\Gamma$ est fini et s'identifie {\`a} $\mathcal{E}$. On peut donc appliquer la
  proposition~\ref{propbasecohomologie}, et on a le r{\'e}sultat suivant :

\begin{prop}

\indent

\begin{enumerate}

\item[$(i)$] La cohomologie $D$-{\'e}quivariante de $Y^{D}$ s'identifie {\`a}
$F_b(\mathcal{E};S(\mathfrak{d}_{\mathbb{C}}^*))$.

\item[$(ii)$] La restriction aux points fixes $i_{D}^*$ : $H_{D}^*(Y)
 \rightarrow F_b(\mathcal{E};S(\mathfrak{d}_{\mathbb{C}}^*))$ est injective.
    
\item[$(iii)$] La cohomologie $D$-{\'e}quivariante de $Y$ est un 
$S(\mathfrak{d}_{\mathbb{C}}^*)$-module libre
qui admet comme base la famille $\{{\hat{\sigma}}^{D}_{\epsilon}\}_{\epsilon \in
\mathcal{E}}$ d'{\'e}l{\'e}ments homog{\`e}nes de degr{\'e} $2l(\epsilon)$ caract{\'e}ris{\'e}s par : 

$$\int_{\overline{Y_{\epsilon'}}}{\hat{\sigma}}^{D}_{\epsilon}=\delta_{\epsilon', \epsilon}.$$
   
\end{enumerate}

\end{prop}

\begin{defi} \label{defsigma}
    
Pour $\epsilon \in \mathcal{E}$, on d{\'e}finit $\sigma^{D}_{\epsilon} \in
F_b(\mathcal{E};S(\mathfrak{d}_{\mathbb{C}}^*))$
par :

$$\left\{ \begin{array}{ll} \sigma^{D}_{\epsilon}(\epsilon') = 
\displaystyle{ \prod_{i \in \pi_{+}(\epsilon)}\lambda_{i}(\epsilon')}
 & si \hspace{0,15 cm} \epsilon \leq \epsilon', \\ 
 \sigma^{D}_{\epsilon}(\epsilon') = 0
 & sinon.
\end{array}\right.$$

\end{defi}

On a alors le th{\'e}or{\`e}me suivant :

\begin{theo} \label{cohomologietours}
    
Pour tout $\epsilon \in \mathcal{E}$, on a : 
$$i_{D}^*({\hat{\sigma}^{D}}_{\epsilon})=\sigma^{D}_{\epsilon}.$$
   
Les restrictions aux points fixes sont donc des produits de formes lin{\'e}aires.

\end{theo}

\begin{proof}

En utilisant les propositions~\ref{formuledelocalisation} et
\ref{decompositiontours} et le lemme~\ref{pointsfixestours} , on
obtient pour tout  ${\hat{\sigma}^{D}} \in H_{D}^*(X)$ : 

\begin{equation} \label{localisation}
\int_{\overline{Y_{\epsilon}}} {\hat{\sigma}^{D}}=  \sum_{\epsilon' \leq
\epsilon}\frac{i_{D}^*({\hat{\sigma}^{D}})(\epsilon')}{\prod_{i
  \in \pi_{+}(\epsilon)}\lambda_{i}(\epsilon')}.  \end{equation}

Soit $\epsilon_{0} \in \mathcal{E}$, et soit
$\tilde{\sigma}^{D}_{\epsilon_{0}}=i_{D}^*({{\hat{\sigma}}^{D}}_{\epsilon_{0}})$. Montrons par
r{\'e}currence sur la longueur de $\epsilon$ que pour tout $\epsilon \in
\mathcal{E}$,
$\tilde{\sigma}^{D}_{\epsilon_{0}}
(\epsilon)=\sigma^{D}_{\epsilon_{0}}(\epsilon)$.

Gr{\^a}ce {\`a} la
formule~\ref{localisation} et {\`a} la caract{\'e}risation de ${\hat{\sigma}^{D}}_{\epsilon_{0}}$,
on d{\'e}montre facilement par r{\'e}currence sur $l(\epsilon)$ que si $\epsilon$ 
n'est pas plus grand que $\epsilon_{0}$, on
a bien $\tilde{\sigma}^{D}_{\epsilon_{0}}(\epsilon)=0$. 

On peut donc se limiter
au cas o{\`u} $\epsilon_{0} \leq \epsilon$. Si $\epsilon=\epsilon_{0}$, la
formule~\ref{localisation} et le fait que
$\int_{\overline{Y_{\epsilon_{0}}}}{\hat{\sigma}^{D}}_{\epsilon_{0}}=1$ nous
donne bien
$\tilde{\sigma}^{D}_{\epsilon_{0}}(\epsilon_{0})=\sigma^{D}_{\epsilon_{0}}(\epsilon_{0})$.

  Soit $\epsilon > \epsilon_{0}$. On suppose
le r{\'e}sultat v{\'e}rifi{\'e} pour tout $\epsilon'$ de longueur strictement plus
petite que $\epsilon$, on applique la formule~\ref{localisation} et le fait
que $\int_{\overline{Y_{\epsilon}}}{\hat{\sigma}^{D}}_{\epsilon_{0}}=0$  pour
obtenir : 

$$\sum_{\epsilon_{0} \leq  \epsilon' <
\epsilon}\!\! \frac{\prod_{i
  \in \pi_{+}(\epsilon_{0})}\lambda_{i}(\epsilon')}{\prod_{i
  \in
  \pi_{+}(\epsilon)}\lambda_{i}(\epsilon')}+\frac{
 \tilde{\sigma}^{D}_{\epsilon_{0}}(\epsilon)}{\prod_{i
  \in \pi_{+}(\epsilon)}\lambda_{i}(\epsilon)}  = 0,$$
d'o{\`u} : 
$$\sum_{\epsilon_{0} \leq \epsilon' <
\epsilon}\!\! \frac{1}{\prod_{i
  \in
  \pi_{+}(\epsilon)\setminus \pi_{+}(\epsilon_{0})}\lambda_{i}(\epsilon')}+
\tilde{\sigma}^{D}_{\epsilon_{0}}(\epsilon){\prod_{i
  \in \pi_{+}(\epsilon)}\lambda_{i}(\epsilon)}  = 0.$$

\bigskip

Si on pose $\tilde{\epsilon}=\epsilon-(j)$, o{\`u} $j$ est le plus grand
{\'e}l{\'e}ment de $\pi_{+}(\epsilon) \setminus  \pi_{+}(\epsilon_{0})$, on a alors : 
$$\sum_{
 \tiny \begin{array}{cc}  \epsilon_{0} \leq \epsilon' <
\epsilon \\
 \epsilon' \neq \tilde{\epsilon} 
\end{array}}\!\!\!\!\!
\frac{1}{\prod_{i
  \in
  \pi_{+}(\epsilon)\setminus \pi_{+}(\epsilon_{0})}\lambda_{i}(\epsilon')}
  \,\, -
\frac{1}{\prod_{i
  \in
  \pi_{+}(\epsilon)\setminus \pi_{+}(\epsilon_{0})}\lambda_{i}(\epsilon)}+
\frac{\tilde{\sigma}^{D}_{\epsilon_{0}}(\epsilon)}{\prod_{i
  \in \pi_{+}(\epsilon)}\lambda_{i}(\epsilon)}
  = 0.$$

En effet, comme $j$ est le plus grand {\'e}l{\'e}ment de $\pi_{+}(\epsilon) \setminus
\pi_{+}(\epsilon_{0})$, pour tout $i \in \pi_{+}(\epsilon) \setminus
\pi_{+}(\epsilon_{0})$, $i \neq j$,  
$\lambda_{i}(\epsilon)=\lambda_{i}(\tilde{\epsilon})$ et 
$\lambda_{j}(\epsilon)=-\lambda_{j}(\tilde{\epsilon})$. Pour les m{\^e}mes raisons, on
s'aper{\c c}oit, en distinguant les termes qui ont un $1$ en $j${\`e}me position
et ceux qui ont un $0$ en $j${\`e}me position, que la premi{\`e}re somme est
nulle, et on obtient alors bien
$\tilde{\sigma}^{D}_{\epsilon_{0}}(\epsilon)=
\prod_{i \in \pi_{+}(\epsilon_{0})}\lambda_{i}(\epsilon)$.

\end{proof}

\begin{exem}

Dans le tableau ci dessous, on explicite ces formules 
 dans le cas de la surface de Hirzebruch $H_{-1}=Y_{\{-1\}}$ :

\begin{figure}[ht]
 \begin{center}
 \setlength{\unitlength}{0.00062500in}
\begingroup\makeatletter\ifx\SetFigFont\undefined%
\gdef\SetFigFont#1#2#3#4#5{%
  \reset@font\fontsize{#1}{#2pt}%
  \fontfamily{#3}\fontseries{#4}\fontshape{#5}%
  \selectfont}%
\fi\endgroup%
{\renewcommand{\dashlinestretch}{30}
\begin{picture}(6099,3039)(0,-10)
\path(612,2412)(6012,2412)
\path(6012,2412)(6012,12)
\path(1212,2412)(1212,3012)
\path(1212,3012)(6012,3012)
\path(6012,3012)(6012,2412)
\path(2412,3012)(2412,12)
\path(3612,3012)(3612,12)
\path(4812,3012)(4812,12)
\path(1212,2412)(1212,12)
\path(612,1812)(6012,1812)
\path(612,1212)(6012,1212)
\drawline(6087,612)(6087,612)
\path(12,2412)(687,2412)
\path(12,1812)(612,1812)
\path(12,1212)(612,1212)
\path(12,612)(6012,612)
\path(12,12)(6012,12)
\path(12,2412)(12,12)
\put(1737,2037){\makebox(0,0)[lb]{\smash{{{\SetFigFont{9}{10.8}{\rmdefault}{\mddefault}{\updefault}$1$}}}}}
\put(2937,2037){\makebox(0,0)[lb]{\smash{{{\SetFigFont{9}{10.8}{\rmdefault}{\mddefault}{\updefault}$1$}}}}}
\put(4137,2037){\makebox(0,0)[lb]{\smash{{{\SetFigFont{9}{10.8}{\rmdefault}{\mddefault}{\updefault}$1$}}}}}
\put(5337,2037){\makebox(0,0)[lb]{\smash{{{\SetFigFont{9}{10.8}{\rmdefault}{\mddefault}{\updefault}$1$}}}}}
\put(1737,1437){\makebox(0,0)[lb]{\smash{{{\SetFigFont{9}{10.8}{\rmdefault}{\mddefault}{\updefault}$0$}}}}}
\put(1737,837){\makebox(0,0)[lb]{\smash{{{\SetFigFont{9}{10.8}{\rmdefault}{\mddefault}{\updefault}$0$}}}}}
\put(1737,237){\makebox(0,0)[lb]{\smash{{{\SetFigFont{9}{10.8}{\rmdefault}{\mddefault}{\updefault}$0$}}}}}
\put(4137,1437){\makebox(0,0)[lb]{\smash{{{\SetFigFont{9}{10.8}{\rmdefault}{\mddefault}{\updefault}$0$}}}}}
\put(2937,1437){\makebox(0,0)[lb]{\smash{{{\SetFigFont{9}{10.8}{\rmdefault}{\mddefault}{\updefault}$\lambda_1$}}}}}
\put(5337,1437){\makebox(0,0)[lb]{\smash{{{\SetFigFont{9}{10.8}{\rmdefault}{\mddefault}{\updefault}$\lambda_1$}}}}}
\put(2937,837){\makebox(0,0)[lb]{\smash{{{\SetFigFont{9}{10.8}{\rmdefault}{\mddefault}{\updefault}$0$}}}}}
\put(1587,2637){\makebox(0,0)[lb]{\smash{{{\SetFigFont{9}{10.8}{\rmdefault}{\mddefault}{\updefault}$(0,0)$}}}}}
\put(2787,2637){\makebox(0,0)[lb]{\smash{{{\SetFigFont{9}{10.8}{\rmdefault}{\mddefault}{\updefault}$(1,0)$}}}}}
\put(3987,2637){\makebox(0,0)[lb]{\smash{{{\SetFigFont{9}{10.8}{\rmdefault}{\mddefault}{\updefault}$(0,1)$}}}}}
\put(5187,2637){\makebox(0,0)[lb]{\smash{{{\SetFigFont{9}{10.8}{\rmdefault}{\mddefault}{\updefault}$(1,1)$}}}}}
\put(2937,237){\makebox(0,0)[lb]{\smash{{{\SetFigFont{9}{10.8}{\rmdefault}{\mddefault}{\updefault}$0$}}}}}
\put(4137,237){\makebox(0,0)[lb]{\smash{{{\SetFigFont{9}{10.8}{\rmdefault}{\mddefault}{\updefault}$0$}}}}}
\put(5037,837){\makebox(0,0)[lb]{\smash{{{\SetFigFont{9}{10.8}{\rmdefault}{\mddefault}{\updefault}\,$\lambda_1+\lambda_2$}}}}}
\put(4887,237){\makebox(0,0)[lb]{\smash{{{\SetFigFont{9}{10.8}{\rmdefault}{\mddefault}{\updefault}$\lambda_1(\lambda_1+\lambda_2)$}}}}}
\put(4137,837){\makebox(0,0)[lb]{\smash{{{\SetFigFont{9}{10.8}{\rmdefault}{\mddefault}{\updefault}$\lambda_2$}}}}}
\put(312,1437){\makebox(0,0)[lb]{\smash{{{\SetFigFont{9}{10.8}{\rmdefault}{\mddefault}{\updefault}$\sigma^D_{(1,0)}$}}}}}
\put(312,2037){\makebox(0,0)[lb]{\smash{{{\SetFigFont{9}{10.8}{\rmdefault}{\mddefault}{\updefault}$\sigma^D_{(0,0)}$}}}}}
\put(312,837){\makebox(0,0)[lb]{\smash{{{\SetFigFont{9}{10.8}{\rmdefault}{\mddefault}{\updefault}$\sigma^D_{(0,1)}$}}}}}
\put(312,237){\makebox(0,0)[lb]{\smash{{{\SetFigFont{9}{10.8}{\rmdefault}{\mddefault}{\updefault}$\sigma^D_{(1,1)}$}}}}}
\end{picture}
} 
 \caption{Restrictions aux points fixes pour la surface de Hirzebruch $H_{-1}$}
 \label{pointsfixeshirz}
 \end{center}
 \end{figure}

\end{exem}



\subsection{Structure multiplicative}

L'alg{\`e}bre $H_D^*(Y)$ est une alg{\`e}bre de polyn{\^o}mes (o{\`u} les ind{\'e}termin{\'e}es
repr{\'e}sentent des classes de degr{\'e} $2$) quotient{\'e}e par des relations
de degr{\'e} $2$ que l'on va expliciter.

 On pose 
$\sigma^{D}_{i}=\sigma^{D}_{(i)}$ et 
${\hat{\sigma}^{D}}_{i}={\hat{\sigma}^{D}}_{(i)}$. On a alors : 

\begin{theo}  \label{produittours}
    
Pour tout $\epsilon \in \mathcal{E}$, on a : 

$${\hat{\sigma}^{D}}_{\epsilon}=\prod_{i \in
  \pi_{+}(\epsilon)}\!\!\!{\hat{\sigma}^{D}}_{i}.$$

De plus, on a les formules de multiplication suivantes : 

$$\left\{ \begin{array}{ll} {\hat{\sigma}^{D}}_{i}\hspace{0,1 cm}{\hat{\sigma}^{D}}_{\epsilon}=
{\hat{\sigma}^{D}}_{\epsilon+(i)}
 & si \hspace{0,15 cm} i\in \pi_{-}(\epsilon), \\ 
\displaystyle{
  {\hat{\sigma}^{D}_{i}\hspace{0,1 cm}{\hat{\sigma}^{D}}_{\epsilon}=\sigma^{D}_{i}
(\epsilon)\hspace{0,1 cm}{\hat{\sigma}^{D}}_{\epsilon}
+\sum_{j<i,j\in
  \pi_{-}(\epsilon)}\!\!\!\!\!\!c_{j,i}(\epsilon)\hspace{0,1 cm}{\hat{\sigma}^{D}}_{\epsilon}}
\hspace{0,1 cm}{\hat{\sigma}^{D}}_{j}}
 & si \hspace{0,15 cm} i \in \pi_{+}(\epsilon).\end{array} \right. $$

\end{theo}

\begin{proof}

Par injectivit{\'e} de la restriction aux points fixes, il suffit de
d{\'e}montrer ces formules pour $\sigma^{D}_{i}$ et $\sigma^{D}_{\epsilon}$.

La premi{\`e}re formule se voit imm{\'e}diatement sur la d{\'e}finition des
$\sigma^{D}_{\epsilon}$.

Pour la deuxi{\`e}me formule, on {\'e}crit le produit 
$\sigma^{D}_{i}\hspace{0,1 cm}\sigma^{D}_{\epsilon}$ 
 sous la forme : 

\begin{equation} \label{prod}
\sigma^{D}_{i}\hspace{0,1 cm}  \sigma^{D}_{\epsilon}=
  \sum_{\epsilon' \in \mathcal{E}}C_{\epsilon'} \sigma^{D}_{\epsilon'}.
\end{equation}

Soit $\epsilon' \in \mathcal{E}$ tel que $\epsilon' \not\geq \epsilon$. On
montre que pour tout $\epsilon'' \leq \epsilon'$, $C_{\epsilon''}=0$. Pour cela,
on proc{\`e}de par r{\'e}currence sur $l(\epsilon'')$ en {\'e}valuant l'{\'e}galit{\'e}~\ref{prod}
en $\epsilon''$. On obtient ainsi que $C_{\epsilon'}=0$ si $\epsilon' \not\geq
\epsilon$. 

De plus, pour des raisons de degr{\'e}, $C_{\epsilon'}=0$ si
$l(\epsilon')>l(\epsilon)+ 1$, et l'{\'e}galit{\'e}~\ref{prod} s'{\'e}crit donc sous la forme : 

$$\sigma^{D}_{i}\hspace{0,1 cm}  \sigma^{D}_{\epsilon}=C_{i}
 \hspace{0,1 cm} \sigma^{D}_{\epsilon}+\sum_{j \neq i,j\in
  \pi_{-}(\epsilon)}\!\!\!\!\!\!\!\!C_{j} \,  \sigma^{D}_{\epsilon+(j)}.$$

En {\'e}valuant en $\epsilon$, on trouve $C_{i}$, puis en {\'e}valuant en
$\epsilon+(j)$, on trouve $C_{j}$.

\end{proof}

En particulier, on a les relations suivantes :

$$  (\hat{\sigma}^{D}_{i})^2=\lambda_{i} \hspace{0,1 cm}  {\hat{\sigma}^{D}}_{i} -  
 \sum_{j<i}c_{j,i} \hspace{0,1 cm}  {\hat{\sigma}^{D}}_{i} \hspace{0,1 cm} 
 {\hat{\sigma}^{D}}_{j}.$$

Si on pose $x_i=\hat{\sigma}^D_i$, on a donc le th{\'e}or{\`e}me suivant :

\begin{theo} \label{carre}

L'alg{\`e}bre $H_D^*(Y)$ est l'alg{\`e}bre de polyn{\^o}mes
$S(\mathfrak{d}_{\mathbb{C}}^*)[x_1, \ldots , x_N]$ (o{\`u} les ind{\'e}termin{\'e}es sont
de degr{\'e} $2$) quotient{\'e}e par les relations :
$$  x_i^2=\lambda_{i} x_i -  
 \sum_{j<i}c_{j,i} x_ix_j.$$

\end{theo}

\begin{exem} \label{exempleA2Bott}

On consid{\`e}re la liste d'entiers $C=\{c_{i,j}\}_{1\leq i<j \leq 3}$ d{\'e}finie par : 
$$\left\{ \begin{array}{ll}
c_{1,2}= c_{2,3}=-1,\\ 
 c_{1,3}=2, 
\end{array}\right.$$ 
et on pose $Y=Y_C$. 
 Alors,
$H_D^*(Y)$ s'identifie au quotient de l'alg{\`e}bre de polyn{\^o}mes 
$\mathbb{C}[\lambda_1,\lambda_2, \lambda_3 ][x_1,x_2,x_3]$ par les relations : 
\begin{eqnarray*} 
x_1^2 & = & \lambda_1 x_1 \\
x_2^2 & = & \lambda_2 x_2 + x_1x_2 \\
x_3^2 & = & \lambda_3 x_3 + x_2x_3  -2 x_1x_3
\end{eqnarray*}

\end{exem}

\begin{exem} \label{exempleB2Bott}

On consid{\`e}re la liste d'entiers $C=\{c_{i,j}\}_{1\leq i<j \leq 4}$ d{\'e}finie par : 
$$\left\{ \begin{array}{lll}
c_{1,2}= c_{1,4}=  c_{3,4}   =-2,\\ 
 c_{1,3}=c_{2,4}  =2,\\
 c_{2,3}  =-1,
\end{array}\right.$$ 
et on pose $Y=Y_C$. 
 Alors, $H_D^*(Y)$ s'identifie au quotient de l'alg{\`e}bre de polyn{\^o}mes 
$\mathbb{C}[ \lambda_1,\lambda_2, \lambda_3 ,\lambda_4 ][x_1,x_2,x_3,x_4]$ 
par les relations : 
\begin{eqnarray*} 
x_1^2 & = & \lambda_1 x_1 \\
x_2^2 & = & \lambda_2 x_2 + 2x_1x_2 \\
x_3^2 & = & \lambda_3 x_3 + x_2x_3  -2 x_1x_3 \\
x_4^2 & = & \lambda_4 x_4 + 2x_3x_4  -2 x_2x_4+2x_1x_4
\end{eqnarray*}

\end{exem}

\subsection{Cohomologie ordinaire}

La d{\'e}composition cellulaire $Y=\coprod_{\epsilon \in \mathcal{E}}Y_{\epsilon}$
fournit une base $(y_{\epsilon})_{\epsilon \in \mathcal{E}}$ de la 
cohomologie ordinaire de $Y$. Si on pose $y_i=y_{(i)}$, la
proposition~\ref{evaluation} et le th{\'e}or{\`e}me~\ref{carre} nous donnent 
la description suivante de la cohomologie ordinaire de $Y$  : 

\begin{prop}

La cohomologie ordinaire de $Y$ est engendr{\'e}e par des {\'e}l{\'e}ments
$(y_i)_{1 \leq i \leq N}$ de degr{\'e} $2$ soumis aux relations : 

$$y_i^2 + \sum_{j<i}c_{j,i}y_iy_j=0.$$

\end{prop}

\section{Cohomologie {\'e}quivariante des vari{\'e}t{\'e}s de Bott-Samelson} \label{sectioncohomologieBS}

On reprend les notations de la section~\ref{BS}. On choisit $N$ racines simples 
 $\mu_1, \ldots, \mu_N$ 
 non n{\'e}cessairement distinctes, et on pose 
$\Gamma=\Gamma( \mu_1, \ldots, \mu_N )$. 
Pour $1 \leq i <j \leq N$, on pose $b_{i,j}=\mu_j(\mu_{i}^{\vee})$ et 
$B=\{b_{i,j} \}_{1 \leq i <j \leq N }$.

\subsection{Restrictions aux points fixes}

Les r{\'e}sultats de cette section ne
d{\'e}pendant pas de la structure complexe sur $\Gamma$, on peut identifier $\Gamma$
avec la tour de Bott $Y=Y_B$, la d{\'e}composition 
$\Gamma=\coprod_{\epsilon \in \mathcal{E}} \Gamma_{\epsilon}$ avec la
  d{\'e}composition $Y=\coprod_{\epsilon \in \mathcal{E}} Y_{\epsilon}$, et le point
  fixe $\epsilon \in \Gamma^T$ avec le point fixe $\epsilon \in Y^D$ {\`a} l'aide de 
  l'isomorphisme $\phi$ de la section~\ref{toursbott}.
 Le tore $T$ agit sur $\Gamma$ via son image  $S$, et  le tore $S$ s'identifie, par 
l'homomorphisme $\gamma$, {\`a} un
sous-tore du tore $D$ qui agit sur $Y$ (voir 
section~\ref{toursbott}).

Pour tout $\epsilon \in \mathcal{E}$, on d{\'e}finit $\hat{\sigma}_{\epsilon}^S
\in H_S^*(\Gamma)$ par $\hat{\sigma}_{\epsilon}^S=p_S^D(
 \hat{\sigma}_{\epsilon}^D )$. D'apr{\`e}s la proposition~\ref{pUV}, la famille 
$\{ \hat{\sigma}_{\epsilon}^S \}_{\epsilon \in \mathcal{E}}$ est une base du 
$S(\mathfrak{s}_{\mathbb{C}}^*)$-module libre $H_S^*(\Gamma)$ v{\'e}rifiant :
$$\forall (\epsilon,\epsilon') \in \mathcal{E}^2, \int_{\overline{\Gamma_{\epsilon'}}}
\hat{\sigma}_{\epsilon}^S=
\delta_{\epsilon', \epsilon}.$$

\medskip

 Comme $T$ agit sur $\Gamma$ via image dans $S$,
$H_T^*(\Gamma)=H_S^*(\Gamma)\otimes_{S(\mathfrak{s}_{\mathbb{C}}^*)}
S(\mathfrak{h}^*)$. 

\smallskip

Pour tout $\epsilon \in \mathcal{E}$, on d{\'e}finit 
alors $\hat{\sigma}_{\epsilon}^T
\in H_T^*(\Gamma)$ par $\hat{\sigma}_{\epsilon}^T=
 \hat{\sigma}_{\epsilon}^S \otimes 1$. Pour tout couple 
$(\epsilon,\epsilon')  \in \mathcal{E}^2$ : 
$$ \int_{\overline{\Gamma_{\epsilon'}}}
\hat{\sigma}_{\epsilon}^T=\int_{\overline{\Gamma_{\epsilon'}}}
\hat{\sigma}_{\epsilon}^S=
\delta_{\epsilon', \epsilon}.$$

\medskip

On obtient ainsi la proposition suivante :

\begin{prop} 

\indent

\begin{enumerate}

\item[$(i)$] La cohomologie $T$-{\'e}quivariante de $\Gamma^T$ s'identifie {\`a}  
$F_b(\mathcal{E};S(\mathfrak{h}^*))$.

\item[$(ii)$] La restriction aux points fixes $i_{T}^*$ : $H_{T}^*(\Gamma)
 \rightarrow F_b(\mathcal{E};S(\mathfrak{h}^*))$ est injective.
    
\item[$(iii)$] La cohomologie $T$-{\'e}quivariante de $\Gamma$ est un 
$S(\mathfrak{h}^*)$-module libre
qui admet comme base la famille $\{\hat{\sigma}_{\epsilon}^T\}_{\epsilon \in
\mathcal{E}}$ d'{\'e}l{\'e}ments homog{\`e}nes de degr{\'e} $2l(\epsilon)$ v{\'e}rifiant : 

$$\int_{\overline{\Gamma_{\epsilon'}}}\hat{\sigma}_{\epsilon}^T=
\delta_{\epsilon', \epsilon}.$$

\end{enumerate}

\end{prop}

\begin{defi}  \label{defcohomologieBS}
    
Pour $\epsilon \in \mathcal{E}$, on d{\'e}finit $\sigma_{\epsilon}^T \in
F_b(\mathcal{E};S(\mathfrak{h}^*))$
par :

$$\left\{ \begin{array}{ll} \sigma_{\epsilon}^T(\epsilon') = 
\displaystyle{ (-1)^{l(\epsilon)}\prod_{i \in \pi_{+}(\epsilon)}\alpha_{i}(\epsilon')}
 & si \hspace{0,15 cm} \epsilon \leq \epsilon', \\ 
 \sigma_{\epsilon}^T(\epsilon') = 0
 & sinon.
\end{array}\right.$$

\end{defi}

On a alors le th{\'e}or{\`e}me suivant :

\begin{theo} \label{cohomologieBS}
    
Pour tout $\epsilon \in \mathcal{E}$, on a : 

$$i_{T}^*(\hat{\sigma}_{\epsilon}^T)=\sigma_{\epsilon}^T.$$
   
\end{theo}

\begin{proof}

 On a le diagramme commutatif suivant :

 $$\xymatrix{
     H_{S}^*(\Gamma)\ar@{^{(}->}[dd]^{i_{S}^*} & &
     \ar[ll]^{p_S^D}H_{D}^*(Y)\ar@{^{(}->}[dd]^{i_{D}^*} \\
     \\
     F_b(\mathcal{E};S(\frak{s}_{\mathbb{C}}^*)) & & 
   \ar[ll]^{\tilde{\tau}}F_b(\mathcal{E};S(\frak{d}_{\mathbb{C}}^*)) } $$
o{\`u} l'application $\tilde{\tau} :  F_b(\mathcal{E};
S(\mathfrak{d}_{\mathbb{C}}^*)) \rightarrow 
 F_b(\mathcal{E};S(\mathfrak{s}_{\mathbb{C}}^*))$ est d{\'e}duite de
$\tau : \mathfrak{d}_{\mathbb{C}}^* \rightarrow \mathfrak{s}_{\mathbb{C}}^* $.

\smallskip

Pour tout couple $(\epsilon, \epsilon') \in \mathcal{E}^2$, on a donc : 
$$i_{T}^*(\hat{\sigma}^{T}_{\epsilon})(\epsilon')=
i_S^*(\hat{\sigma}^{S}_{\epsilon})(\epsilon')
= i_S^* p_S^D({\hat{\sigma}^{D}}_{\epsilon})(\epsilon')=
\tilde{\tau}(\sigma^{D}_{\epsilon})(\epsilon').$$

D'apr{\`e}s l'expression de $\sigma^{D}_{\epsilon}$ (voir d{\'e}finition~\ref{defsigma}), 
il suffit alors de prouver que pour tout $\epsilon \in \mathcal{E}$
et tout $i \in \{ 1, 2, \ldots, N \}$ :
$$\alpha_{i}(\epsilon)=-\tau(\lambda_{i}(\epsilon)).$$

 Comme 
$ s_{\nu_{s}}\nu_{s}=-\nu_{s}$, il faut v{\'e}rifier
que pour toute suite $\nu_{1}, \nu_{2}, \ldots ,\nu_{s}$ de racines simples, on a
la relation suivante : 
\begin{equation} \label{alphalambda}
s_{\nu_{1}}s_{\nu_{2}} \cdots s_{\nu_{s-1}}\nu_{s}=\nu_{s} +\sum_{1 \leq i
  <s}b_{i}\nu_{i}, \end{equation}
avec $\displaystyle{ b_{i} = 
\sum_{ {\tiny \begin{array}{cc} {i_{0}=i<i_{1}< \cdots < i_{m}=s},\\
 m>0
\end{array} }}}\!\!\!\!\!\!\!\!\!\!\!\!\!\!\!\!(-1)^m
a(\nu_{i_{0}},\nu_{i_{1}})a(\nu_{i_{1}},\nu_{i_{2}})
\cdots a(\nu_{i_{m-1}},\nu_{i_{m}}), $\\ 
o{\`u} $a(\nu, \nu')$ est le nombre de Cartan
associ{\'e} aux racines simples $\nu$ et $\nu'$ (i.e. $a(\nu, \nu')=
\nu'(\nu^{\vee}))$. Cette relation est une cons{\'e}quence
imm{\'e}diate de la d{\'e}finition des r{\'e}flexions simples {\`a} l'aide des nombres de Cartan
 et se d{\'e}montre par r{\'e}currence
sur $s$ : pour $s=1$, cela traduit la relation $s_{\nu}\nu'=\nu'-a(\nu,
\nu')\nu$.

\end{proof}

\subsection{Structure multiplicative et cohomologie ordinaire}

 On pose 
$\sigma_{i}^T=\sigma_{(i)}^T$ et 
$\hat{\sigma}_{i}^T=\hat{\sigma}_{(i)}^T$. La structure multiplicative de
$H_{T}^{*}(\Gamma)$ est donn{\'e}e par le th{\'e}or{\`e}me suivant : 

\begin{theo} \label{produitBS}
    
Pour tout $\epsilon \in \mathcal{E}$, on a : 

$$\hat{\sigma}_{\epsilon}^T=\prod_{i \in
  \pi_{+}(\epsilon)}\!\!\!\hat{\sigma}_{i}^T.$$

De plus, on a les formules de multiplication suivantes : 

$$\left\{ \begin{array}{ll} \hat{\sigma}_{i}^T\hat{\sigma}_{\epsilon}^T=
\hat{\sigma}_{\epsilon+(i)}^T
 & si \hspace{0,15 cm} i\in \pi_{-}(\epsilon), \\ 
\displaystyle{
  \hat{\sigma}_{i}^T\hat{\sigma}_{\epsilon}^T=\sigma_{i}^T(\epsilon)
\hat{\sigma}_{\epsilon}^T
+\sum_{j<i,j\in
  \pi_{-}(\epsilon)}\!\!\!\!\!\!\alpha_{j}^{i}(\epsilon)(\mu_j^{\vee}  )\,
\hat{\sigma}_{\epsilon}^T}\hat{\sigma}_{j}^T
 & si \hspace{0,15 cm} i\in \pi_{+}(\epsilon),\end{array}\right.$$

o{\`u} on a pos{\'e} $\alpha_{j}^{i}(\epsilon)=v_{j+1}^{i}(\epsilon)(\mu_{i})$.
   
\end{theo}

\begin{proof}

Pour d{\'e}montrer ces formules, on peut remplacer $\hat{\sigma}_i^T$ 
et $\hat{\sigma}_{\epsilon}^T$ par $\hat{\sigma}_i^S$
 et $\hat{\sigma}_{\epsilon}^S$. On utilise alors les r{\'e}sultats analogues pour
 les tours de Bott et l'homomorphisme de restriction $p_S^D$.

\medskip

La proposition~\ref{pUV} et le th{\'e}or{\`e}me~\ref{produittours} donnent imm{\'e}diatement la
premi{\`e}re formule.

\smallskip

Pour prouver la deuxi{\`e}me formule, il faut montrer : 
$$s_{j}\alpha_{j}^{i}(\epsilon)-\alpha_{j}^{i}(\epsilon)
=-\tilde{\tau}(c_{j,i}(\epsilon))\mu_{j}=-c_{j,i}(\epsilon)\mu_{j}.$$
 
Cette formule est une cons{\'e}quence imm{\'e}diate de l'{\'e}quation~\ref{alphalambda}.

\end{proof}

En particulier, on a l'expression suivante des carr{\'e}s des {\'e}l{\'e}ments de degr{\'e} $2$
:

$$  (\hat{\sigma}_{i}^T)^2=\mu_{i}\hat{\sigma}_{i}^T -  
 \sum_{j<i}b_{j,i}\hat{\sigma}_{i}^T\hat{\sigma}_{j}^T.$$

Si on pose $x_i=\hat{\sigma}_{i}^T$, on a donc le th{\'e}or{\`e}me suivant : 

\begin{theo} \label{carreBS}

L'alg{\`e}bre $H_T^*(\Gamma)$ est l'alg{\`e}bre de polyn{\^o}mes
$S(\mathfrak{h}^*)[x_1, \ldots , x_N]$ (o{\`u} les ind{\'e}termin{\'e}es sont de degr{\'e} $2$) 
quotient{\'e}e par les relations : 
$$  x_i^2=\mu_{i} x_i -  
 \sum_{j<i}b_{j,i} x_ix_j.$$

\end{theo}

De plus, on retrouve alors le r{\'e}sultat suivant
prouv{\'e} dans \cite{bs} : 

\begin{prop}

La cohomologie ordinaire de $\Gamma$ est engendr{\'e}e par des {\'e}l{\'e}ments
$(y_{i})_{1 \leq i \leq N}$ de degr{\'e} $2$ soumis aux relations : 

$$y_{i}^2 + \sum_{j<i}b_{j,i}y_{i}y_{j}=0.$$

\end{prop}

\begin{exem} \label{exempleA2BS}

Dans le cas $A_2$, on prend $\Gamma=\Gamma(\alpha_1,\alpha_2,\alpha_1)$. Alors,
$\Gamma$ s'identifie {\`a} $Y_C$ o{\`u} $C$ est d{\'e}finie dans
l'exemple~\ref{exempleA2Bott}. 

On a alors les m{\^e}mes relations dans $H_T^*(\Gamma)$ que dans $H_D^*(Y)$ 
en rempla{\c c}ant $\lambda_1$ par $\alpha_1$, 
$\lambda_2$ par $\alpha_2$, et $\lambda_3$ par $\alpha_1$.

L'alg{\`e}bre   
$H_T^*(\Gamma)$ s'identifie donc {\`a} l'alg{\`e}bre de polyn{\^o}mes 
$\mathbb{C}[ \alpha_1,\alpha_2 ][x_1,x_2,x_3]$ quotient{\'e}e par les relations : 
\begin{eqnarray*}
x_1^2 & = & \alpha_1 x_1 \\
x_2^2 & = & \alpha_2 x_2 + x_1x_2 \\
x_3^2 & = & \alpha_1 x_3 + x_2x_3  -2 x_1x_3
\end{eqnarray*}

\end{exem}

\begin{exem} \label{exempleB2BS}

Dans le cas $B_2$, on prend $\Gamma=\Gamma(\alpha_1,\alpha_2,\alpha_1,\alpha_2)$. Alors,
$\Gamma$ s'identifie {\`a} $Y_C$ o{\`u} $C$ est d{\'e}finie dans
l'exemple~\ref{exempleB2Bott}. 

On a alors les m{\^e}mes relations  dans $H_T^*(\Gamma)$ que dans $H_D^*(Y)$ 
en rempla{\c c}ant $\lambda_1$ par $\alpha_1$, 
$\lambda_2$ par $\alpha_2$, $\lambda_3$ par $\alpha_1$, et $\lambda_4$ par
$\alpha_2$.

L'alg{\`e}bre $H_T^*(\Gamma)$ s'identifie donc {\`a} l'alg{\`e}bre de polyn{\^o}mes 
$\mathbb{C}[ \alpha_1,\alpha_2 ]  [x_1,x_2,x_3,x_4]$ quotient{\'e}e par les relations :
\begin{eqnarray*}
x_1^2 & = & \alpha_1 x_1 \\
x_2^2 & = & \alpha_2 x_2 + 2x_1x_2 \\
x_3^2 & = & \alpha_1 x_3 + x_2x_3  -2 x_1x_3 \\
x_4^2 & = & \alpha_2 x_4 + 2x_3x_4  -2 x_2x_4+2x_1x_4
\end{eqnarray*}

\end{exem}

\section{Cohomologie {\'e}quivariante des vari{\'e}t{\'e}s de drapeaux} \label{cohomologievd}

 La d{\'e}composition
$X=\coprod_{w\in W}X_{w}$ munit $X$ d'une structure de $CW$-complexe orient{\'e} 
$T$-{\'e}quivariant ne comportant aucune cellule de dimension impaire et seulement
un nombre fini de cellules en chaque dimension paire. De plus,
l'ensemble $X^T \approx W$ des points fixes de $T$ dans $X$ {\'e}tant discret, on
peut {\`a} nouveau appliquer la proposition~\ref{propbasecohomologie} pour
obtenir le r{\'e}sultat suivant {\'e}tabli dans \cite{aa} et \cite{kkc} :

\begin{prop}

\indent

\begin{enumerate}

\item[$(i)$] La cohomologie $T$-{\'e}quivariante de $X^T$ s'identifie {\`a}
 $F_{b}(W;S(\mathfrak{h}^*))$.

\item[$(ii)$] La restriction aux points fixes $i_{T}^*$ : $H_{T}^*(X)
 \rightarrow F_{b}(W;S(\mathfrak{h}^*))$ est injective.
    
\item[$(iii)$] La cohomologie $T$-{\'e}quivariante de $X$ est un $S(\mathfrak{h}^*)$-module libre
qui admet comme base la famille $\{\hat{\xi}^{w}\}_{w \in W}$ d'{\'e}l{\'e}ments
homog{\`e}nes de degr{\'e} $2l(w)$
caract{\'e}ris{\'e}s par : 

$$\int_{\overline{X_{w'}}}\hat{\xi}^{w}=\delta_{w',w}.$$
   
\end{enumerate}

\end{prop}

On pose  $\xi^{w}=i_{T}^*(\hat{\xi}^{w})$. Soit $(w,v) \in W^2$ et soit 
$v=s_{i_{1}}\cdots s_{i_{l}}$ une d{\'e}composition non n{\'e}cessairement r{\'e}duite de $v$, 
pour $1 \leq j \leq l$, on d{\'e}finit un {\'e}l{\'e}ment $\beta_{j}$ de 
$S(\mathfrak{h}^*)$ par $\beta_{j}=s_{i_{1}}\cdots 
s_{i_{j-1}}\alpha_{i_{j}}$. La formule suivante est prouv{\'e}e par
Sara Billey dans \cite{sb} : 

\begin{theo} \label{cohomologieVD}

    Soient $w$ et $v$ deux {\'e}l{\'e}ments de $ W$ tels que $w\leq v$ et $m=l(w)$, on a alors : 
   $$ \xi^{w}(v) = \sum\beta_{j_{1}}\cdots \beta_{j_{m}},$$
    o{\`u} la somme porte sur l'ensemble des entiers $1\leq j_{1} 
   < \cdots < j_{m} \leq l$ tels que $w=s_{i_{j_{1}}}\cdots 
   s_{i_{j_{m}}}$.
    
\end{theo}

\begin{exem}
Pla{\c c}ons nous dans le cas $A_4$ et calculons $\xi^{w}(v)$
avec $w=s_{3}s_{2}$ et $v=s_{2}s_{3}s_{2}s_{1}s_{2}$. Il y a $2$ fa{\c c}ons de
``trouver $w$ en dessous de $v$'' :
$w=s_{i_{2}}s_{i_{3}}$, 
$w=s_{i_{2}} s_{i_{5}}$, 
 et on trouve donc :
$$\xi^{w}(v)=(\alpha_{2}+\alpha_{3})\alpha_{3}+
(\alpha_{2}+\alpha_{3})(\alpha_{1}+\alpha_{2})=\alpha_1\alpha_2 +
\alpha_1\alpha_3+\alpha_2^2
+2\alpha_2\alpha_3 +\alpha_3^2 . $$

\end{exem}

Retrouvons le th{\'e}or{\`e}me~\ref{cohomologieVD} gr{\^a}ce aux r{\'e}sultats pr{\'e}c{\'e}dents.

Soit $v=s_{\mu_{1}} \cdots
s_{\mu_{N}}$ une d{\'e}composition non n{\'e}cessairement r{\'e}duite d'un {\'e}l{\'e}ment
$v$ de $W$. On pose $\Gamma=
\Gamma(\mu_{1}, \ldots, \mu_{N})$ et $g=g_{\mu_{1}, \ldots, \mu_{N}}$. La
proposition suivante va nous permettre de retrouver le th{\'e}or{\`e}me \ref{cohomologieVD} :

\begin{prop} \label{decompositionVBS}

Soit $w \in W$, on a : 

$$g^*(\hat{\xi}^w)=\!\!\!\!\!\!\!\!
\sum_{ \tiny \begin{array}{cc}\epsilon \in \mathcal{E}, l(\epsilon)=l(w) \\ 
 { \rm et } \,  v(\epsilon)=w
\end{array}}\! \! \! \! \! \! \! \! \!  \! \! \! \hat{\sigma}_{\epsilon}^T.$$

\end{prop}

\begin{proof}

Pour d{\'e}montrer la formule, il faut montrer : 

$$\int_{\overline{\Gamma_{\epsilon}}}g^*(\hat{\xi}^w)=\delta_{v(\epsilon),w}
\delta_{l(\epsilon),l(w)}.$$

Distinguons deux cas :  

\smallskip

Si $\epsilon$ correspond {\`a} une d{\'e}composition r{\'e}duite de
$v(\epsilon)$, comme $g_{\big| \overline{\Gamma_{\epsilon}}} :  \overline{\Gamma_{\epsilon}}
\rightarrow X$ s'identifie {\`a} l'application $g_{\mu_{i} , i\in
  \pi_{+}(\epsilon)}$, d'apr{\`e}s la proposition~\ref{desingularisation} : 
$$\int_{\overline{\Gamma_{\epsilon}}}g^*(\hat{\xi}^w)=
\int_{\overline{X_{v(\epsilon)}}}\hat{\xi}^w=\delta_{v(\epsilon) , w}.$$

\smallskip

Si $\epsilon$ n'est pas une d{\'e}composition r{\'e}duite de $v(\epsilon)$, 
alors, d'apr{\`e}s le lemme~\ref{image}, $g$ envoie $\overline{\Gamma_{\epsilon}}$
dans
$\overline{X_{\underline{v}(\epsilon)}}$ qui est de
dimension strictement plus petite que $\overline{\Gamma_{\epsilon}}$, et donc
$\int_{\overline{\Gamma_{\epsilon}}}g^*(\hat{\xi})=0$ pour tout {\'e}l{\'e}ment 
$ \hat{\xi}$ de $H_T^*(X)$, et en particulier :
$$\int_{\overline{\Gamma_{\epsilon}}}g^*(\hat{\xi}^w)=0.$$


\end{proof}

\begin{proof}[D{\'e}monstration du th{\'e}or{\`e}me~\ref{cohomologieVD}]

Si on note $\tilde{u} : F_b(W;S(\mathfrak{h}^*))  \rightarrow 
 F_b(\mathcal{E};S(\mathfrak{h}^*))$ l'application induite par $\overline{u} : 
 \mathcal{E}\rightarrow W$ d{\'e}finie par
$\overline{u}(\epsilon)=v(\epsilon)$, on a le diagramme commutatif suivant : 
$$\xymatrix{
     H_{T}^*(\Gamma)\ar[dd]^{i_{T}^*} & &
     \ar[ll]^{g^{*}}H_{T}^*(X)\ar[dd]^{i_{T}^*} \\
     \\
     F_b(\mathcal{E};S(\frak{h}^*)) & & 
   \ar[ll]^{\tilde{u}}F_b(W;S(\frak{h}^*)) }$$

 Ce diagramme nous donne en particulier pour tout {\'e}l{\'e}ment $w \in W$ : 
 
$$\xi^w(v)=(\tilde{u}i_T^*(\xi^w))((\mathbf{1}))=i_T^*g^*(\hat{\xi}^w)((\mathbf{1})),$$
et donc gr{\^a}ce {\`a} la proposition~\ref{decompositionVBS} : 

$$\xi^w(v)=\!\!\!\!\!\!\!\!\!\!\sum_{ \tiny \begin{array}{ll}\epsilon' \in \mathcal{E}, 
l(\epsilon')=l(w) \\ 
\hspace{0,25 cm} {\rm et } \hspace{0,15 cm}  v(\epsilon')=w
\end{array}}\!\!\!\!\!\!\!\!\!\!\!\!\!\!
\sigma_{\epsilon'}^T((\mathbf{1}\rm)),$$
ce qui nous redonne bien le th{\'e}or{\`e}me~\ref{cohomologieVD} {\`a} l'aide de 
l'expression de $\sigma_{\epsilon'}^T((\mathbf{1}))$ 
(d{\'e}finition~\ref{defcohomologieBS}), en remarquant que pour tout entier $i$
compris entre $1$ et $N$, $\alpha_{i}((\bf{1}\rm))=-\beta_{i}$.

\end{proof}

\begin{exem}

On consid{\`e}re le cas $A_2$, $\Gamma=\Gamma(\alpha_1,\alpha_2,
 \alpha_1)$ la vari{\'e}t{\'e} de Bott-Samelson associ{\'e}e
 {\`a} la d{\'e}composition $w_0=s_1s_2s_1$ du plus grand {\'e}l{\'e}ment du groupe de
 Weyl $W$. On reprend les notations du
th{\'e}or{\`e}me~\ref{carreBS}, et on a alors : 
\begin{eqnarray*} 
g^*(\hat{\xi}^{s_1}) & = & x_1 + x_3 \\
g^*(\hat{\xi}^{s_2}) & = & x_2 \\
g^*(\hat{\xi}^{s_1s_2}) & = & x_1x_2 \\
g^*(\hat{\xi}^{s_2s_1}) & = & x_2x_3 \\
g^*(\hat{\xi}^{s_1s_2s_1}) & = & x_1x_2x_3 
\end{eqnarray*}

On peut v{\'e}rifier alors facilement, gr{\^a}ce aux relations explicit{\'e}es dans
l'exemple~\ref{exempleA2BS}, que $g^*(H_T^*(X))$ est bien une sous-alg{\`e}bre de $H_T^*
(\Gamma)$.

\end{exem}

\begin{exem}

On consid{\`e}re le cas $B_2$, $\Gamma=\Gamma(\alpha_1,\alpha_2,
 \alpha_1,\alpha_2)$  la vari{\'e}t{\'e} de Bott-Samelson associ{\'e}e
 {\`a} la d{\'e}composition $w_0=s_1s_2s_1s_2$ du plus grand {\'e}l{\'e}ment du groupe de
 Weyl $W$. On reprend les notations du
th{\'e}or{\`e}me~\ref{carreBS}, et on a alors : 
\begin{eqnarray*} 
g^*(\hat{\xi}^{s_1}) & = & x_1 + x_3 \\
g^*(\hat{\xi}^{s_2}) & = & x_2 +x_4\\
g^*(\hat{\xi}^{s_1s_2}) & = & x_1x_2+x_1x_4+x_3x_4 \\
g^*(\hat{\xi}^{s_2s_1}) & = & x_2x_3 \\
g^*(\hat{\xi}^{s_1s_2s_1}) & = & x_1x_2x_3 \\
g^*(\hat{\xi}^{s_2s_1s_2}) & = & x_2x_3x_4 \\
g^*(\hat{\xi}^{s_1s_2s_1s_2}) & = & x_1x_2x_3x_4
\end{eqnarray*}

Comme pr{\'e}c{\'e}demment, on peut v{\'e}rifier,  gr{\^a}ce aux relations explicit{\'e}es dans
l'exemple~\ref{exempleB2BS}, que $g^*(H_T^*(X))$ est bien une sous-alg{\`e}bre de $H_T^*
(\Gamma)$.

\medskip

Illustrons le th{\'e}or{\`e}me~\ref{cohomologieVD} sur cet exemple, et calculons
$ \xi^{s_1}(s_1s_2s_1s_2)$ : 
$$\xi^{s_1}(s_1s_2s_1s_2)=\sigma^T_1((\mathbf{1}))+\sigma^T_3((\mathbf{1}))=
 \alpha_1+ s_1s_2  \alpha_1   =
\alpha_1+(\alpha_1+\alpha_2)=2\alpha_1+\alpha_2.$$

\bigskip

 On peut utiliser ce plongement de
$H_T^*(X)$ dans $H_T^*(\Gamma)$ pour calculer des produits. En effet, si on 
veut calculer  par exemple $\hat{\xi}^{s_1s_2}\hat{\xi}^{s_2s_1}$, on a :
 $$g^*(\hat{\xi}^{s_1s_2}\hat{\xi}^{s_2s_1})=(x_1x_2+x_1x_4+x_3x_4)x_2x_3=$$
$$x_1x_3(\alpha_2x_2+2x_1x_2)+x_1x_2x_3x_4+x_2(\alpha_1x_3+x_2x_3-2x_1x_3)x_4=$$
$$(2\alpha_1+\alpha_2)x_1x_2x_3
+x_1x_2x_3x_4+\alpha_1x_2x_3x_4+(\alpha_2x_2+2x_1x_2)x_3x_4-
2x_1x_2x_3x_4=$$
$$(2\alpha_1+\alpha_2)x_1x_2x_3+(\alpha_1+\alpha_2)x_2x_3x_4+x_1x_2x_3x_4,$$
et donc $\hat{\xi}^{s_1s_2}\hat{\xi}^{s_2s_1}=
(2\alpha_1+\alpha_2)\hat{\xi}^{s_1s_2s_1}+(\alpha_1+\alpha_2)\hat{\xi}^{s_2s_1s_2}
+\hat{\xi}^{s_1s_2s_1s_2}$.
\medskip

On g{\'e}n{\'e}ralisera cette m{\'e}thode dans le chapitre~\ref{schubertequivariant}.

\end{exem}

\chapter{K-th{\'e}orie {\'e}quivariante} \label{ktheorie}

\section{Pr{\'e}liminaires}

\subsection{D{\'e}finitions}

 Soit $U_{\mathbb{C}}$ un tore complexe 
 d'alg{\`e}bre de Lie $\mathfrak{u}_{\mathbb{C}}$,
et soit $U   \subset U_{\mathbb{C}}$ le tore compact maximal de
 $U_{\mathbb{C}}$.
  On note 
 $\mathfrak{u} \subset \mathfrak{u}_{\mathbb{C}} $ l'alg{\`e}bre de Lie de $U$.
 On note $X[U]$ 
le groupe des caract{\`e}res de $U$, et on pose $R[U]=\mathbb{Z}[X[U]]$. On note
 $Q[U]$ le corps des fractions de $R[U]$.

Pour tout poids entier $\alpha \in i\mathfrak{u}^* \subset 
\mathfrak{u}_{\mathbb{C}}^*$, on note $e^{\alpha} : U \rightarrow S^1$ le
caract{\`e}re correspondant.

\bigskip

Soit $Z$ un espace topologique compact muni d'une action continue de $U$. 
On d{\'e}finit la $K$-th{\'e}orie $U$-{\'e}quivariante de $Z$ comme le groupe construit {\`a}
partir du semi-groupe des classes
d'isomorphisme de fibr{\'e}s vectoriels complexes de dimension finie
$U$-{\'e}quivariants au 
dessus de $Z$. On munit ce groupe
d'une structure d'anneau d{\'e}finie {\`a} l'aide du produit tensoriel. De
plus, comme la  $K$-th{\'e}orie $U$-{\'e}quivariante du point s'identifie {\`a} $R[U]$, on
obtient une structure de 
$R[U]$-alg{\`e}bre qu'on notera $K_{U}(Z)$.

\medskip

Toute application $g : Z_1 \rightarrow Z_2$ continue et $U$-{\'e}quivariante d{\'e}finit
un morphisme de $R[U]$-alg{\`e}bre $g^* : K_{U}(Z_2) \rightarrow K_{U}(Z_1)$. En
particulier, l'inclusion $Z^U \subset Z$ d{\'e}finit un morphisme 
$i_U^*  : K_{U}(Z) \rightarrow K_{U}(Z^U)$ appel{\'e} restriction aux points fixes. 
 Si l'ensemble des points fixes $Z^U$ est discret, $  K_{U}(Z^U) $ 
s'identifie de mani{\`e}re {\'e}vidente {\`a} $F(Z^U; R[U])$ la
$R[U]$-alg{\`e}bre
des fonctions de $Z^U$ {\`a} valeurs dans $R[U]$ munie de l'addition et de la
multiplication point par point. On obtient alors un morphisme $i_U^* : 
  K_{U}(Z)  \rightarrow  F(Z^U; R[U]) $.

\subsection{Formule de localisation}

On suppose que $Z$ est une vari{\'e}t{\'e} complexe projective lisse
de dimension 
$n$ munie d'une action de 
$U_{\mathbb{C}}$. Cette action induit alors une action de $U$ sur $Z$. 

La vari{\'e}t{\'e} $Z$ {\'e}tant lisse, le groupe 
$K_{0}(U_{\mathbb{C}},Z)$ construit {\`a} partir du semi-groupe des classes 
d'isomorphisme de faisceaux 
$U_{\mathbb{C}}$-{\'e}quivariants coh{\'e}rents sur $Z$ est isomorphe au groupe 
$K^{0}(U_{\mathbb{C}}, Z)$ construit {\`a} partir du semi-groupe des classes 
d'isomorphisme de faisceaux $U_{\mathbb{C}}$-{\'e}quivariants localement libres sur
$Z$ (voir \cite{ginzburg}, chapitre $5$). On identifie donc ces deux groupes.

On a un morphisme canonique : $K^{0}(U_{\mathbb{C}}, Z) \rightarrow K_U(Z)$. On
suppose que ce morhisme est un isomorphisme (c'est le cas pour les tours de Bott
et les vari{\'e}t{\'e}s de drapeaux dans le cas fini), et on identifie ces deux groupes.

 Pour toute sous-vari{\'e}t{\'e} $U_{\mathbb{C}}$-invariante  $Z'$ et tout 
 faisceau $\mathcal{F} \in  K^{0}(U_{\mathbb{C}}, Z) $, l'action de
 $U_{\mathbb{C}}$ sur $Z$ induit une action de $U_{\mathbb{C}}$ sur 
${\rm H}^k(Z',\mathcal{F}_{/Z'})$, et on d{\'e}finit 
$\chi(Z', \mathcal{F}) \in R[U]$ par :  
$$\forall u \in U, \hspace{0,1 cm} \chi(Z',
\mathcal{F})(u)=\sum_{k}(-1)^k{ \rm Tr }
 (u;  { \rm H}^k(Z',\mathcal{F}_{/Z'})).$$

On suppose de plus que $Z^U$ est fini. 
En chaque point fixe $m \in Z^U$, on note
$(\alpha_1^m, \ldots ,\alpha_n^m) \in (i\mathfrak{u}^*)^n 
 \subset (\mathfrak{u}_{\mathbb{C}}^*)^n$ les poids de la
repr{\'e}sentation de $U$ dans $T_mZ$, l'espace tangent {\`a} $Z$ en $m$. 
Dans ce cas, la formule $5.11.9$ de \cite{ginzburg} 
s'{\'e}crit de la mani{\`e}re suivante : 

\begin{prop} \label{pointsfixesab}
 
Pour tout faisceau $\mathcal{F}$ localement libre et
$U_{\mathbb{C}}$-{\'e}quivariant au dessus de $Z$, $\chi(Z,\mathcal{F})  $ se
calcule gr{\^a}ce {\`a} la formule suivante  : 

$$\chi(Z,\mathcal{F})=\sum_{m \in Z^U}\frac{i_{U}^*(\mathcal{F})(m) }
{\prod_{1 \leq i \leq n}(1-e^{-\alpha_i^m}) }.$$

\end{prop}

\section{K-th{\'e}orie {\'e}quivariante des tours de Bott}

On reprend les notations de la section~\ref{tours}. Soit $N\geq 1$ un entier
naturel, et soit $C=\{c_{i,j}\}_{1 \leq i<j \leq N}$ une liste d'entiers. On pose
$Y=Y_{C}$. 

Pour un poids entier $\lambda \in \oplus_{1 \leq k \leq N}\mathbb{Z}\lambda_k
\subset i\mathfrak{d}^*$, on note $e^{\lambda} : D \rightarrow S^1$ le
caract{\`e}re correspondant. 

On montre par r{\'e}currence sur la dimension de $Y$ (voir \cite{kkk} o{\`u} le 
r{\'e}sultat est d{\'e}montr{\'e} dans le cas particulier des vari{\'e}t{\'e}s de 
Bott-Samelson) que le morphisme canonique : $K^{0}(D_{\mathbb{C}}, Y)
\rightarrow K_D(Y)$ est un isomorphisme. 
Dans la suite, on identifie donc ces deux groupes.

\medskip

Tout comme dans le cas de la cohomologie {\'e}quivariante, on a la structure
suivante de la K-th{\'e}orie $D$-{\'e}quivariante de $Y$ :  

\begin{prop} \label{propbasektheorie}

\indent

\begin{enumerate}

\item[$(i)$] La  $K$-th{\'e}orie $D$-{\'e}quivariante de $Y^D$ s'identifie {\`a}
$F(\mathcal{E};R[D])$.

\item[$(ii)$] La restriction aux points fixes $i_{D}^*$ : $K_{D}(Y)
 \rightarrow F(\mathcal{E};R[D])$ est injective.

\item[$(iii)$] La $K$-th{\'e}orie $D$-{\'e}quivariante de $Y$ est un $R[D]$-module
  libre de rang $2^N$.
\end{enumerate}

\end{prop}

\begin{proof}

Le point $(i)$ est imm{\'e}diat.

\smallskip

Le point $(ii)$ est une cons{\'e}quence de $(iii)$. En effet, d'apr{\`e}s le th{\'e}or{\`e}me de
localisation (voir \cite{kequivariante}), le morphisme : 
$K_{D}(Y)\otimes_{R[D]}Q[D] \rightarrow F(\mathcal{E}; Q[D])$ induit par
$i_D^*$ est un isomorphisme. De plus, comme $K_{D}(Y)$ est un $R[D]$-module
libre, $K_{D}(Y)$ s'injecte dans $  K_{D}(Y)\otimes_{R[D]}Q[D] $, et on a donc
le diagramme commutatif suivant qui prouve que $i_{D}^*$ est injective : 
$$\xymatrix{
     K_{D}(Y)\ar[dd]^{i_D^*} \ar@{^{(}->}[rr] & &
    K_{D}(Y)\otimes_{R[D]}Q[D]\ar[dd]^{\simeq} \\
     \\
     F(\mathcal{E};R[D]) \ar[rr] & & 
   F(\mathcal{E};Q[D]) }
  \\
$$

\bigskip

Pour d{\'e}montrer $(iii)$, on va expliciter une base 
$\{ \hat{\mu}_{\epsilon}^D \}_{\epsilon \in \mathcal{E}}$
du $R[D]$-module $K_D(Y)$. On proc{\`e}de par r{\'e}currence sur $N \geq 1$.

\smallskip

 Pour $N=1$,
$K_{S^1}(\mathbb{P}^1)$ est un $R[S^1]$-module libre engendr{\'e} par le fibr{\'e} en
droites 
trivial $\mathbf{1}$ et par le fibr{\'e} $\mathbf{E}$
d{\'e}fini comme le fibr{\'e} en droites tautologique sur  $\mathbb{P}^1$ (voir
\cite{atiyah}, corollaire 2.2.2). On pose $ \hat{\mu}_{(0)}^D= \mathbf{E}$ et 
$ \hat{\mu}_{(1)}^D= \mathbf{1} - \mathbf{E}$.

\smallskip

On suppose le r{\'e}sultat v{\'e}rifi{\'e} pour toute tour de Bott de dimension $N-1$. Soit
$Y=Y_{C}$ et $Y'=Y_{C_{N-1}}$. Alors $Y=\mathbb{P}(\mathbf{1}\oplus \mathbf{L}_N)$, o{\`u} 
 $\mathbf{L}_N$ est un fibr{\'e} en droites au dessus de $Y'$ 
(voir la section~\ref{def} pour la d{\'e}finition de
$Y_{C_{N-1}}$ et  $\mathbf{L}_N$). On note 
$\pi_N : Y=\mathbb{P}(\mathbf{1}\oplus \mathbf{L}_N)
\rightarrow Y'$ la projection de $Y$ sur $Y'$ et $\mathbf{E} \in
K_D(Y)$ le fibr{\'e} tautologique au dessus de 
$Y=\mathbb{P}(\mathbf{1}\oplus \mathbf{L}_N)$.

\smallskip

Soit $D'=(S^1)^{N-1}$, par hypoth{\`e}se de r{\'e}currence, $K_{D'}(Y')$ est un 
$R[D']$-module libre engendr{\'e} par une base 
$\{ \hat{\mu}_{f}^{D'} \}_{f \in \{0,1\}^{N-1}}$. 
On d{\'e}finit une action de
 $D=(S^1)^N$
sur $Y'$ en faisant agir la derni{\`e}re composante trivialement sur $Y'$. On a alors 
$K_{D}(Y')=K_{D'}(Y')\otimes_{\mathbb{Z}}R[S^1]$ (voir \cite{kequivariante}),
et $K_{D}(Y')$ est donc
un $R[D]$-module libre qui admet comme base la famille
$\{ \hat{\mu}_{f}^D \}_{f \in \{0,1\}^{N-1}}$, o{\`u} pour tout $f \in
\{0,1\}^{N-1}$, on a pos{\'e} $ \hat{\mu}_{f}^D=\hat{\mu}_{f}^{D'}\otimes 1$.

\smallskip

 La projection $\pi_N$ munit $K_D(Y)$ d'une structure de $K_D(Y')$-module. Comme
$Y=\mathbb{P}(\mathbf{1}\oplus \mathbf{L}_N)$, 
 $K_{D}(Y)$ est un $K_{D}(Y')$-module 
libre de rang $2$ engendr{\'e} par le fibr{\'e} trivial $\mathbf{1}$ et le fibr{\'e} tautologique  
$\mathbf{E}$ (voir \cite{atiyah} th{\'e}or{\`e}me 2.2.1, ou \cite{kequivariante}). 
Soit $p : \mathcal{E}=\{0,1\}^N \rightarrow
\{0,1\}^{N-1}$ la projection selon les $N-1$ premi{\`e}res coordonn{\'e}es, 
on obtient donc une base du  
$R[D]$-module $K_D(Y)$ en posant :

$$\left\{ \begin{array}{ll} \hat{\mu}_{\epsilon}^D = 
\pi_N^*(\hat{\mu}_{p(\epsilon)}^D)\mathbf{E}
 &\,\,\, { \rm si } \,\,\, \epsilon_N =0, \\ 
  \hat{\mu}_{\epsilon}^D = 
\pi_N^*(\hat{\mu}_{p(\epsilon)}^D)(\mathbf{1}-\mathbf{E})
 &\,\,\, { \rm si } \,\,\, \epsilon_N =1 .
\end{array}\right.$$
 
\end{proof}

On va expliciter la base $\{ \hat{\mu}_{\epsilon}^D \}_{\epsilon \in
  \mathcal{E}}$ d{\'e}finie dans la d{\'e}monstration de la proposition pr{\'e}c{\'e}dente. On
  reprend le diagramme de la section~\ref{def} : 
$$ \begin{array}{clccccc}
 & &  & & & \mathbb{P}(\mathbf{1}\oplus \mathbf{L}_N)=&Y_C \\ 
 &  &  & &  &  \downarrow \pi_N &  \\
  &  & & & &  Y_{C_{N-1}}&  \\
 & &   & &\Ddots &  \\
  & & \mathbb{P}(\mathbf{1}\oplus  \mathbf{L}_2)=& Y_{C_2}&   & &  \\
  & & \downarrow \pi_2 & &  & &  \\
  & \mathbb{C}P^1 \,\,\, = & \!\! Y_{C_1} & & & &  \\
  &\,\, \downarrow \pi_1  & & &  & &  \\
 \{un \, point\}= &\,\,  Y_0 & & &  & &  
\end{array}$$

\medskip

Pour $1 \leq i \leq N-1$, chaque vari{\'e}t{\'e} $Y_{C_i}$ {\'e}tant munie d'une action de 
$ (S^1)^i $, on d{\'e}finit une action de  $D=(S^1)^N$ sur $Y_{C_{i}}$ en faisant agir
trivialement les derni{\`e}res composantes de $D$.
Pour $1 \leq i \leq N$, on pose $\Pi_i=\pi_{i+1} \pi_{i+2} \cdots \pi_N  : Y 
\rightarrow Y_{C_{i}}$ ($\Pi_N = { \rm Id }_{Y}$). 
De plus, on note $\mathbf{E_i} \in K_{D}(Y_{C_i})
\simeq K_{(S^1)^i}(Y_{C_i} )\otimes_{\mathbb{Z}}R[(S^1)^{N-i}]$ le fibr{\'e}
tautologique sur $Y_{C_{i}}$, et 
on pose $\mathbf{F_i}=\mathbf{1}-\mathbf{E_i} $.

Pour tout $\epsilon \in \mathcal{E}$, on v{\'e}rifie alors que 
$\hat{\mu}_{\epsilon}^D \in K_D(Y)$ est donn{\'e} par la formule :

$$  \hat{\mu}_{\epsilon}^D =\prod_{i\in \pi_+(\epsilon)}\Pi_i^*(\mathbf{F_i})
\prod_{i\in \pi_-(\epsilon)}\Pi_i^*(\mathbf{E_i}).$$

\medskip

Le th{\'e}or{\`e}me suivant donne la valeur des restrictions aux points fixes des
classes $\hat{\mu}_{\epsilon}^D $.
Si on pose $ \mu_{\epsilon}^D = i_{D}^*( \hat{\mu}_{\epsilon}^D)$, on a la
formule suivante : 

\begin{theo}  \label{restrictionkb}
    
Pour $(\epsilon,\epsilon') \in \mathcal{E}^2$ :

$$\left\{ \begin{array}{ll} \mu_{\epsilon}^D(\epsilon') = 
\displaystyle{\prod_{i \in \pi_{+}(\epsilon')}e^{-\lambda_{i}(\epsilon')}
\prod_{i \in \pi_{+}(\epsilon)}(e^{\lambda_{i}(\epsilon')}-1) }
 & { \rm si } \hspace{0,15 cm} \epsilon \leq \epsilon', \\ 
 \mu_{\epsilon}^D(\epsilon') = 0
 & { \rm  sinon }.
\end{array}\right.$$

\end{theo}

\begin{proof}

Pour d{\'e}montrer ce th{\'e}or{\`e}me, il suffit de calculer les restrictions aux
points fixes des fibr{\'e}s $\Pi_i^*(\mathbf{E_i})$. Pour $1 \leq i \leq N$, on pose
$p_i : \mathcal{E} =\{0,1 \}^N \rightarrow \{0,1 \}^{i}$ la projection selon les
$i$ premi{\`e}res coordonn{\'e}es. On a alors pour tout $\epsilon' \in \mathcal{E}$, 
$i_D^*( \Pi_i^*(\mathbf{E_i}))(\epsilon')=i_D^*(\mathbf{E_i})(p_i(\epsilon'))$. 

Or, pour $e \in \{0,1\}^i$, $i_D^*(\mathbf{E_i})(e)=1$ si $e_i=0$, et
$i_D^*(\mathbf{E_i})(e)=e^{-\lambda_i(\epsilon')}$ si $e_i=1$. On obtient donc
pour tout $\epsilon' \in \mathcal{E}$ : 

$$\left\{ \begin{array}{ll} i_D^*( \Pi_i^*(\mathbf{E_i}))(\epsilon')=1

 & { \rm si } \hspace{0,15 cm} \epsilon'_i = 0,  \\ 
i_D^*( \Pi_i^*(\mathbf{E_i}))(\epsilon')=e^{-\lambda_i(\epsilon')}
 & { \rm si } \hspace{0,15 cm} \epsilon'_i = 1,
\end{array}\right.$$

et : 

$$\left\{ \begin{array}{ll} i_D^*( \Pi_i^*(\mathbf{F_i}))(\epsilon')=0

 & { \rm si } \hspace{0,15 cm} \epsilon'_i = 0,  \\ 
i_D^*( \Pi_i^*(\mathbf{F_i}))(\epsilon')=1-e^{-\lambda_i(\epsilon')}
 & { \rm si } \hspace{0,15 cm} \epsilon'_i = 1 .
\end{array}\right.$$

\medskip

Comme pour tout $\epsilon \in \mathcal{E}$, 
$\hat{\mu}_{\epsilon}^D =\prod_{i\in \pi_+(\epsilon)}\Pi_i^*(\mathbf{F_i})
\prod_{i\in \pi_-(\epsilon)}\Pi_i^*(\mathbf{E_i})$, on obtient bien
 $\mu_{\epsilon}^D(\epsilon')=0$ si $\epsilon' \not\geq
\epsilon$, et si  $\epsilon' \geq \epsilon$ : 

$$\mu_{\epsilon}^D(\epsilon')=\prod_{i\in
  \pi_+(\epsilon)}(1-e^{-\lambda_i(\epsilon')})  \prod_{i\in
  \pi_-(\epsilon) \cap \pi_+(\epsilon')}e^{-\lambda_i(\epsilon')}=
\prod_{i \in \pi_{+}(\epsilon')}e^{-\lambda_{i}(\epsilon')}
\prod_{i \in \pi_{+}(\epsilon)}(e^{\lambda_{i}(\epsilon')}-1).
$$

\end{proof}

\begin{exem}

On consid{\`e}re la surface de Hirzebruch $H_{-1}=Y_{\{-1\}}$. On pose : 

 $$\left\{ \begin{array}{lll} \epsilon_{1}=(0,0), &   \\ 
 \epsilon_{2}=(1,0), & \epsilon_{3}=(0,1) ,
 \\ \epsilon_{4}=(1,1). & 
 & \end{array} \right.$$

Si on d{\'e}finit la matrice $\mathcal{M}=\{ \mu_{i,j}  \}_{1 \leq i< j \leq 4}$ par 
$\mu_{i,j}=\mu_{\epsilon_i}^D(\epsilon_j)$, alors : 

$$\mathcal{M}=\begin{pmatrix} 1 & e^{-\lambda_1} & e^{-\lambda_2} &
  e^{-2\lambda_1-\lambda_2} \\
 0 &1- e^{-\lambda_1} & 0 &
  e^{-\lambda_1-\lambda_2}(1- e^{-\lambda_1}) \\
 0 & 0 &1- e^{-\lambda_2} &
  e^{-\lambda_1}(1- e^{-\lambda_1-\lambda_2}) \\
0 & 0 &0  &
  (1- e^{-\lambda_1})(1-e^{-\lambda_1-\lambda_2}) 
\end{pmatrix}.$$

\end{exem}

\bigskip

La base $\{  \hat{\mu}_{\epsilon}^D \} _{\epsilon \in \mathcal{E}}$ est reli{\'e}e {\`a}
la d{\'e}composition cellulaire $Y=\coprod_{\epsilon \in \mathcal{E}}Y_{\epsilon}$ par le
th{\'e}or{\`e}me suivant :

\begin{theo} \label{ktheorietours}
    
La famille $\{  \hat{\mu}_{\epsilon}^D \} _{\epsilon \in \mathcal{E}}$ est une
base du R[D]-module $K_{D}(Y)$ caract{\'e}ris{\'e}e par : 

$$\forall
(\epsilon, \epsilon') \in \mathcal{E}^2, \,\, \chi(\overline{Y_{\epsilon'}}, 
\hat{\mu}_{\epsilon}^D)= \delta_{\epsilon, \epsilon'} . $$
   
\end{theo}

\begin{proof}

On sait d{\'e}j{\`a} que la famille  $\{  \hat{\mu}_{\epsilon}^D \} _{\epsilon \in
  \mathcal{E}}$ est une base de $K_{D}(Y)$. Pour $(\epsilon, \epsilon') \in
  \mathcal{E}^2$, on va calculer  
$\chi(\overline{Y_{\epsilon'}}, \hat{\mu}_{\epsilon}^D)$ gr{\^a}ce {\`a} la formule de
  localisation.

En utilisant la proposition~\ref{pointsfixesab}, 
 et le lemme~\ref{pointsfixestours}, on
obtient pour tout  $\hat{\mu}^D \in K_{D}(Y)$ et tout $\epsilon \in \mathcal{E}$: 

\begin{equation} \label{ab} \chi(\overline{Y_{\epsilon}}, \hat{\mu}^D)=  
\sum_{\epsilon' \leq
\epsilon}\frac{i_{D}^*(\hat{\mu}^D)(\epsilon')}{\prod_{i \in \pi_{+}(\epsilon)}
(1-e^{-\lambda_{i}(\epsilon')})}.  \end{equation}

Cette formule et le th{\'e}or{\`e}me~\ref{restrictionkb} nous montrent imm{\'e}diatement
que $\chi(\overline{Y_{\epsilon}}, \hat{\mu}_{\epsilon}^D)=1$ et 
$\chi(\overline{Y_{\epsilon}}, \hat{\mu}_{\epsilon_0}^D)=0$ si $\epsilon_0 \not\leq
\epsilon$.

Soit $\epsilon_0 \in \mathcal{E}$ tel que $\epsilon_0 \leq \epsilon$ et $\epsilon_0 \neq
\epsilon$. Alors, la formule~\ref{ab} et le th{\'e}or{\`e}me~\ref{restrictionkb} nous
donnent :

$$ \chi(\overline{Y_{\epsilon}}, \hat{\mu}_{\epsilon_0}^D) = \sum_{\epsilon_{0} \leq
  \epsilon' \leq
\epsilon}\frac{\displaystyle{\prod_{i \in \pi_{+}(\epsilon')}e^{-\lambda_{i}(\epsilon')}
\prod_{i \in \pi_{+}(\epsilon_{0})}(e^{\lambda_{i}(\epsilon')}-1)}}
{\displaystyle{\prod_{i \in
  \pi_{+}(\epsilon)}(1-e^{-\lambda_{i}(\epsilon')})}},$$

d'o{\`u} : 

$$\chi(\overline{Y_{\epsilon}}, \hat{\mu}_{\epsilon_0}^D) = 
\sum_{\epsilon_{0} \leq  \epsilon' \leq
\epsilon}\frac{\displaystyle{\prod_{i \in \pi_{+}(\epsilon') \setminus \pi_{+}(\epsilon_{0})}
e^{-\lambda_{i}(\epsilon')}} }
{\displaystyle{\prod_{i \in \pi_{+}(\epsilon)\setminus \pi_{+}(\epsilon_{0})}
(1-e^{-\lambda_{i}(\epsilon')})}}.$$

\medskip

Soit $j$ le plus grand
{\'e}l{\'e}ment de $\pi_{+}(\epsilon) \setminus  \pi_{+}(\epsilon_{0})$, on a alors :

$$\mbox{ \small { \mbox{ $ \chi(\overline{Y_{\epsilon}},
      \hat{\mu}_{\epsilon_0}^D)=$}}}\!\!\!\!\!\sum_{
 \tiny \begin{array}{cc}  \epsilon_{0} \leq \epsilon' \leq
\epsilon \\
 \epsilon'_j=0 
\end{array}}\!\!\!\!\!\!\!\!\mbox{ \small { \mbox{ $ 
\frac{ \displaystyle{ \prod_{i \in \pi_{+}(\epsilon') \setminus
  \pi_{+}(\epsilon_{0})} e^{-\lambda_{i}(\epsilon') }}}
{\displaystyle{\prod_{i \in \pi_{+}(\epsilon)\setminus \pi_{+}(\epsilon_{0})}
(1-e^{-\lambda_{i}(\epsilon')})}}$}}} +
\!\!\!\!\!\sum_{
 \tiny \begin{array}{cc}  \epsilon_{0} \leq \epsilon' \leq
\epsilon \\
 \epsilon'_j=1 
\end{array}}\!\!\!\!\!\!\!\!\mbox{ \small { \mbox{ $ 
\frac{ \displaystyle{ \prod_{i \in \pi_{+}(\epsilon') \setminus
  \pi_{+}(\epsilon_{0})} e^{-\lambda_{i}(\epsilon') }}}{\displaystyle{\prod_{i 
\in \pi_{+}(\epsilon)\setminus \pi_{+}(\epsilon_{0})}
(1-e^{-\lambda_{i}(\epsilon')})}}$}}}, $$
d'o{\`u} :

$$\mbox{ \small { \mbox{ $ \chi(\overline{Y_{\epsilon}},
      \hat{\mu}_{\epsilon_0}^D)=$}}}\!\!\!\!\!\sum_{
 \tiny \begin{array}{cc}  \epsilon_{0} \leq \epsilon' \leq
\epsilon \\
 \epsilon'_j=0 
\end{array}}\!\!\!\!\!\!\!\! \mbox{ \small { \mbox{ $ \left[ \frac{\displaystyle{ 
\prod_{i \in \pi_{+}(\epsilon') \setminus
  \pi_{+}(\epsilon_{0})} e^{-\lambda_{i}(\epsilon') }}}{
\displaystyle{\prod_{i \in \pi_{+}(\epsilon)\setminus \pi_{+}(\epsilon_{0})}
(1-e^{-\lambda_{i}(\epsilon')})}}+\frac{ \displaystyle{ \prod_{i \in
  \pi_{+}(\epsilon'+(j)) \setminus
  \pi_{+}(\epsilon_{0})} e^{-\lambda_{i}(\epsilon'+(j)) }}}
{\displaystyle{\prod_{i \in \pi_{+}(\epsilon)\setminus \pi_{+}(\epsilon_{0})}
(1-e^{-\lambda_{i}(\epsilon'+(j))})}} \right]$}}}.
$$

Chaque terme de cette somme est nulle.
En effet, soit $\epsilon'$ un {\'e}l{\'e}ment de la sommation, 
comme $j$ est le plus grand {\'e}l{\'e}ment de $\pi_{+}(\epsilon) \setminus
\pi_{+}(\epsilon_{0})$, pour tout $i \in \pi_{+}(\epsilon) \setminus
\pi_{+}(\epsilon_{0})$, 
$\lambda_{i}(\epsilon'+(j))=\lambda_{i}(\epsilon')$ si $i\neq j$, et 
$\lambda_{j}(\epsilon'+(j))=-\lambda_{j}(\epsilon')$. Le terme de la somme
associ{\'e} {\`a} $\epsilon'$ est donc :
$$\frac{ \displaystyle{ \prod_{i \in \pi_{+}(\epsilon') \setminus
  \pi_{+}(\epsilon_{0})} e^{-\lambda_{i}(\epsilon') }}}{\displaystyle{\prod_{i 
\in \pi_{+}(\epsilon-(j))\setminus \pi_{+}(\epsilon_{0})}
(1-e^{-\lambda_{i}(\epsilon')})}} \,\,\,\, \left[ \frac{1}{ 1-e^{-\lambda_{j}(\epsilon')} }+
 \frac{e^{\lambda_{j}(\epsilon')}  }{ 1-e^{\lambda_{j}(\epsilon')} }
 \right].$$
 
Ce terme est bien nul d'apr{\`e}s la relation 
$\frac{1}{1-e^{-x}}+\frac{e^{x}}{1-e^{x}}=0$, et on obtient donc $\chi(
\overline{Y_{\epsilon}}, \hat{\mu}_{\epsilon_0}^D) =0$.

\end{proof}

\section{K-th{\'e}orie {\'e}quivariante des vari{\'e}t{\'e}s de Bott-Samelson}

On reprend les notations de la section~\ref{BS}. On choisit $N$ racines simples
$\mu_1, \ldots , \mu_N$ non n{\'e}cessairement distinctes, et on pose $\Gamma = 
\Gamma( \mu_1, \ldots , \mu_N  )$. Pour $1 \leq i<j \leq N$, on pose
$b_{i,j}= \mu_j(\mu^{\vee}_i)$ et $B=\{b_{i,j}\}_{1\leq i<j \leq N}$.

Comme dans le cas de la cohomologie, les r{\'e}sultats de cette section ne d{\'e}pendent
pas de la structure complexe de $\Gamma$, et on identifie donc
$\Gamma$
avec la tour de Bott $Y=Y_B$, la d{\'e}composition 
$\Gamma=\coprod_{\epsilon \in \mathcal{E}} \Gamma_{\epsilon}$ avec la
  d{\'e}composition $Y=\coprod_{\epsilon \in \mathcal{E}} Y_{\epsilon}$, et le point
  fixe $\epsilon \in \Gamma^T$ avec le point fixe $\epsilon \in Y^D$.

Le tore $T$  agit sur $\Gamma$ via son image  $S$, et l'action 
 de $S$ sur $\Gamma$ s'identifie {\`a} celle d'un sous-tore de
$D$ sur $Y$. Pour tout $\epsilon \in \mathcal{E}$, on note alors
$\hat{\mu}^S_{\epsilon}$ l'{\'e}l{\'e}ment de $K_S(\Gamma)$ obtenu {\`a} partir de 
$\hat{\mu}^D_{\epsilon}$ par restriction {\`a} $S$ de l'action de $D$. Ces {\'e}l{\'e}ments
s'obtiennent {\`a} l'aide des fibr{\'e}s de Hopf de la m{\^e}me mani{\`e}re que les classes 
$\hat{\mu}^D_{\epsilon}$, et ils forment donc une base du $R[S]$-module 
$K_S(\Gamma)$.

De plus, pour tout couple $(\epsilon,\epsilon') \in \mathcal{E}^2$ :
$$\chi(\overline{\Gamma_{\epsilon'}}, \hat{\mu}_{\epsilon}^S)
=\delta_{\epsilon', \epsilon}.
$$

Comme le tore $T$ agit sur $\Gamma$ via son image $S$, on a $R[S] \subset
R[T]$, et  
$K_T(\Gamma)$ s'identifie {\`a} 
$K_S(\Gamma)\! \otimes_{R[S]}\!\!R[T]$. Si on pose 
$\hat{\mu}^T_{\epsilon}=\hat{\mu}^S_{\epsilon}\otimes 1 $, on a donc la proposition
suivante :

\begin{prop} \label{propbasekbs}

\indent

\begin{enumerate}

\item[$(i)$] La  $K$-th{\'e}orie $T$-{\'e}quivariante de $\Gamma^T$ s'identifie {\`a}
$F(\mathcal{E};R[T])$.

\item[$(ii)$] La restriction aux points fixes $i_{T}^*$ : $K_{T}(\Gamma)
 \rightarrow F(\mathcal{E};R[T])$ est injective.
    
\item[$(iii)$] La $K$-th{\'e}orie $T$-{\'e}quivariante de $\Gamma$ est un $R[T]$-module libre
 admettant comme base la famille $\{\hat{\mu}_{\epsilon}^T\}_{\epsilon \in
\mathcal{E}}$ qui v{\'e}rifie :

$$\chi(\overline{\Gamma_{\epsilon'}}, \hat{\mu}_{\epsilon}^T)=\delta_{\epsilon',
  \epsilon}.$$
   
\end{enumerate}

\end{prop}





\medskip

Pour tout $\epsilon \in \mathcal{E}$, on pose  $\mu_{\epsilon}^T=
i_{T}^*(\hat{\mu}_{\epsilon}^T)$. Le th{\'e}or{\`e}me~\ref{restrictionkb} 
et une d{\'e}monstration analogue {\`a} celle du
th{\'e}or{\`e}me~\ref{cohomologieBS} nous donnent alors le r{\'e}sultat suivant :

\begin{theo} \label{ktheorieBS}
    
Pour $\epsilon \in \mathcal{E}$, on a :

$$\left\{ \begin{array}{ll} \mu_{\epsilon}^T(\epsilon') = 
\displaystyle{\prod_{i \in \pi_{+}(\epsilon')}e^{\alpha_{i}(\epsilon')}
\prod_{i \in \pi_{+}(\epsilon)}(e^{-\alpha_{i}(\epsilon')}-1) }
 & { \rm si } \hspace{0,15 cm} \epsilon \leq \epsilon', \\ 
 \mu_{\epsilon}^T(\epsilon') = 0
 & { \rm  sinon }.
\end{array}\right.$$
   
\end{theo}

\section{K-th{\'e}orie {\'e}quivariante des vari{\'e}t{\'e}s de drapeaux}

\subsection{D{\'e}finitions}

La vari{\'e}t{\'e} de drapeaux $X$ n'{\'e}tant pas compacte en g{\'e}n{\'e}ral, on d{\'e}finit $K_T(X)$
de la mani{\`e}re suivante : 

Pour tout entier $n \geq 0$, on d{\'e}finit $\displaystyle{X_{n}=\!\!\!\!\!
 \bigcup_{ \tiny 
\begin{array}{ll} \hspace{0,15 cm} w \in W \\ 
 l(w) \leq n
\end{array}}\!\!\!\!\!X_{w}}$. Soit $\mathcal{F}$ la filtration : 
$$\mathcal{F} :
\emptyset=X_{-1} \subset X_{0} \subset X_{1} \subset \cdots , $$ 
alors : 

\begin{enumerate}

\item[$(1)$] chaque $X_{n}$ est un sous espace compact $T$-stable de $X$ et,

\item[$(2)$] la topologie de $X$ est la topologie limite induite par la filtration
    $\mathcal{F}$.

\end{enumerate}

\smallskip

Gr{\^a}ce {\`a} cette filtration, on d{\'e}finit alors la K-th{\'e}orie $T$-{\'e}quivariante de $X$, 
not{\'e}e $K_{T}(X)$, par
$\displaystyle{K_{T}(X)=\lim_{  \leftarrow}}_{ 
 \tiny n \rightarrow +\infty
}\!\!\!\!\! \!\!\!\!\!\!\!\! K_{T}(X_{n})$. Cette d{\'e}finition est ind{\'e}pendante de la filtration $\mathcal{F}$
v{\'e}rifiant $(1)$ et $(2)$.

\bigskip

On note $F(W; Q[T])$) la 
$R[T]$-alg{\`e}bre des fonctions de $W$ {\`a}
valeurs dans $Q[T]$ munie de
l'addition et de la multiplication point par point. Pour un poids entier
$\lambda \in i\mathfrak{t}^*$, on note $e^{\alpha}
: T \rightarrow S^1$ le caract{\`e}re correspondant. Pour tout $1 \leq i \leq r$,
on d{\'e}finit alors un op{\'e}rateur de Demazure $D_{i}$ sur $F(W; Q[T])$ par : 

$$(D_{i}f)(v)=\frac{f(v)-f(vs_{i})e^{-v\alpha_{i}}}{1-e^{-v
\alpha_{i}}}.$$
   
Dans \cite{kkk}, Kostant et Kumar montrent que ces op{\'e}rateurs de Demazure 
v{\'e}rifient les relations de tresses de $W$. Pour tout $w \in W$,
on peut donc d{\'e}finir un op{\'e}rateur $D_{w}$ sur $F(W; Q[T])$ d{\'e}fini par
$D_w=D_{i_1}D_{i_2} \cdots D_{i_l}$ si $w=s_{i_1}s_{i_2} \cdots s_{i_l}$ est une
d{\'e}composition r{\'e}duite de $w$. 

De plus, pour tout 
$1 \leq i \leq r$, $D_i^2=D_i$. Donc, si pour $\underline{u} \in 
\underline{W}$, on note $D_{\underline{u}}=D_u$, alors pour tout couple 
$(\underline{v},\underline{w})\in \underline{W}^2$, $ D_{\underline{v}}D_{\underline{w}}
=D_{\underline{v}\, \underline{w}}$.

\medskip

On note $\Psi$ la
sous-alg{\`e}bre de $F(W; R[T])$ d{\'e}finie par : 
$$ \Psi = \{ f \in F(W; R[T]), \hspace{0,2 cm} 
{\rm telles \hspace{0,2 cm} que } \hspace{0,2 cm}
\forall w \in W, \hspace{0,2 cm} D_{w}f \in F(W; R[T]) \}.$$

\medskip

  L'ensemble des points fixes $X^T \approx  W$ 
{\'e}tant discret, on
peut identifier $K_{T}(X^T)$ avec $F(W; R[T])$ et on obtient ainsi un morphisme 
$i_{T}^* : K_{T}(X) \rightarrow F(W; R[T])$. On notera~$*$ l'involution de
$K_{T}(X)$ d{\'e}finie par la dualit{\'e} des fibr{\'e}s, et on notera de la m{\^e}me fa{\c c}on
l'involution de $R[T]$ d{\'e}finie sur les caract{\`e}res par
$*(e^{\lambda})=e^{-\lambda}$, ce qui induit une involution de $F(W; R[T])$.
 Pour tout {\'e}l{\'e}ment $\tau \in K_{T}(X)$, 
$*i_{T}^*(\tau)=i_{T}^*(*\tau)$. Le r{\'e}sultat suivant est prouv{\'e} dans \cite{kkk} :

\begin{prop} \label{propositionKKK}

L'application $i_{T}^*$ est injective, et l'image de $K_{T}(X)$ par cette
application est {\'e}gale {\`a} $\Psi$. De plus, $\Psi=\prod_{w \in W}R[T]\psi^w$, o{\`u}
les fonctions $\psi^w$ sont uniquement d{\'e}termin{\'e}es par les relations : 
$$\forall (v,w) \in W^2, \, D_v(\psi^w)(1)=\delta_{v,w}.$$

De plus, les fonctions  $\psi^w$
v{\'e}rifient les propri{\'e}t{\'e}s suivantes :  

\smallskip

\begin{enumerate}

 \item[$(i)$] $\psi^{w}(v)=0$ sauf si $w \leq v$,
   
  \item[$(ii)$] $\psi^{w}(w) = \prod_{\beta \in \Delta(w^{-1})}
(1-e^{\beta})$,

 \item[$(iii)$]  $ \left\{ \begin{array}{ll} D_{i}\psi^{w} = \psi^{w} + \psi^{ws_{i}}
 & {\rm si } \hspace{0,2 cm} ws_{i}<w, \\ 
D_{i}\psi^{w}=0
 & {\rm si } \hspace{0,2 cm} ws_{i}>w,  \end{array} \right.$

\item[$(iv)$] $ \forall v \in W, \psi^{1}(v) =e^{\rho-v\rho} $.
    
\end{enumerate}

\end{prop}

\begin{rema}

Un {\'e}l{\'e}ment $f=(a_{w})_{w \in W}$ de $\prod_{w \in W}R[T]\psi^{w}$ est bien une 
fonction de $W$ {\`a} valeurs dans $R[T]$. En effet soit $v \in W$, d'apr{\`e}s la
propri{\'e}t{\'e}~$(i)$, 
$\sum_{w \in W}a_{w}\psi^{w}(v)$ est une somme finie o{\`u} les termes
{\'e}ventuellement non nuls correspondent aux {\'e}l{\'e}ments $u$ de $W$ qui v{\'e}rifient $u
\leq v$.

\end{rema}

\begin{rema}

Les fonctions $\psi^w$ sont uniquement d{\'e}termin{\'e}es 
par les propri{\'e}t{\'e}s $(i)$, $(ii)$, $(iii)$ et $(iv)$ de la proposition
pr{\'e}c{\'e}dente. 

\end{rema}

On pose $\hat{\psi}^w=(i_{T}^*)^{-1}(\psi^{w})$, et pour $v \in W$ on note
$\hat{D}_{v} : K_{T}(X) \rightarrow K_{T}(X)$ l'application induite par $D_v :
\Psi \rightarrow \Psi$.

\begin{rema}

Dans le cas fini, $K_T(X)$ s'identifie {\`a} $K^0(H,X)$ (voir \cite{kkk}), et
Kostant et Kumar montrent dans \cite{kkk} que la base $\{\hat{\psi}^w\}_{w\in
  W}$ de $K_T(X) \simeq K^0(H,X)$ est r{\'e}li{\'e}e aux vari{\'e}t{\'e}s de Schubert par les
relations : 
$$ \forall (v,w) \in W^2, \chi(\overline{X_v}, *\hat{\psi}^w)=\delta_{v,w}.
$$
\end{rema}

\bigskip

 Dans \cite{kkk},  Kostant et Kumar composent $i_{T}^{*}$ avec $\phi : 
F(W; Q[T]) \rightarrow
F(W; Q[T])$ d{\'e}finie par 
$\phi(f)(w)=f(w^{-1})$ pour tout {\'e}lement $f$ de $F(W; Q[T])$ et tout $w \in
W$. Ils trouvent alors la sous alg{\`e}bre $\Psi'$ (not{\'e}e $\Psi$ dans \cite{kkk}) 
de $F(W;R[T])$ : 
$$\Psi' = \{ f \in F(W;R[T]), \hspace{0,2 cm} {\rm telles \hspace{0,2 cm} que } \hspace{0,2 cm}
\forall w \in W, \hspace{0,2 cm} D'_{w}f \in F(W;R[T]) \},$$
o{\`u} les op{\'e}rateurs $D_{w}'$ sont d{\'e}finis {\`a} partir des op{\'e}rateurs $D_{i}'$ donn{\'e}s
par : 
$$(D_{i}'f)(v)=\frac{f(v)-f(s_{i}v)e^{-v^{-1}\alpha_{i}}}{1-e^{-v^{-1}
\alpha_{i}}}.$$

Ils consid{\`e}rent la base $\psi_{w} '$ (not{\'e}e $\psi^{w}$ dans \cite{kkk}) 
de $\Psi '$ reli{\'e}e {\`a} la base $\psi^{w}$ de la proposition~\ref{propositionKKK} par la relation 
$\psi_{w}'=\phi(\psi^{w^{-1}})$. Pour tout couple $(w,v) \in W^2$, $\psi_{w}'(v)=
\psi^{w^{-1}}(v^{-1})$.

\subsection{Restrictions aux points fixes}

Soit $v \in W$ et soit $v=s_{i_{1}} \cdots s_{i_{l}}$ une d{\'e}composition non
n{\'e}cessairement r{\'e}duite de
$v$. Comme dans la section~\ref{cohomologievd}, pour $1\leq j \leq l$, on notera
$\beta_{j}$ l'{\'e}l{\'e}ment de  
$\mathfrak{h}^*$ d{\'e}fini par $\beta_{j}=s_{i_{1}} \cdots s_{i_{j-1}}
\alpha_{i_{j}}$.

\begin{theo} \label{ktheorieVD}

 Si $w \in W$ est tel que $w \leq v$, on a la formule suivante :    

   $$ \psi^{w}(v) =e^{\rho - v\rho}\sum_{l(w) \leq m \leq l}
 \sum
 (e^{-\beta_{j_{1}}}-1) \cdots (e^{-\beta_{j_{m}}}-1),
$$
    o{\`u} la deuxi{\`e}me somme porte sur l'ensemble des entiers $1\leq j_{1} 
   < \cdots < j_{m} \leq l$ tels que $\underline{s_{i_{j_{1}}}} \ldots   
  \underline{s_{i_{j_{m}}}}=\underline{w}$.
    
\end{theo}

Donnons quelques exemples de calculs pour expliciter cette 
formule.

\smallskip

Tout d'abord pour tout $v\in
W$, on retrouve bien $\psi^{1}(v)=e^{\rho-v\rho}$, puisque la seule fa{\c c}on
de ``trouver $1$ en dessous de $v$'' est de prendre la suite vide.

\medskip

Pla{\c c}ons nous dans le cas $A_4$ et calculons $\psi^{w}(v)$
avec $w=s_{3}s_{2}$ et $v=s_{2}s_{3}s_{2}s_{1}s_{2}$. Il y a $3$ fa{\c c}ons de
``trouver $\underline{w}$ en dessous de $v$'' :
$\underline{w}=\underline{s_{i_{2}}}\hspace{0,1 cm} \underline{s_{i_{3}}}$, 
$\underline{w}=\underline{s_{i_{2}}}\hspace{0,1 cm} \underline{s_{i_{5}}}$, 
$\underline{w}=\underline{s_{i_{2}}}\hspace{0,1 cm}
\underline{s_{i_{3}}}\hspace{0,1 cm}\underline{s_{i_{5}}}$, et on trouve donc :
$$\psi^{w}(v)=e^{\alpha_{2}+(\alpha_{2}+\alpha_{3})+\alpha_{3}+(\alpha_{1}+\alpha_{2}+\alpha_{3})+(\alpha_{1}+\alpha_{2})}[(e^{-(\alpha_{2}+\alpha_{3})}-1)(e^{-\alpha_{3}}-1)$$
$$+
(e^{-(\alpha_{2}+\alpha_{3})}-1)(e^{-(\alpha_{1}+\alpha_{2})}-1)  +
(e^{-(\alpha_{2}+\alpha_{3})}-1)(e^{-\alpha_{3}}-1)(e^{-(\alpha_{1}+\alpha_{2})}-1)]$$
$$=e^{2\alpha_{1}+4\alpha_{2}+3\alpha_{3}}[1+e^{-(\alpha_{1}+2\alpha_{2}+2\alpha_{3})}-
e^{-(\alpha_{2}+\alpha_{3})}-
e^{-(\alpha_{1}+\alpha_{2}+\alpha_{3})}]$$
$$=e^{2\alpha_{1}+4\alpha_{2}+3\alpha_{3}}+e^{\alpha_{1}+2\alpha_{2}+\alpha_{3}}-
e^{2\alpha_{1}+3\alpha_{2}+2\alpha_{3}}-e^{\alpha_{1}+3\alpha_{2}+2\alpha_{3}}.$$

\medskip

Pour d{\'e}montrer le th{\'e}or{\`e}me, nous aurons besoin du lemme suivant : 

\begin{lemm} \label{*chi}

Soit $v=s_{\mu_{1}} \cdots s_{\mu_{N}}$ une d{\'e}composition non n{\'e}cessairement
r{\'e}duite d'un {\'e}l{\'e}ment $v$ de $W$. On pose $\Gamma=
\Gamma(\mu_{1}, \ldots, \mu_{N})$ et $g=g_{\mu_{1}, \ldots, \mu_{N}} : \Gamma
\rightarrow X$. Pour tout
$\epsilon \in \mathcal{E}$, on a alors : 
$$\forall \tau \in K_{T}(X) \hspace{0,1 cm}, \hspace{0,2 cm} \chi(
\overline{\Gamma_{\epsilon}},g^*(*\tau))=*(D_{\underline{v}(\epsilon)}i_{T}^*(\tau))(1).$$

\end{lemm}

\begin{proof}
 Soit $\tau \in K_{T}(X)$, on proc{\`e}de par r{\'e}currence sur $l(\epsilon)$. 
Le r{\'e}sultat est trivial si $l(\epsilon)=0$. 

Supposons le r{\'e}sultat v{\'e}rifi{\'e} pour
tout $\epsilon'$ de longueur strictement inf{\'e}rieure {\`a} $p$, et soit $\epsilon$ de
longueur $p$. En utilisant la formule~\ref{ab}, on obtient :

$$ \chi(
\overline{\Gamma_{\epsilon}},g^*(*\tau))= \sum_{\epsilon' \leq
\epsilon}\frac{i_{T}^*(g^*(*\tau))(\epsilon')}{\prod_{i \in \pi_{+}(\epsilon)}
(1-e^{\alpha_{i}(\epsilon')})}= \sum_{\epsilon' \leq
\epsilon}\frac{*i_{T}^*(\tau)(v(\epsilon'))}{\prod_{i \in \pi_{+}(\epsilon)}
(1-e^{\alpha_{i}(\epsilon')})}.$$

Soit $j$ le plus grand {\'e}l{\'e}ment de $\pi_+(\epsilon)$, et soit $\tilde{\epsilon}=
\epsilon -(j)$. En distinguant les {\'e}l{\'e}ments $\epsilon'$ tels que $\epsilon'_j=0$ et
ceux tels que $\epsilon'_j=1$, on obtient : 

$$ \chi(
\overline{\Gamma_{\epsilon}},g^*(*\tau))= \sum_{\epsilon' \leq
\tilde{\epsilon}}\frac{*i_{T}^*(\tau)(v(\epsilon'))}
{(1-e^{\alpha_{j}(\epsilon')})\prod_{i \in
  \pi_{+}(\tilde{\epsilon})}
(1-e^{\alpha_{i}(\epsilon')})}+ \sum_{\epsilon' \leq
\tilde{\epsilon}}\frac{*i_{T}^*(\tau)(v(\epsilon')s_{\mu_j})}
{(1-e^{-\alpha_{j}(\epsilon')})\prod_{i \in
  \pi_{+}(\tilde{\epsilon})}
(1-e^{\alpha_{i}(\epsilon')})},$$
d'o{\`u} : 
$$ \chi(
\overline{\Gamma_{\epsilon}},g^*(*\tau))= \sum_{\epsilon' \leq
\tilde{\epsilon}}\frac{1}{\prod_{i \in
  \pi_{+}(\tilde{\epsilon})}
(1-e^{\alpha_{i}(\epsilon')})}
\left[\frac{i_{T}^*(\tau)(v(\epsilon'))-e^{-v(\epsilon')\alpha_{j}}
i_{T}^*(\tau)(v(\epsilon')s_{\mu_j})}
{1-e^{-v(\epsilon')\alpha_{j}}}\right]^*$$
$$= \sum_{\epsilon' \leq
\tilde{\epsilon}}\frac{*i_{T}^*(\hat{D}_{s_{\mu_j}}\tau)
(v(\epsilon'))}{\prod_{i \in \pi_{+}(\epsilon)}
(1-e^{\alpha_{i}(\epsilon')})},
$$
i.e., d'apr{\`e}s la formule initiale : 
$$ \chi(
\overline{\Gamma_{\epsilon}},g^*(*\tau))=\chi(
\overline{\Gamma_{\tilde{\epsilon}}},g^*(*\hat{D}_{s_{\mu_j}}\tau)),$$
et donc par hypoth{\`e}se de r{\'e}currence ($\tilde{\epsilon}$ {\'e}tant de longueur $p-1$) : 
$$\chi(
\overline{\Gamma_{\epsilon}},g^*(*\tau))=*(D_{\underline{v}(\tilde{\epsilon})}
i_{T}^*(\hat{D}_{s_{\mu_j}}\tau))(1)=*(D_{\underline{v}(\tilde{\epsilon})}
D_{s_{\mu_j}}i_{T}^*(\tau))(1)=*(D_{\underline{v}(\epsilon)}i_{T}^*(\tau))(1).$$

\end{proof}

\begin{proof}[D{\'e}monstration du th{\'e}or{\`e}me~\ref{ktheorieVD}]

Soit $v=s_{\mu_{1}} \cdots s_{\mu_{N}}$ une d{\'e}composition non n{\'e}cessairement
r{\'e}duite d'un {\'e}l{\'e}ment $v$ de $W$. On pose $\Gamma=
\Gamma(\mu_{1}, \ldots, \mu_{N})$ et $g=g_{\mu_{1}, \ldots, \mu_{N}}$. Soit 
$w$ un {\'e}l{\'e}ment quelconque de $W$.

 D'apr{\`e}s le lemme~\ref{*chi} :
$$\forall \tau \in K_{T}(X) \hspace{0,1 cm}, \hspace{0,2 cm} \chi(
\overline{\Gamma_{\epsilon}},g^*(*\tau))=*(D_{\underline{v}(\epsilon)}i_{T}^*(\tau))(1).$$

 Or, d'apr{\`e}s la caract{\'e}risation de la base
$\{\psi^w\}_{w \in W}$ (proposition~\ref{propositionKKK}) : 
$$\forall (u,w) \in W^2, (D_{u}(\psi^{w}))(1)=\delta_{u,w}.$$

On d{\'e}duit des deux formules pr{\'e}c{\'e}dentes que pour tout $\epsilon \in \mathcal{E}$, on a
\begin{equation} \label{chi}
  \chi(\overline{\Gamma_{\epsilon}},g^*(*\hat{\psi}^w))=
\delta_{\underline{v}(\epsilon)
 ,\underline{w}}. \hspace{2 cm} \end{equation}

D'apr{\`e}s la caract{\'e}risation de la base
$\{\hat{\mu}_{\epsilon}^T\}_{\epsilon \in \mathcal{E}}$, on a donc :

$$g^*(\hat{\psi}^{w})=\!\!\!\!\!\!\!\sum_{\epsilon \in \mathcal{E},
 \, \underline{v}(\epsilon)=\underline{w}}\!\!\!\!\!\!\!*
\hat{\mu}_{\epsilon}^T \hspace{0,1 cm}.$$

Or, comme en cohomologie, on a le diagramme commutatif suivant : 
$$\xymatrix{
     K_{T}(\Gamma)\ar[dd]^{i_{T}^*} & &
     \ar[ll]^{g^{*}}K_{T}(X)\ar[dd]^{i_{T}^*} \\
     \\
     F(\mathcal{E};R[T]) & & 
   \ar[ll]^{\hat{u}}F(W;R[T]) }
  \\
$$
o{\`u} l'application $\hat{u}$ est d{\'e}duite de l'application 
$\overline{u} : \mathcal{E} \rightarrow W$ d{\'e}finie par
$\overline{u}(\epsilon)=v(\epsilon)$.
 On en d{\'e}duit donc pour tout $w \in W$ :    
 $$\displaystyle{\psi^{w}(v)=(\hat{u}i_{T}^*(
\hat{\psi}^{w}))((\bold{1})) =i_T^*g^*(\hat{\psi}^{w})((\bold{1}))
 =\!\!\!\!\!\!\sum_{\epsilon \in \mathcal{E},
  \, \underline{v}(\epsilon)= 
\underline{w}}\!\!\!\!\!\!\!*\mu_{\epsilon}^T((\bf{1}\rm))},$$
 ce qui nous donne
bien le th{\'e}or{\`e}me~\ref{ktheorieVD} {\`a} l'aide du th{\'e}or{\`e}me~\ref{ktheorieBS} et de la
formule suivante : 
\begin{equation} \label{rho1}
\displaystyle{\rho - v\rho = \sum_{j=1}^{k}\beta_{j}}.
\end{equation}

Cette formule est une cons{\'e}quence imm{\'e}diate de la proposition $29$ de la section
VI$.1.10$ de \cite{bou} qui affirme que pour tout $1 \leq i \leq r$ : 
$$s_i(\rho)=\rho - \alpha_i.$$

\end{proof}

\section{Une autre base de $K_{T}(X)$}

Dans ce paragraphe, on se place dans le cas fini (i.e. $W$ fini $\Leftrightarrow
\mathfrak{g}$ de dimension finie). On note $w_{0}$ le plus grand
{\'e}l{\'e}ment de $W$. On choisit $w_{0} = s_{\mu_{1}} \cdots  s_{\mu_{N}}$ une
d{\'e}composition r{\'e}duite de $w_{0}$ et on pose $\Gamma=
\Gamma(\mu_{1}, \ldots, \mu_{N})$ et $g=g_{\mu_{1}, \ldots, \mu_{N}} : \Gamma
\rightarrow X$.

\smallskip

Comme $X$ est une vari{\'e}t{\'e} complexe projective lisse, $K_{0}(H,X)$ s'identifie 
{\`a} $K^{0}(H,X)$. De plus, dans \cite{kkk}, Kostant et Kumar montrent que
le morphisme canonique $K^{0}(H,X) \rightarrow K_T(X)$ est un isomorphisme. Dans
la suite, on identifie donc ces trois groupes.

 On sait alors que la d{\'e}composition en cellules de Schubert $X=\coprod_{w \in
  W}X_{w}$ fournit une base $\{[\mathcal{O}_{\overline{X_w}}]\}_{w \in W}$ de
$K_{0}(H,X)$. Les classes $[\mathcal{O}_{\overline{X_w}}]$ sont d{\'e}finies par le faisceau
structural de $\overline{X_{w}}$ prolong{\'e} par $0$ sur $X \setminus
\overline{X_{w}}$. Pour $w \in W$, on pose
$\hat{\gamma}^w=*[\mathcal{O}_{\overline{X_w}}] \in K_{T}(X)$, et
$\gamma^{w}=i_{T}^*(\hat{\gamma}^w)$.

\smallskip

Le but de cette section est de d{\'e}terminer les {\'e}l{\'e}ments $a_{w}^v \in R[T]$ d{\'e}finis
par :
$$\hat{\gamma}^{w}=\sum_{v \in W}a_{w}^v\hat{\psi}^v.
$$

Pour cela, on a besoin du r{\'e}sultat suivant prouv{\'e} dans \cite{kkk} : 

\begin{prop}

Pour tout $w \in W$ et tout entier $1 \leq i \leq r$,

$$D_{i}(\gamma^w)=\left\{ \begin{array}{ll} \gamma^w 

 & {\rm si } \hspace{0.2 cm} ws_{i}<w, \\ 
\gamma^{ws_{i}}
 & {\rm si } \hspace{0.2 cm} ws_{i}>w.  \end{array} \right.$$

\end{prop}

 Soit $w \in W$ et soit $s_{i}$ 
une r{\'e}flexion simple telle que $ws_{i}>w$. Si on
applique l'op{\'e}rateur $D_{i}$ {\`a} la d{\'e}composition 
$\displaystyle{\gamma^{w}=\sum_{v \in W}a_{w}^v\psi^v}$, on
obtient : 
$$\gamma^{ws_{i}}=D_{i}(\gamma^w)=
\sum_{v \in W}a_{w}^v D_{i}\psi^{v}.$$

En utilisant les relations v{\'e}rifi{\'e}es par les fonctions $\psi^{v}$, on trouve
alors : 
$$\sum_{v \in W}a_{ws_{i}}^{v}\psi^{v}=\sum_{v \in W, vs_{i}<v}a_{w}^{v}(\psi^{v} + 
\psi^{vs_{i}})=\sum_{v \in W, vs_{i}<v}a_{w}^{v}\psi^{v}
+ \sum_{v \in W, vs_{i}>v}a_{w}^{vs_{i}}\psi^{v}. $$

On obtient donc la relation de r{\'e}currence suivante sur les coefficients
$a_{w}^{v}$ : 

$$a_{ws_{i}}^{v}=\left\{ \begin{array}{ll} a_{w}^{v}
 & {\rm si } \hspace{0.2 cm} vs_{i}<v, \\ 
a_{w}^{vs_{i}}
 & {\rm si } \hspace{0.2 cm} vs_{i}>v.  \end{array} \right.$$

De ces relations, on d{\'e}duit en utilisant les relations~\ref{hecke1} et
\ref{hecke2} : 

\begin{equation} \label{rec} 
\forall (w,v) \in W^2, \quad
a_{w}^{v}=a_{1}^{\underline{v} \hspace{0,1 cm} \underline{w^{-1}}} ,
 \end{equation}
 o{\`u} pour $\underline{u} \in \underline{W}$, on a pos{\'e} 
$a_{1}^{\underline{u}}=a_{1}^{u}$. Il suffit donc de trouver la d{\'e}composition de
$\gamma^{1}$. 

Pour cela, on aura besoin des valeurs de $\gamma^{1}$ : 
$$\left\{ \begin{array}{l} \gamma^{1}(1)=\prod_{\alpha \in
      \Delta_{+}}
(1-e^{-\alpha}),
   \\ 
\gamma^{1}(v)=0
 \hspace{0,5 cm} {\rm si } \hspace{0.2 cm} v \neq 1.  \end{array} \right.$$

La valeur de $ \gamma^{1}(1)$ est calcul{\'e}e {\`a} l'aide de la formule
d'auto-intersection (voir \cite{ginzburg}, proposition $5.4.10$), 
et les autres valeurs sont nulles par le th{\'e}or{\`e}me de localisation.

Comme on a $\hat{\gamma}^{1}=\sum_{v \in
    W}a_{1}^{v}\hat{\psi}^{v}$, soit $v \in W$ et $\epsilon \in \mathcal{E}$ 
tel que $u(\epsilon)=v$, d'apr{\`e}s la formule~\ref{chi}, le coefficient $a_{1}^{v}$
est donn{\'e} par :
$$a_{1}^{v}=*\chi(\overline{\Gamma_{\epsilon}},g^*(*\hat{\gamma}^{1})).$$

En utilisant la formule~\ref{ab} et les valeurs de
$*\gamma^{1}=*i_{T}^*\hat{\gamma}^{1}$, on obtient alors : 
$$a_{1}^{v}=  \sum_{\epsilon' \leq
\epsilon, u(\epsilon')=1}\frac{\prod_{\alpha \in
      \Delta_{+}}
(1-e^{-\alpha})}{\prod_{i \in \pi_{+}(\epsilon)}
(1-e^{-\alpha_{i}(\epsilon')})}. $$

On a donc la proposition suivante : 

\begin{prop}

Soit $v \in W$ et soit $v=s_{i_{1}} \cdots s_{i_{l}}$ une d{\'e}composition r{\'e}duite de
$v$. Pour tout sous-ensemble $I$ de $\{1, \ldots , l \}$ et tout entier $1 \leq i
\leq l$, on pose $\displaystyle{\beta_{i}(I)=(\prod_{j \in I , j \leq
  i}s_{i_{j}})\alpha_{i}}$ (\:$\beta_{i}(I)=\alpha_{i}$ si $I \cap \{ 1, \ldots ,
i \} = \emptyset$). Alors 
$$\displaystyle{\sum_{k=1}^{l}\sum
  \frac{\prod_{\alpha \in
  \Delta_{+}}(1-e^{-\alpha})
}{
\prod_{i=1}^{l}
(1-e^{-\beta_{i}(\{j_{1},\ldots , j_{k}\}  )})}},$$
 o{\`u} la deuxi{\`e}me somme porte sur l'ensemble des indices 
$1 \leq j_{1} < \cdots < j_{k}\leq l$ tels que $s_{i_{j_{1}}} \cdots s_{i_{j_{k}}}=1$, est
un {\'e}l{\'e}ment de $R[T]$ qui ne d{\'e}pend pas du choix d'une d{\'e}composition r{\'e}duite de
$v$. Si on note $b^{v}$ cet {\'e}l{\'e}ment, alors $a_{1}^{v}=b^{v}$.
 
\end{prop}

\begin{rema}

La proposition pr{\'e}c{\'e}dente est encore valable si on prend une d{\'e}composition non
r{\'e}duite de $v$.

\end{rema}

Si on utilise la relation~\ref{rec}, et si pour $\underline{v} \in \underline{W}$, on pose
$b^{\underline{v}}=b^v$, on obtient alors le th{\'e}or{\`e}me suivant :

\begin{theo} \label{changementdebases} 

Soit $w \in W$, alors : 
$$*[\mathcal{O}_{\overline{X_{w}}}]=\sum_{v \in W}
b^{\underline{v}\hspace{0,1 cm} \underline{w^{-1}} }\hat{\psi}^{v}.$$

\end{theo}

Soit $Q_{W}$ le $Q[T]$-module libre qui admet pour base la famille $\{\delta_{w}\}_{w
  \in W}$ et qu'on munit d'une structure d'anneau d{\'e}finie par : 
$$(q_{1}\delta_{w_{1}}).(q_{2}\delta_{w_{2}})=q_{1}(w_{1}q_{2})\delta_{w_{1}w_{2}}
\hspace{0,1 cm}, \hspace{0,1 cm} \forall (q_{1}, q_{2}) \in Q[T]^2 , \hspace{0,1 cm} { \rm et }
\hspace{0,1 cm}  (w_{1},w_{2}) \in W^2,$$
o{\`u} l'action de $W$ sur $Q[T]$ est d{\'e}duite de celle de $W$ sur $T$. Dans
\cite{kkk}, B. Kostant et S. Kumar introduisent des {\'e}l{\'e}ments $\{y_{i} \}_{1 \leq
  i \leq r}$ de $Q[T]$ d{\'e}finis par : 
$$y_{i}=\frac{1}{1-e^{-\alpha_{i}}}(\delta_{1}-e^{-\alpha_{i}}\delta_{s_{i}}).$$

Les $y_{i}$ v{\'e}rifiant les relations de tresses, on peut
d{\'e}finir un {\'e}l{\'e}ment $y_{w} \in Q_{W}$ pour tout $w \in W$. 
On d{\'e}finit alors des {\'e}l{\'e}ments $\{ b_{v,w}
\}_{(v,w) \in W^2}$ de $Q[T]$ par  : 
$$y_{v^{-1}}=\sum_{w \in W}b_{v,w}\delta_{w^{-1}}.$$
 
D'apr{\`e}s l'expression combinatoire de $b_{v,w}$ donn{\'e}e par le lemme $3.5$ de
\cite{ku} :  
$$b^v=b_{v^{-1},1}\prod_{\alpha \in \Delta_{+}}(1-e^{-\alpha}).$$ 
 
De plus, dans \cite{ku} S. Kumar montre que quand  $\overline{X_{v}}$
est lisse : $$\displaystyle{b_{v,1}=\prod_{\gamma \in
    S(v)}(1-e^{-\gamma})^{-1}},$$ 
o{\`u} pour $u \in W$, $S(u)=\{ \alpha \in R^{+}, s_{\alpha} \leq u \}$. En
particulier $b^{w_{0}}=1$, et donc d'apr{\`e}s le th{\'e}or{\`e}me~\ref{changementdebases} : 
$$*[\mathcal{O}_{\overline{X_{w_{0}}}}]=\sum_{w \in W}\hat{\psi}_{w}.$$

\section{Lien avec les alg{\`e}bres de Hecke} \label{hecke}

Dans cette section, on va donner une d{\'e}monstration purement combinatoire du
th{\'e}or{\`e}me~\ref{ktheorieVD}. On se limite au cas d'une d{\'e}composition
r{\'e}duite de $v$.

Cette d{\'e}monstration est similaire {\`a} celle donn{\'e}e par Sarah Billey
dans \cite{sb} dans le cas de la cohomologie {\'e}quivariante. Le th{\'e}or{\`e}me~\ref{independance}, dont
on va donner une id{\'e}e de la d{\'e}monstration, est d{\'e}montr{\'e} pour le type $A$ dans
\cite{llt}. Contrairement {\`a} ce qu'on a fait pr{\'e}c{\'e}demment, on doit d'abord
d{\'e}montrer que l'expression du th{\'e}or{\`e}me~\ref{ktheorieVD} est ind{\'e}pendante du choix d'une
d{\'e}composition de $v \in W$. C'est pour cela qu'on utilise les alg{\`e}bres de Hecke.

\medskip

Soit $A$ un anneau commutatif. On d{\'e}finit $\mathcal{H}$  comme la
$A$-alg{\`e}bre engendr{\'e}e par $\{u_{i}\}_{1 \leq i \leq r}$ soumis aux relations de
tresses d{\'e}finissant $W$ et aux relations $u_{i}^2=u_{i}$. C'est une alg{\`e}bre de
Hecke (voir \cite{hum}). Soit $w \in W$, on
peut d{\'e}finir $u_{w} \in \mathcal{H}$ par $u_{w}=u_{i_{1}} \cdots u_{i_{l}}$ o{\`u}
$w=s_{i_{1}}
\cdots s_{i_{l}}$ est une d{\'e}composition r{\'e}duite quelconque de $w$. Les {\'e}l{\'e}ments 
$\{u_{w}\}_{w \in W}$ forment une base du $A$-module $\mathcal{H}$.

Soit $s_{i_{1}}, \ldots , s_{i_{k}}$ une suite de r{\'e}flexions simples et soit
$\underline{w}=\underline{s}_{i_{1}} \cdots \underline{s}_{i_{k}} \in
\underline{W}$. 
D'apr{\`e}s les relations v{\'e}rifi{\'e}es par les $u_{i}$, $u_{i_{1}} \cdots u_{i_{k}}=u_{w}$.

 Pour tout $1 \leq i \leq r$, on d{\'e}finit la fonction :
 $$h_{i} : \begin{array}{l}A \rightarrow \mathcal{H} \\
x \mapsto 1 + (x-1)u_{i}.\end{array}$$

 On v{\'e}rifie que ces fonctions $h_{i}$
satisfont les relations suivantes ({\'e}nonc{\'e}es sous une forme diff{\'e}rente dans
\cite{fk}) : 

\begin{prop} \label{YangBaxter}

Soient $1 \leq i,j \leq r$ des entiers distincts, deux {\'e}l{\'e}ments quelconques $x$
et $y$ de $A$ v{\'e}rifient les {\'e}quations suivantes :  

$$\left\{ \begin{array}{llll} h_{i}(x)h_{j}(y)=h_{j}(y)h_{i}(x)
 & {\rm si } \hspace{0.2 cm} (s_{i}s_{j})^2 = 1, \\ 
 h_{i}(x) h_{j}(xy)  h_{i}(y)=h_{j}(y)h_{i}(xy)h_{j}(x)
 & {\rm si } \hspace{0.2 cm}(s_{i}s_{j})^3 = 1, \\ 
 h_{i}(x) h_{j}(xy) h_{i}(xy^2) h_{j}(y)=h_{j}(y)h_{i}(xy^2) h_{j}(xy) h_{i}(x)
 & {\rm si } \hspace{0.2 cm} (s_{i}s_{j})^4 = 1, \\ 
h_{i}(x) h_{j}(x^3y) h_{i}(x^2y) h_{j}(x^3y^2)h_{i}(xy) h_{j}(y)
 &   \\ 
=h_{j}(y) h_{i}(xy) h_{j}(x^3y^2) h_{i}(x^2y) h_{j}(x^3y) h_{i}(x)
 & {\rm si } \hspace{0.2 cm}  (s_{i}s_{j})^6 = 1.

\end{array} \right.$$

\end{prop}

\medskip

Dans la suite, on prendra pour $A$ l'anneau $R[T]$.
Soit $w \in W$ et $w=s_{i_{1}}\cdots s_{i_{l}}$ une d{\'e}composition r{\'e}duite de
$w$. On d{\'e}finit un {\'e}l{\'e}ment de $\mathcal{H}$ par : 

$$\mathcal{R}_{i_{1},\ldots , i_{l}} = \prod_{j=1}^{l}h_{i_{j}}(e^{-\beta_{i_{j}}}).$$

A l'aide de la proposition~\ref{YangBaxter}, une d{\'e}monstration analogue {\`a} celle donn{\'e}e par
S. Billey dans \cite{sb} dans le cas de l'alg{\`e}bre nil-Coxeter (voir la
d{\'e}finition dans \cite{sb}) nous donne le
r{\'e}sultat suivant :

\begin{theo} \label{independance}

Soit $w \in W$. L'{\'e}l{\'e}ment $\mathcal{R}_{i_{1},\ldots , i_{l}}$ de $\mathcal{H}$ est
ind{\'e}pendant du choix d'une d{\'e}composition r{\'e}duite $w=s_{i_{1}} \cdots s_{i_{l}}$ de
$w$. Il ne d{\'e}pend que de $w$ et on le notera donc $\mathcal{R}_{w}$.

\end{theo}

Donnons une id{\'e}e de la d{\'e}monstration. D'apr{\`e}s la d{\'e}finition de
$\mathcal{R}_{i_{1},\ldots , i_{l}}$ et d'apr{\`e}s la connexit{\'e} du graphe des
d{\'e}compositions r{\'e}duites de $w$ (qui est une cons{\'e}quence imm{\'e}diate de la
propri{\'e}t{\'e} d'{\'e}change {\'e}nonc{\'e}e par exemple dans \cite{bou}, IV$.1.5$), 
on peut se contenter de regarder ce qui se passe
pour un {\'e}l{\'e}ment $w$ correspondant {\`a} une relation de tresses. Prenons par
exemple $w=s_{i}s_{j}s_{i}=s_{j}s_{i}s_{j}$. Alors $\mathcal{R}_{i, j ,
  i}=h_{i}(e^{-\alpha_{i}})h_{j}(e^{-\alpha_{i}-\alpha_{j}})h_{i}(e^{-\alpha_{j}})$ et 
$\mathcal{R}_{j, i ,
  j}=h_{j}(e^{-\alpha_{j}})h_{i}(e^{-\alpha_{i}-\alpha_{j}})h_{j}(e^{-\alpha_{i}})$, et en
utilisant la deuxi{\`e}me relation de la proposition~\ref{YangBaxter}, on obtient le
r{\'e}sultat. Les autres cas se traitent de la m{\^e}me mani{\`e}re.

\smallskip

Le terme $\sum_{l(w) \leq m \leq l(v)} \sum (e^{-\beta_{j_{1}}}-1) 
\cdots (e^{-\beta_{j_{m}}}-1)$ du th{\'e}or{\`e}me~\ref{ktheorieVD} est le coefficient
de $\mathcal{R}_{v}$ sur
$u_{w}$ dans la base $\{u_{w}\}_{w \in W}$ de $\mathcal{H}$ et est donc bien ind{\'e}pendant de la
d{\'e}composition r{\'e}duite de $v$ choisie.

\medskip

Notons $\tilde{\psi}^{w}$ l'{\'e}l{\'e}ment de $F(W,R[T])$ d{\'e}fini par la formule du
th{\'e}or{\`e}me~\ref{ktheorieVD} (on pose $\tilde{\psi}^{w}(v)=0$ si $v$ n'est pas plus 
grand que $w$). Pour d{\'e}montrer ce th{\'e}or{\`e}me, il suffit de montrer que 
les fonctions $(\tilde{\psi}^{w})_{w \in W}$ v{\'e}rifient les quatre propri{\'e}t{\'e}s
de la proposition~\ref{propositionKKK}. Les propri{\'e}t{\'e}s $(i)$ et 
$(iv)$ sont imm{\'e}diates.

\medskip

Pour d{\'e}montrer la propri{\'e}t{\'e} $(ii)$, rappelons tout d'abord le lemme
suivant (voir \cite{bou}, VI.$1.6$, corollaire $2$) : 

\begin{lemm}

Soit $v \in W$ et $v=s_{i_{1}} \cdots s_{i_{k}}$ une d{\'e}composition r{\'e}duite de $v$,
alors $ \Delta(v^{-1})=\{ \beta_{j}, 1 \leq j
  \leq k \}$.

\end{lemm}




D'apr{\`e}s ce lemme et la formule~\ref{rho1}, on a donc :

$$\tilde{\psi}^{w}(w)=\prod_{\beta \in \Delta(w^{-1})}e^{\beta}\prod_{\beta
\in \Delta(w^{-1})}(e^{-\beta}-1)=\prod_{\beta \in
\Delta(w^{-1})}(1-e^{\beta})=\psi^{w}(w),$$
ce qui nous donne la propri{\'e}t{\'e} $(ii)$.

\bigskip

Montrons maintenant que les $(\tilde{\psi}^{w})_{w \in W}$ 
v{\'e}rifient la propri{\'e}t{\'e} $(iii)$ de la
proposition~\ref{propositionKKK}.

Soit $w \in W$ et $s_{i}$ une r{\'e}flexion simple. 

\smallskip

Supposons tout d'abord
$ws_{i}>w$. Il faut alors montrer que pour tout $v \in W$, on a : 
$$\tilde{\psi}^{w}(v)=\tilde{\psi}^{w}(vs_{i})e^{-v\alpha_{i}}.$$

On peut supposer $vs_{i}>v$. 

Si $v$ n'est pas plus grand que $w$, $vs_{i}$ non
plus. En effet, $w$ n'a pas de d{\'e}composition r{\'e}duite qui finit par $s_{i}$ car
$ws_{i}>w$. Si $vs_{i}$ est plus grand que $w$, toute d{\'e}composition r{\'e}duite de
$v$ admet donc une sous-d{\'e}composition r{\'e}duite {\'e}gale {\`a} $w$, ce qui contredit le fait
que $v$ n'est pas plus grand que $w$.

On peut donc supposer $w \leq v <vs_{i}$. Comme $w$ n'a aucune
d{\'e}composition qui finit par $s_{i}$, la somme est la m{\^e}me {\`a} gauche et {\`a}
droite de l'{\'e}galit{\'e}. Il suffit donc de v{\'e}rifier $e^{\rho - v\rho}=
e^{\rho - vs_{i}\rho}e^{-v\alpha_{i}}$, ce qui est une cons{\'e}quence
imm{\'e}diate de la formule~\ref{rho1}.

\medskip

Supposons maintenant $ws_{i}<w$. Il faut montrer que pour tout $v \in 
W$, on a :

 \begin{equation} \label{dem} \frac{\tilde{\psi}^{w}(v)-\tilde{\psi}^{w}
(vs_{i})e^{-v\alpha_{i}}}
    {1-e^{-v\alpha_{i}}}=\tilde{\psi}^{w}(v)+\tilde{\psi}^{ws_{i}}(v). 
    \end{equation}

Supposons tout d'abord $vs_{i}>v$. On se place dans le cas o{\`u} $w \leq 
vs_{i}$ (sinon le r{\'e}sultat est trivial). On choisit une d{\'e}composition 
r{\'e}duite $v=s_{i_{1}} \cdots s_{i_{l}}$ de $v$. On prend pour 
$vs_{i}$ la d{\'e}composition $vs_{i}=s_{i_{1}} \cdots s_{i_{l}}s_{i}$. 
On trouve alors (en utilisant la formule~\ref{rho1}) : 

$$\tilde{\psi}^{w}(vs_{i})e^{-v\alpha_{i}}=\tilde{\psi}^{w}(v)+
(e^{-v\alpha_{i}}-1)(\tilde{\psi}^{w}(v)+\tilde{\psi}^{ws_{i}}(v)),$$
le premier terme venant des sous-d{\'e}compositions de $v$ ``{\'e}gales'' {\`a} 
$w$, le deuxi{\`e}me des m{\^e}mes sous-d{\'e}compositions de $v$ auxquelles 
on rajoute $s_{i}$ {\`a} la fin et qui redonnent donc $w$ (car 
$ws_{i}<w$ et donc $\underline{w}\, \underline{s}_i=\underline{w}$), 
et le troisi{\`e}me des sous-d{\'e}compositions de $v$ 
``{\'e}gales'' {\`a} $ws_{i}$. On trouve alors bien la formule~\ref{dem}.

Supposons maintenant $vs_{i}<v$. On peut appliquer ce qui pr{\'e}c{\`e}de 
{\`a} $v'=vs_{i}$ car $v's_{i}>v'$ et on trouve : 

$$\frac{\tilde{\psi}^{w}(vs_{i})-\tilde{\psi}^{w}(v)e^{v\alpha_{i}}}
     {1-e^{v\alpha_{i}}}=\tilde{\psi}^{w}(vs_{i})+
     \tilde{\psi}^{s_{i}w}(vs_{i}). $$

De plus, on peut appliquer le cas $ws_{i}>w$ {\`a} $w'=ws_{i}$ et on 
obtient : $\tilde{\psi}^{ws_{i}}(vs_{i}) = e^{v\alpha_{i}}
\tilde{\psi}^{ws_{i}}(v)$. En substituant ainsi 
$\tilde{\psi}^{ws_{i}}(vs_{i})$ dans l'expression pr{\'e}c{\'e}dente, on 
obtient la formule~\ref{dem}.

\chapter{Calcul de Schubert {\'e}quivariant} \label{schubertequivariant}

On reprend les notations de la section~\ref{cohomologievd}. 
Le calcul de Schubert {\'e}quivariant cherche {\`a} comprendre la structure
multiplicative de $H_T^*(X)$ en calculant les polyn{\^o}mes $p_{u,v}^w \in
S(\mathfrak{h}^*)$ qui v{\'e}rifient :
 
$$\hat{\xi}^u \hat{\xi}^v = \sum_{w \in W} p_{u,v}^w \hat{\xi}^w.$$

On voit facilement, par r{\'e}currence sur $l(w)$ et gr{\^a}ce aux relations 
$\xi^u(w)=0$ si $w \not\geq u$, que $ p_{u,v}^w=0$ sauf si 
$w \geq u$ et $w \geq v$.

Dans le cas o{\`u} $u=s_i$ est une r{\'e}flexion simple, la formule de Pieri-Chevalley (voir
\cite{kkc}) donne les valeurs de $p_{s_i,v}^w$ : 

$$\hat{\xi}^{s_i} \hat{\xi}^v = \xi^{s_i}(v)\hat{\xi}^v +
\sum_{v \rightarrow w} \rho_i(\beta^{\vee}(v,w)) \hat{\xi}^w.$$

Dans \cite{robinson}, Shawn Robinson g{\'e}n{\'e}ralise cette formule pour le type $A$ 
dans le cas o{\`u} $u=\prod_{i\leq k \leq j}s_k$ (pour $i\leq j$) est un produit de r{\'e}flexions
simples successives.

\bigskip

Dans \cite{kkc}, Kostant et Kumar donnent une formule g{\'e}n{\'e}rale pour ces
coefficients $p_{u,v}^w$. On note $Q(\mathfrak{h}^*)$ le corps des fractions de
$S(\mathfrak{h}^*)$ et $F(W;Q(\mathfrak{h}^*))$ la $Q(\mathfrak{h}^*)$-alg{\`e}bre
des fonctions de $W$ {\`a} valeurs dans $Q(\mathfrak{h}^*)$ munie de l'addition et
de la multiplication point par point. Pour tout $1 \leq i \leq r$, on d{\'e}finit un
op{\'e}rateur $A_i : F(W;Q(\mathfrak{h}^*)) \rightarrow F(W;Q(\mathfrak{h}^*))$ 
par :
$$\forall u \in W,  A_i(f)(u) = \frac{f(us_i)-f(u)}{u\alpha_i}.$$

Les op{\'e}rateurs $A_i$ v{\'e}rifiant les relations de tresse, on peut d{\'e}finir un
op{\'e}rateur $A_w$ pour tout $w \in W$. De plus, on d{\'e}finit l'op{\'e}rateur 
$s_i  : F(W;Q(\mathfrak{h}^*)) \rightarrow F(W;Q(\mathfrak{h}^*))$ 
par :
$$\forall u \in W,  s_i(f)(u) = f(us_i).$$

Les fonctions $\{ \xi^w \}_{w \in W} $ v{\'e}rifient alors les relations 
suivantes : 
 $$ \left\{ \begin{array}{ll} A_{i}\xi^{w} =  \xi^{ws_{i}}
 & {\rm si } \hspace{0,2 cm} ws_{i}<w, \\ 
A_{i}\xi^{w}=0
 & {\rm si } \hspace{0,2 cm} ws_{i}>w.  \end{array} \right.$$

\medskip

Soit $w \in W$ et soit $w=s_{i_1} \cdots s_{i_n}$ une d{\'e}composition r{\'e}duite de
$w$. Alors, pour tout couple $(u,v) \in W^2$, les polyn{\^o}mes 
$p_{u,v}^w$ sont donn{\'e}s par la formule suivante : 

\begin{equation} \label{*kk}
p_{u,v}^w=\!\!\!\!\!\!\!\!\!\!\!\!\!\!\! \sum_{
 \tiny \begin{array}{cc}  1 \leq j_1 < \cdots <j_m\leq  n \\
 {\rm tels \, que \, } s_{i_{j_1}} \cdots s_{i_{j_m}}=u
\end{array}}\! \!\!\!\!\!\!\!\!\!\!\!\!\!\!\! A_{i_1} \circ \cdots \circ \hat{\hat{A}}_{i_{j_1}}
\circ \cdots \circ \hat{\hat{A}}_{i_{j_m}} \circ \cdots \circ A_{i_n}(\xi^v)(1),
\end{equation}
o{\`u} $m=l(u)$, et o{\`u} la notation  $\hat{\hat{A}}_{i_{j_1}}$ signifie qu'on remplace
l'op{\'e}rateur $A_{i_{j_1}}$ par l'op{\'e}rateur $s_{i_{j_1}}$.

On va donner une formule un peu plus explicite pour calculer ces
coefficients. Cette formule g{\'e}n{\'e}ralise celle donn{\'e}e par Haibao Duan pour la
cohomologie ordinaire \cite{duan}. 
Pour trouver cette formule, il faut mieux comprendre la structure multiplicative de
la cohomologie {\'e}quivariante des vari{\'e}t{\'e}s de Bott-Samelson.

\section{G{\'e}n{\'e}ralit{\'e}s}

Soit $A$ un anneau commutatif unitaire, et soit $N \geq 1$ un entier naturel. On consid{\`e}re
une liste $D=\{d_{i,j} \}_{1 \leq i \leq j \leq N}$ d'{\'e}l{\'e}ments de $A$. Pour 
$1 \leq k \leq N$, on 
d{\'e}finit le polyn{\^o}me $Q_k \in A[X_1, \ldots, X_N]$ par : 
$$Q_k= X_k^2 -d_{k,k}X_k  -\sum_{l < k} d_{l,k}X_kX_l,
$$
et on d{\'e}finit alors la $A$-alg{\`e}bre $\mathcal{A}_{D}$ par : 
$$\mathcal{A}_{D}=A[X_1, \ldots, X_N] /<Q_1, \ldots, Q_N>,$$
o{\`u} $<Q_1, \ldots, Q_N>$ d{\'e}signe l'id{\'e}al de $A[X_1, \ldots, X_N]$ engendr{\'e} par 
$Q_1, \ldots, Q_N$. 
On note $x_i \in \mathcal{A}_D$ l'image de $X_i$ dans $\mathcal{A}_D$ et
pour $\epsilon \in \mathcal{E} =\{0,1 \}^N$, on pose $x^{\epsilon}=\prod_{i\in
  \pi_+(\epsilon)}x_i$.

\begin{prop}

La famille $\{x^{\epsilon} \}_{\epsilon \in \mathcal{E}}$ est une base du
$A$-module $\mathcal{A}_D$ qui est donc un $A$-module libre de rang $2^N$.

\end{prop}

\begin{proof} 

On proc{\`e}de par r{\'e}currence sur $N\geq 1$. Pour $N=1$, le r{\'e}sultat est imm{\'e}diat.

Supposons le r{\'e}sultat v{\'e}rifi{\'e} au rang $N-1$, et soit $D=\{d_{i,j} \}_{1 \leq i
  \leq j \leq N}$ une liste d'{\'e}l{\'e}ments de $A$. Alors, comme pour $1 \leq k \leq
  N-1$, $Q_k \in A[X_1, \ldots, X_{N-1}]$ : 

$$\mathcal{A}_D \simeq
 \big(A[X_1, \ldots, X_{N-1}]/<Q_1, \ldots, Q_{N-1}>\big)[X_N]/\overline{Q}_{N},$$
o{\`u} $\overline{Q}_{N}$ d{\'e}signe l'image de $Q_N$ dans $ \left(A[X_1, \ldots,
X_{N-1}]/\!\!<Q_1, \ldots, Q_{N-1}>\right)[X_N]$.

On conclut alors en utilisant le cas $N-1$ pour l'anneau $A$, puis le cas $N=1$
pour l'anneau $ A[X_1, \ldots, X_{N-1}]/<Q_1, \ldots, Q_{N-1}>$.

\end{proof}

On va expliciter la structure multiplicative de $\mathcal{A}_D$. 
On note $q_{\epsilon_1, \epsilon_2}^{\epsilon}$ les {\'e}l{\'e}ments de $A$ d{\'e}finis 
par : 
$$x^{\epsilon_1}x^{\epsilon_2}=\sum_{\epsilon \in \mathcal{E}} 
q_{\epsilon_1, \epsilon_2  }^{\epsilon} x^{\epsilon}.$$

Plus g{\'e}n{\'e}ralement, pour tout polyn{\^o}me $P \in A[x_1, \ldots, x_N]$, on note
$P^{\epsilon}$ les {\'e}l{\'e}ments de $A$ d{\'e}finis par : 
$$P=\sum_{\epsilon \in \mathcal{E}} P^{\epsilon} x^{\epsilon},$$
o{\`u} on continue {\`a} noter $P$ l'{\'e}l{\'e}ment de $\mathcal{A}_D$ d{\'e}fini par $P$.

\begin{defi}

Soit $\epsilon \in \mathcal{E}$. On note $\{i_1 < \cdots <i_l \}$ les {\'e}l{\'e}ments
 de $\pi_+(\epsilon)$. On d{\'e}finit alors l'application $T^{\epsilon} : 
 A[x_{1}, x_{2},\ldots,  x_{N}] \rightarrow A$ de la mani{\`e}re suivante :

\begin{enumerate}

\item[$(i)$] $T^{\epsilon}$ est $A$-lin{\'e}aire,

\item[$(ii)$] si $P$ est un mon{\^o}me qui n'est pas dans $ A[x_{i_1}, 
x_{i_2},\ldots,  x_{i_l}]$, alors 
$T^{\epsilon}( P)=0$,

\item[$(iii)$] si $P \in   A[x_{i_1}, x_{i_2},\ldots,  x_{i_{l-1}}]$, alors 
$T^{\epsilon}( P)=0$,

\item[$(iv)$] $T^{(i_{1})}( x_{i_{1}}^s)=d_{i_1,i_1}^{s-1}$,

\item[$(v)$] Si $Q \in A[x_{i_1}, x_{i_2},\ldots,  x_{i_{l-1}}]$, alors
pour $s \geq 1$ : 
$$ T^{\epsilon}(Qx_{i_l}^s)=T^{\epsilon-(i_l)}\Big[Q(d_{i_l,i_l}+
\sum_{j<l}d_{i_j,i_l}x_{i_j})^{s-1} \Big] .$$

\end{enumerate}

Ces cinq relations d{\'e}finissent compl{\`e}tement (r{\'e}cursivement) les applications
 $T^{\epsilon}$.

\end{defi}

\begin{exem}

 Prenons $N=2$, $d_{1,2}=1$ et $\epsilon=(\mathbf{1})=(1,1)$. Alors : 

\medskip

$T^{(\mathbf{1})}(x_{1}^s)=0$ pour tout $s$, \, \, \,
$T^{(\mathbf{1})}(x_{2})=0$,

\smallskip

si $t \geq 1$ et $s \geq 1$,
$T^{(\mathbf{1})}(x_1^sx_2^t)=T^{(1,0)}\big(x_1^s(d_{2,2}+x_1)^{t-1}\big)=
 d_{1,1}^{s-1}(d_{1,1}+d_{2,2})^{t-1} $,

\smallskip

et si $t \geq 2$,
$T^{(\mathbf{1})}(x_2^t)=T^{(1,0)}\big((d_{2,2}+x_1)^{t-1}\big)=\sum_{k=0}^{t-2}C_{t-1}^k
d_{2,2}^kd_{1,1}^{t-2-k}$.

\end{exem}

\begin{prop}  \label{multiplicationA}

Pour tout $\epsilon \in \mathcal{E}$ et tout polyn{\^o}me $P \in A[x_{1}, x_{2},
\ldots,  x_{N}]$ :
$$P_{\epsilon}=T^{\epsilon}(P).$$

En particulier, pour tout couple $(\epsilon_1, \epsilon_2) \in \mathcal{E}^2$ : 
$$q_{\epsilon_1, \epsilon_2 }^{
  \epsilon}=T^{\epsilon}(x^{\epsilon_1}x^{\epsilon_2}).$$

\end{prop}

\begin{proof}

Il faut v{\'e}rifier les cinq relations qui d{\'e}finissent les op{\'e}rateurs
$T^{\epsilon}$.

 La premi{\`e}re est imm{\'e}diate.

Pour tout $1 \leq k \leq N$, on a la formule suivante :

\begin{equation}  \label{carr{\'e}sA}
x_k^2=d_{k,k}x_k+\sum_{l<k}d_{l,k}x_lx_k.
\end{equation}

Les relations $(ii)$ et $(iii)$ se d{\'e}duisent de l'{\'e}quation~\ref{carr{\'e}sA}.

De plus, on d{\'e}montre par r{\'e}currence sur $s \geq 1$, gr{\^a}ce {\`a} la formule~\ref{carr{\'e}sA}, 
 que pour tout $1 \leq k \leq N$, et
tout $s \geq 1$ :
 
$$x_{k}^s=\Big(\sum_{i=0}^{s-1}C_{s-1}^i d_{k,k}^i
\big(\sum_{l<k}d_{l,k}x_l\big)^{s-1-i}\Big)x_k=
\big(d_{k,k}+\sum_{l<k}d_{l,k}x_l\big)^{s-1}x_k.$$

Cette formule nous permet alors de montrer les relations $(iv)$ et $(v)$.

\end{proof}

\section{Structure multiplicative de ${H}_T^*(\Gamma)$}

Soit $\mu_1, \ldots, \mu_N$ une suite de $N$ racines simples non
n{\'e}cessairement distinctes. On pose $\Gamma=\Gamma(\mu_1, \ldots, \mu_N)$. On
peut appliquer les r{\'e}sultats de la section pr{\'e}c{\'e}dente {\`a} ${H}_T^*(\Gamma)$ (plus
g{\'e}n{\'e}ralement {\`a} $H_D^*(Y)$, o{\`u} $Y$ est une tour de Bott).

En effet, si on pose pour $1 \leq l < k \leq N$, $d_{l,k}=-\mu_k(\mu_l^{\vee})$,
et $d_{k,k}=\mu_k$, alors, d'apr{\`e}s le th{\'e}or{\`e}me~\ref{carreBS}, 
$H_T^*(\Gamma)$ s'identifie {\`a} l'alg{\`e}bre
$\mathcal{A}_D$, o{\`u} on prend pour $A$ l'anneau $S(\mathfrak{h}^*)$.
De plus, pour tout $\epsilon \in \mathcal{E}$, 
$\hat{\sigma}_{\epsilon}^T$ s'identifie {\`a} $x^{\epsilon}$.

\begin{defi}

Soit $\epsilon \in \mathcal{E}$. On note $\{i_1 < \cdots <i_l \}$ les {\'e}l{\'e}ments
 de $\pi_+(\epsilon)$. On d{\'e}finit alors l'application $T^{\epsilon}_{\mu_1, \ldots, \mu_N} : 
S(\mathfrak{h}^*) [x_{1}, x_{2},\ldots,  x_{N}] \rightarrow S(\mathfrak{h}^*)$ 
de la mani{\`e}re suivante :

\begin{enumerate}

\item[$(i)$] $T^{\epsilon}_{\mu_1, \ldots, \mu_N}$ est $S(\mathfrak{h}^*)$-lin{\'e}aire,

\item[$(ii)$] si $P$ est un mon{\^o}me qui n'est pas dans $ S(\mathfrak{h}^*)[x_{i_1}, 
x_{i_2},\ldots,  x_{i_l}]$, alors 
$T^{\epsilon}_{\mu_1, \ldots, \mu_N}( P)=0$,

\item[$(iii)$] si $P \in   S(\mathfrak{h}^*)[x_{i_1}, x_{i_2},\ldots,  x_{i_{l-1}}]$, alors 
$T^{\epsilon}_{\mu_1, \ldots, \mu_N}( P)=0$,

\item[$(iv)$] $ T^{(i_{1})}_{\mu_1, \ldots, \mu_N} ( x_{i_{1}}^s)=\mu_{i_1}^{s-1}$,

\item[$(v)$]  Si $Q \in S(\mathfrak{h}^*)[x_{i_1}, x_{i_2},\ldots,  x_{i_{l-1}}]$, alors
pour $s \geq 1$ : 
$$ T^{\epsilon}_{\mu_1, \ldots,
  \mu_N}(Qx_{i_l}^s)=T^{\epsilon-(i_l)}\Big[Q(\mu_{i_l}
-\sum_{j<l}\mu_{i_l}(\mu_{i_j}^{\vee})x_{i_j})^{s-1} \Big] .$$

\end{enumerate}

De plus, on pose $T_{\mu_1, \ldots, \mu_N}=T^{(\mathbf{1})}_{\mu_1, \ldots, \mu_N}$.

\end{defi}

\smallskip

D'apr{\`e}s la proposition~\ref{multiplicationA}, la structure multiplicative de
$H_T^*(\Gamma)$ est donn{\'e}e par le th{\'e}or{\`e}me suivant : 

\begin{theo} \label{multiplicationBS}

Pour tout couple $(\epsilon_1, \epsilon_2) \in \mathcal{E}^2$ : 
$$\hat{\sigma}_{\epsilon_1}^T \hat{\sigma}_{\epsilon_2}^T = 
\sum_{\epsilon \in \mathcal{E}} T^{\epsilon}_{\mu_1, \ldots, \mu_N}
(x^{\epsilon_1}x^{\epsilon_2}) \, \hat{\sigma}_{\epsilon}^T.$$

\end{theo}

\section{Calcul de Schubert {\'e}quivariant}

Soit $w \in W$ et soit $w=s_{\mu_1} \cdots s_{\mu_N}$ une d{\'e}composition
quelconque de $w$. La
proposition~\ref{decompositionVBS} et le th{\'e}or{\`e}me~\ref{multiplicationBS} nous
donnent la formule suivante : 

\begin{theo} \label{formulefinale}
 
Pour tout couple $(u,v) \in W^2$ : 
\begin{equation} \label{formule}
p_{u,v}^w=T_{\mu_1, \ldots, \mu_N} \Big[
\big(\! \! \! \!  \! \sum_{ \tiny \begin{array}{cc}\epsilon \in 
\mathcal{E}, l(\epsilon)=l(u) \\ 
  { \rm et } \,  v(\epsilon)=u
\end{array}}\! \! \! \!  \!   \! x^{\epsilon} \big)
\big(\! \! \! \! \! \sum_{ \tiny \begin{array}{cc}\epsilon' 
\in \mathcal{E}, l(\epsilon')=l(v) \\ 
  { \rm et } \,  v(\epsilon')=v
 \end{array}} \! \! \! \! \!   \!   x^{\epsilon'} \big) \Big].
\end{equation}

\end{theo}

Dans cette formule, on somme a priori sur plus de termes que dans la
formule~\ref{*kk}, mais chaque terme est beaucoup plus facile {\`a} calculer. De
plus, cette formule exprime directement les polyn{\^o}mes $p_{u,v}^w$ en fonction  des
racines simples et des nombres de Cartan.

\begin{exem}

On prend le cas $A_5$, $u=s_{5}s_{2}$, $v=s_4s_5s_3s_4$ et
$w=s_4s_5s_2s_3s_4$. Alors : 
$$p_{u,v}^w=T_{\alpha_4,\alpha_5, \alpha_2,\alpha_3,\alpha_4}
\big[(x_2x_3)(x_1x_2x_4x_5)\big]$$
$$=T_{\alpha_4,\alpha_5, \alpha_2,\alpha_3,\alpha_4}
(x_1x_2^2x_3x_4x_5)
=T_{\alpha_4,\alpha_5,\alpha_2,\alpha_3,\alpha_4}^{(1,1,0,0,0)}
(x_1x_2^2)=\alpha_4+\alpha_5.$$

\end{exem}

\begin{exem}

On consid{\`e}re le cas $G_2$, et on prend $u=s_2s_1s_2$,
$v=s_1s_2s_1$, et $w=s_1s_2s_1s_2$. Alors : 
$$p_{u,v}^w=T_{\alpha_1,\alpha_2,\alpha_1,\alpha_2}
\big[(x_2x_3x_4)(x_1x_2x_3)\big]$$
$$=T_{\alpha_1,\alpha_2,\alpha_1,\alpha_2}
(x_1x_2^2x_3^2x_4)
=T_{\alpha_1,\alpha_2,\alpha_1,\alpha_2}^{(1,1,0,0)}\big[x_1x_2^2(\alpha_1+3x_2-2x_1)\big]
$$
$$=\alpha_1T_{\alpha_1,\alpha_2,\alpha_1,\alpha_2}^{(1,0,
0,0)}\big[x_1(\alpha_2+x_1)\big]+3T_{\alpha_1,\alpha_2,\alpha_1,\alpha_2}^{(1,0,
0,0)}\big[x_1(\alpha_2+x_1)^2\big]
-2T_{\alpha_1,\alpha_2,\alpha_1,\alpha_2}^{(1,0,0,0)}\big[x_1^2(\alpha_2+x_1)\big]
$$
$$=\alpha_1(\alpha_2+\alpha_1)+3(\alpha_2^2+2\alpha_1\alpha_2+\alpha_1^2)
-2(\alpha_1\alpha_2+\alpha_1^2)$$
$$=2\alpha_1^2 +5 \alpha_1\alpha_2 +3\alpha_2^2.
$$

\end{exem}

\begin{exem}

Pour calculer un produit dans le cas fini, au lieu de calculer chaque
coefficient $p_{u,v}^w$ avec la formule~\ref{formule}, 
on peut aussi utiliser le plongement
$g^* : H_T^*(X) \rightarrow H_T^*(\Gamma)$, o{\`u} $\Gamma$ est la vari{\'e}t{\'e} de Bott-Samelson
associ{\'e}e {\`a} une suite de racines simples correspondant {\`a} une d{\'e}composition r{\'e}duite
de $w_0$, le plus grand {\'e}l{\'e}ment du groupe de Weyl.

On se place par exemple dans le cas $A_3$, et on prend $w_0 = s_3s_2s_1s_3s_2s_3$
pour d{\'e}composition r{\'e}duite de $w_0$. 

Alors $H_T^*(\Gamma)$ est l'alg{\`e}bre de polyn{\^o}mes
 $S(\mathfrak{h}^*)[x_1,x_2,x_3,x_4,x_5,x_6]$
 quotient{\'e}e par les relations : 
\begin{eqnarray*}
x_1^2 & = & \alpha_3x_1 \\
  x_2^2 & = & \alpha_2x_2 + x_1x_2 \\
  x_3^2 & = & \alpha_1x_3 + x_2x_3+x_1x_3 \\
x_4^2 & = & \alpha_3x_4 +x_2x_4-2x_1x_4 \\
x_5^2 & = & \alpha_2x_5 + x_4x_5+x_3x_5-2x_2x_5+x_1x_5 \\
x_6^2 & = & \alpha_3x_6 + x_5x_6-2x_4x_6+x_2x_6-2x_1x_6
\end{eqnarray*}

De plus, on a en particulier : 
\begin{eqnarray*}
g^*(\hat{\xi}^{s_3s_2s_1}) & = &x_1x_2x_3 \\
 g^*(\hat{\xi}^{s_3s_2}) & = & x_1x_2+x_1x_5+x_4x_5 \\
g^*(\hat{\xi}^{s_3s_2s_1s_2}) & = &  x_1x_2x_3x_5\\
g^*(\hat{\xi}^{s_3s_2s_1s_3s_2}) &= & x_1x_2x_3x_4x_5
\end{eqnarray*}

\medskip

Si on veut calculer $\hat{\xi}^{s_3s_2s_1}\hat{\xi}^{s_3s_2}$, on utilise $g^*$
: 
$$g^*(\hat{\xi}^{s_3s_2s_1}\hat{\xi}^{s_3s_2})=x_1x_2x_3(x_1x_2+x_1x_5+x_4x_5)$$
$$= \alpha_3x_1( \alpha_2x_2 + x_1x_2)x_3+\alpha_3x_1x_2x_3x_5+x_1x_2x_3x_4x_5$$
$$=(\alpha_3^2+\alpha_2\alpha_3)x_1x_2x_3+\alpha_3x_1x_2x_3x_5+x_1x_2x_3x_4x_5,$$
et donc : 
$$\hat{\xi}^{s_3s_2s_1}\hat{\xi}^{s_3s_2}=(\alpha_3^2+\alpha_2\alpha_3)
\hat{\xi}^{s_3s_2s_1}+\alpha_3\hat{\xi}^{s_3s_2s_1s_2}+\hat{\xi}^{s_3s_2s_1s_3s_2}.$$

\end{exem}

\smallskip

\begin{rema}

En {\'e}valuant ces polyn{\^o}mes $p_{u,v}^w$ en $0$, on retrouve la formule donn{\'e}e par
Haibao Duan dans \cite{duan} pour la cohomologie ordinaire. Cette formule est
utilis{\'e}e par Duan et Zhao dans \cite{duan2}
pour proposer un algorithme de calcul de Schubert ordinaire.

\end{rema}

Le nombre de termes {\`a} calculer dans la formule~\ref{formule} d{\'e}pend du choix de la
d{\'e}composition de $w$, et il arrive que certains termes s'annulent. En effet, un
th{\'e}or{\`e}me de William Graham \cite{gra2} affirme que les polyn{\^o}mes $p_{u,v}^w$
sont ``positifs'' dans le sens o{\`u} ils sont combinaisons lin{\'e}aires 
{\`a} coefficients positifs de termes de la forme
$\alpha_I= \prod_{1 \leq i \leq r}\alpha_{i}^{n_i}$ pour $I \in
\mathbb{N}^r$. Or, la formule~\ref{formule} peut comporter des termes
``n{\'e}gatifs''.

 Donnons un exemple tr{\`e}s simple pour illustrer ce probl{\`e}me : on se place dans le
 cas $A_2$, et on veut calculer $p_{s_2, s_1s_2}^{w}$ o{\`u} $w=s_1s_2s_1=s_2s_1s_2$
 est le plus grand {\'e}l{\'e}ment de $W$. Si on prend la premi{\`e}re d{\'e}composition de $w$,
 on voit imm{\'e}diatement que $p_{s_2, s_1s_2}^{w}=0$. En revanche, si on prend la
 deuxi{\`e}me d{\'e}composition, on doit faire le calcul suivant : 

$$p_{s_2,s_1s_2}^{s_2s_1s_2}=T_{\alpha_2,\alpha_1,\alpha_2}
\big[(x_1+x_3)x_2x_3\big]=1+T_{\alpha_2,\alpha_1,\alpha_2}(x_2x_3^2)$$
$$=1+ T_{\alpha_2,\alpha_1,\alpha_2}^{(1,1,0)}\big[x_2(\alpha_2+x_2-2x_1)\big]
=1+1-2=0. $$ 

\bigskip

\begin{exem}

Donnons un dernier exemple de calcul. On consid{\`e}re le cas $A_6$, et on prend
$u=s_1s_3s_5s_6$, $v=s_2s_5s_6$ et $w=s_1s_2s_3s_4s_5s_6$. Alors : 
$$p_{u,v}^w=T_{\alpha_1,\alpha_2,\alpha_3,\alpha_4,\alpha_5,\alpha_6}\big[
(x_1x_3x_5x_6)(x_2x_5x_6)\big]
=T_{\alpha_1,\alpha_2,\alpha_3,\alpha_4,\alpha_5,\alpha_6}
(x_1x_2x_3x_5^2x_6^2)$$
$$=T_{\alpha_1,\alpha_2,\alpha_3,\alpha_4,\alpha_5,\alpha_6}^{(1,1,1,1,1,0)}
\big[ x_1x_2x_3x_5^2 (\alpha_6+x_5)\big]
=\alpha_6+T_{\alpha_1,\alpha_2,\alpha_3,\alpha_4,\alpha_5,\alpha_6}^{(1,1,1,1,0,0)}
\big[ x_1x_2x_3(\alpha_5+x_4)^2\big]$$
$$=\alpha_6+2\alpha_5+T_{\alpha_1,\alpha_2,\alpha_3,\alpha_4,\alpha_5,\alpha_6}^{(1,1,1
 ,1,
0,0) }( x_1x_2x_3x_4^2)
=\alpha_1+\alpha_2+\alpha_3+\alpha_4+2\alpha_5+\alpha_6.$$


\end{exem}

 \bibliography{Ktheory}

\providecommand{\bysame}{\leavevmode ---\ }
\providecommand{\og}{``}
\providecommand{\fg}{''}
\providecommand{\smfandname}{et}
\providecommand{\smfedsname}{\'eds.}
\providecommand{\smfedname}{\'ed.}
\providecommand{\smfmastersthesisname}{M\'emoire}
\providecommand{\smfphdthesisname}{Th\`ese}
\begin{thebibliography}{10}

\bibitem{aa}
{\scshape A.~{\sc Arabia}} -- {\og Cohomologie {T}-\'equivariante de la
  vari\'et\'e de drapeaux d'un groupe de {K}ac-{M}oody\fg}, \emph{Bulletin de
  la Soci\'et\'e Math\'ematique de France} \textbf{117} (1989), p.~129--165.

\bibitem{atiyah}
{\scshape M.~{\sc Atiyah}} -- \emph{{K}-theory}, Benjamin, 1967.

\bibitem{livreaudin}
{\scshape M.~{\sc Audin}} -- \emph{The topology of torus actions on symplectic
  manifolds}, Progress in Mathematics, vol.~93, Birkh{\"a}user, 1991.

\bibitem{bv}
{\scshape N.~{\sc Berline} {\normalfont \smfandname} M.~{\sc Vergne}} -- {\og
  Classes caract\'eristiques \'equivariantes. {F}ormules de localisation en
  cohomologie \'equivariante\fg}, \emph{C.R.A.S.} \textbf{295} (1982),
  p.~539--541.

\bibitem{sb}
{\scshape S.~{\sc Billey}} -- {\og Kostant polynomials and the cohomology of
  {G/B}\fg}, \emph{Duke Mathematical Journal} \textbf{96} (1999), p.~205--224.

\bibitem{bs}
{\scshape R.~{\sc Bott} {\normalfont \smfandname} H.~{\sc Samelson}} -- {\og
  Applications of the theory of {M}orse to symmetric spaces\fg}, \emph{American
  Journal of Mathematics} \textbf{70} (1958), p.~964--1028.

\bibitem{bou}
{\scshape N.~{\sc Bourbaki}} -- \emph{Groupes et alg\`ebres de {L}ie, chap.
  4-6}, Hermann, Paris, 1968.

\bibitem{hbs}
{\scshape H.~{\sc C. Hansen}} -- {\og On cycles in flag manifolds\fg},
  \emph{Mathematica Scandinavica} \textbf{33} (1973), p.~269--274.

\bibitem{ginzburg}
{\scshape N.~{\sc Chriss} {\normalfont \smfandname} V.~{\sc Ginzburg}} --
  \emph{Representation {T}heory and {C}omplex {G}eometry}, Birkh{\"a}user,
  1997.

\bibitem{cox}
{\scshape D.~A. {\sc Cox}} -- {\og The homogeneous coordinate ring of a toric
  variety\fg}, \emph{Journal of Algebraic Geometry} \textbf{4} (1995), no.~1,
  p.~17--50.

\bibitem{demazure}
{\scshape M.~{\sc Demazure}} -- {\og D{\'e}singularisation des vari{\'e}t{\'e}s
  de {S}chubert g{\'e}n{\'e}ralis{\'e}es\fg}, \emph{Annales scientifiques de
  l'Ecole normale sup{\'e}rieure} \textbf{7} (1974), p.~53--88.

\bibitem{duan}
{\scshape H.~{\sc Duan}} -- {\og Multiplicative rule of {S}chubert classes\fg},
  \emph{Math.{AG}/0306227} (2003).

\bibitem{duan2}
{\scshape H.~{\sc Duan} {\normalfont \smfandname} X.~{\sc Zhao}} -- {\og A
  program for multiplying {S}chubert classes\fg}, \emph{Math.{AG}/0309158}
  (2003).

\bibitem{fk}
{\scshape S.~{\sc Fomin} {\normalfont \smfandname} A.~{\sc N. Kirillov}} --
  {\og Universal exponential solution of the {Y}ang-{B}axter equation\fg},
  \emph{Letters in mathematical Physics} \textbf{37} (1996), p.~273--284.

\bibitem{gauss}
{\scshape S.~{\sc Gaussent}} -- {\og The fiber of the {B}ott-{S}amelson
  resolution\fg}, \emph{Indagationes Mathematicae} \textbf{12} (2001),
  p.~453--468.

\bibitem{km}
{\scshape V.~{\sc G. Kac}} -- \emph{Infinite dimensional {L}ie algebras},
  Cambridge University Press, 1985.

\bibitem{kp}
{\scshape V.~{\sc G. Kac} {\normalfont \smfandname} D.~{\sc H. Peterson}} --
  {\og Regular functions on certain infinite dimensional groups\fg},
  \emph{Arithmetic and Geometry-II}, Birkh{\"a}user, 1983, p.~141--166.

\bibitem{gra2}
{\scshape W.~{\sc Graham}} -- {\og Positivity in equivariant {S}chubert
  calculus\fg}, \emph{Duke Mathematical Journal} \textbf{109} (2001),
  p.~599--614.

\bibitem{gra}
\bysame , {\og Equivariant {K}-theory and {S}chubert varieties\fg},
  \emph{preprint} (2002).

\bibitem{tours}
{\scshape M.~{\sc Grossberg} {\normalfont \smfandname} Y.~{\sc Karshon}} --
  {\og Bott towers, complete integrability, and the extended character of
  representations\fg}, \emph{Duke Mathematical Journal} \textbf{76} (1994),
  p.~23--58.

\bibitem{hum}
{\scshape J.~E. {\sc Humphreys}} -- \emph{Reflection {G}roups and {C}oxeter
  {G}roups}, Cambridge studies in advanced mathematics, vol.~29, Cambridge
  University Press, 1990.

\bibitem{kkc}
{\scshape B.~{\sc Kostant} {\normalfont \smfandname} S.~{\sc Kumar}} -- {\og
  {T}he {N}il {H}ecke ring and cohomology of {G/P} for a {K}ac-{M}oody group
  {G}\fg}, \emph{Advances in {M}athematics} \textbf{68} (1986), p.~187--237.

\bibitem{kkk}
\bysame , {\og {T}-equivariant {K}-theory of generalized flag varieties\fg},
  \emph{Journal of Differential Geometry} \textbf{32} (1990), p.~549--603.

\bibitem{ku}
{\scshape S.~{\sc Kumar}} -- {\og The nil-{H}ecke ring and singularities of
  {S}chubert varieties\fg}, \emph{Inventiones Mathematicae} \textbf{123}
  (1996), p.~471--506.

\bibitem{livrekumar}
\bysame , \emph{Kac {M}oody {G}roups, their {F}lag {V}arieties and
  {R}epresentation {T}heory}, Progress in Mathematics, vol. 204,
  Birkh{\"a}user, 2002.

\bibitem{llt}
{\scshape A.~{\sc Lascoux}, B.~{\sc Leclerc} {\normalfont \smfandname} J.-Y.
  {\sc Thibon}} -- {\og Flag varieties and the {Y}ang-{B}axter equation\fg},
  \emph{Letters in mathematical Physics} \textbf{40} (1997), p.~75--90.

\bibitem{robinson}
{\scshape S.~{\sc Robinson}} -- {\og A {P}ieri-type formula for
  ${H}_{{T}}^{*}({SL}_{n}(\mathbb{C})/{B}$\fg}, \emph{Journal of Algebra}
  \textbf{249} (2002), no.~1, p.~38--58.

\bibitem{kequivariante}
{\scshape G.~{\sc Segal}} -- {\og {T}-equivariant {K}-theory\fg},
  \emph{Publications Math{\'e}matiques de l'Institut des Hautes Etudes
  Scientifiques} \textbf{34} (1968), p.~129--151.

\bibitem{mw3}
{\scshape M.~{\sc Willems}} -- {\og Cohomologie et {K}-th\'eorie
  \'equivariantes des vari\'et\'es de {B}ott-{S}amelson et des vari\'et\'es de
  drapeaux\fg}, \emph{Bulletin de la Société Mathématique de France} (à
  paraitre).

\bibitem{mw}
\bysame , {\og Cohomologie \'equivariante des vari\'et\'es de
  {B}ott-{S}amelson\fg}, \emph{Math.GR/0201050} (2002).

\bibitem{mw2}
\bysame , {\og K-th\'eorie \'equivariante des vari\'et\'es de drapeaux et des
  vari\'et\'es de {B}ott-{S}amelson\fg}, \emph{Math.AG/0204265} (2002).

\end{thebibliography}
 \bibliographystyle{smfplain}

\end{document}